# Feigenbaum-Coullet-Tresser universality and Milnor's Hairiness Conjecture

By Mikhail Lyubich

*To Yakov Grigorievich Sinai on his* 60<sup>th</sup> *birthday*

**Abstract**

We prove the Feigenbaum-Coullet-Tresser conjecture on the hyperbolicity of the renormalization transformation of bounded type. This gives the first computer-free proof of the original Feigenbaum observation of the universal parameter scaling laws. We use the Hyperbolicity Theorem to prove Milnor's conjectures on self-similarity and "hairiness" of the Mandelbrot set near the corresponding parameter values. We also conclude that the set of real infinitely renormalizable quadratics of type bounded by some $N > 1$ has Hausdorff dimension strictly between 0 and 1. In the course of getting these results we supply the space of quadratic-like germs with a complex analytic structure and demonstrate that the hybrid classes form a complex codimension-one foliation of the connectedness locus.

## Contents

# 1. Introduction

1.1. *The universality phenomenon.* In the 1970's Feigenbaum and independently Coullet and Tresser discovered a "universal scaling law" of transition from regular to chaotic dynamics through cascades of doubling bifurcations (see Figure 1). The meaning of this discovery is that the geometry of the bifurcation loci in generic one-parameter families of certain dynamical systems is independent of the specific family. The importance of this discovery for dynamical systems theory and physical applications (fluid dynamics, statistical physics etc.) was realised shortly.

To explain the universality phenomenon, the authors introduced a renormalization transformation $R$ in an appropriate space of dynamical systems, and conjectured that this transformation had a unique fixed point $f_*$, and that this point was hyperbolic, with one-dimensional unstable manifold [F1], [F2], [CT], [TC]. Originally stated only for the period-doubling case, this conjecture was later extended to a wider class of combinatorics of "bounded type," real as well as complex ones; see [DGP], [GSK]. In this paper we will prove this conjecture for the renormalization operator of real bounded type acting in the space of quadratic-like germs.

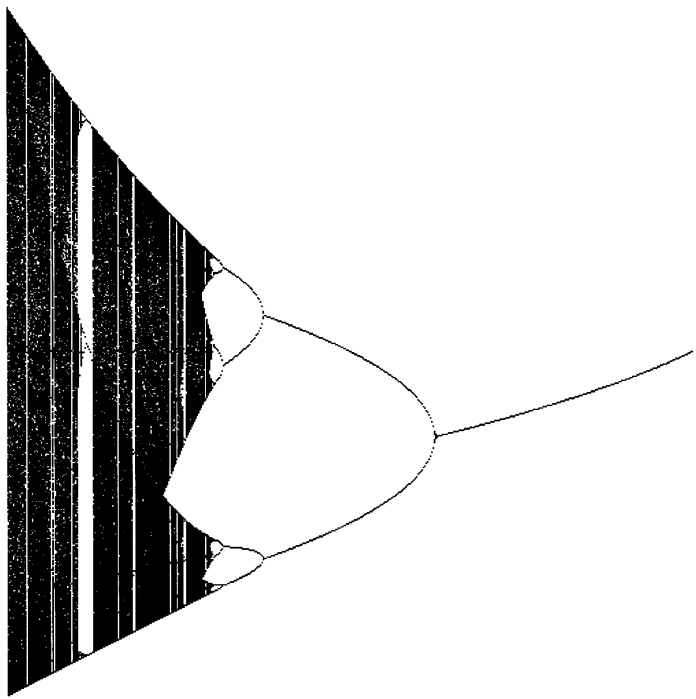

Figure 1. Cascade of doubling bifurcations.
This picture became symbolic for one-dimensional dynamics.

1.2. *Statement of the results.* The concepts used below (quadratic-like maps, hybrid classes, renormalization, Mandelbrot copies, combinatorial type) are basic in holomorphic dynamics and will be precisely defined in Section 3. In the next section, Section 4, we will supply the space $\mathcal{QG}$ of quadratic-like germs (considered up to rescaling) with topology and complex analytic structure.

Let $M_0$ stand for the Mandelbrot set, and $\mathcal{N}$ stand for the full family of Mandelbrot copies $M \subset M_0$ different from $M_0$ itself. To each $M \in \mathcal{N}$ corresponds the renormalization operator $R_M : \mathcal{T}_M \to \mathcal{QG}$ defined on the "renormalization strip" $\mathcal{T}_M \subset \mathcal{QG}$. This operator admits an analytic continuation to a neighborhood of $\mathcal{T}_M$.

Given a family $\mathcal{L}$ of disjoint Mandelbrot copies, we can consider the corresponding piecewise defined renormalization operator

$$(1.1) \qquad R_\mathcal{L} : \bigcup_{M \in \mathcal{L}} \mathcal{T}_M \to \mathcal{QG}.$$

If $\mathcal{L}$ is a finite family then $R_\mathcal{L}$ is called a renormalization operator of *bounded type*.

By a "real" quadratic-like map we mean a quadratic-like map preserving the real line. By a "real" Mandelbrot copy $M \in \mathcal{N}$ we mean a Mandelbrot copy centered on the real line. If the family $\mathcal{L}$ consists of real Mandelbrot copies then one says that the operator $R_\mathcal{L}$ has *real combinatorics*.

Let $\Sigma_d$ stand for the space of bi-infinite sequences in $d$ symbols, and let $\omega : \Sigma_d \to \Sigma_d$ be the shift transformation on this space.

HYPERBOLICITY THEOREM. *If there exists a renormalization operator $R = R_\mathcal{L}$ of real bounded type defined on the union of d renormalization strips, then there is a compact R-invariant set $\mathcal{A}$ (the "renormalization horseshoe") with the following properties*:

- *The restriction $R|\mathcal{A}$ is topologically conjugate to $\omega|\Sigma_d$ and is uniformly hyperbolic*;

- *Any stable leaf $\mathcal{W}^s(f)$, $f \in \mathcal{A}$, coincides with the hybrid class of $f$ and has codimension* 1;

- *Any unstable leaf $\mathcal{W}^u(f)$ is an analytic curve which transversally passes through all real hybrid classes except the cusp one (corresponding to $c = 1/4$).*

*Remark.* "Stable/unstable" leaves above mean the connected components containing $f$ of the sets of point whose forward/backward orbits are exponentially asymptotic to the corresponding orbit of $f$.

By a Feigenbaum quadratic $P_c : z \mapsto z^2 + c$ (or a Feigenbaum parameter value) we will mean an infinitely renormalizable map (or the corresponding parameter value) of bounded type. The following result was conjectured by Milnor [M]:

HAIRINESS THEOREM. *Let $c \in [-2, 1/4]$ be a real Feigenbaum parameter value. Then the rescalings of the Mandelbrot set near $c$ converge in the Hausdorff metric on compact sets to the whole complex plane.*

Everyone who saw computer pictures of the Mandelbrot set realized that it was not self-similar: Otherwise wandering around it would not be so fascinating. However some self-similar features are still observable. In particular, it was conjectured by Milnor that the little Mandelbrot sets around the Feigenbaum point of stationary type have asymptotically the same shape ([M, Conjs. 3.1 and 3.3]). The following result proves some of these conjectures. Here we state it in the case of stationary combinatorics, postponing the statement for bounded combinatorics until Section 9.

SELF-SIMILARITY THEOREM. *Let $M$ be a real Mandelbrot copy and $\sigma : M \to M_0$ be the homeomorphism of $M$ onto the whole Mandelbrot set $M_0$. Then $\sigma$ has a unique real fixed point $c$. Moreover, $\sigma$ is $C^{1+\alpha}$-conformal at $c$, with the derivative at $c$ equal to the Feigenbaum universal scaling constant $\lambda = \lambda_M > 1$.*

*Remark.* The Feigenbaum universal constant can actually be defined as the above derivative. However, the logic of our discussion makes it more natural to introduce it first as the unstable eigenvalue of the renormalization operator $R_M$ at its fixed point.

Any real Feigenbaum parameter value $c$ of stationary type is a limit of superattracting points $c_n$ of periods $p^n$ (where $p = p(M)$) obtained from the center of $M$ by $n$-fold "tuning" (see Section 5.1 for the definition). The following theorem gives the first computer independent proof of the Feigenbaum parameter scaling law in the quadratic family and the universal nature of this law.

UNIVERSALITY THEOREM. *Let $\mathcal{S} = \{f_\mu\}$ be a real analytic one-parameter family of quadratic-like maps transversally intersecting the hybrid class $\mathcal{H}_c$ at $\mu_*$. Then for all sufficiently big $n$, $\mathcal{S}$ has a unique intersection point $\mu_n$ near $\mu_*$ with the hybrid class $\mathcal{H}_{c_n}$, and*

$$|\mu_n - \mu_*| \sim a\lambda^{-n},$$

*where $\lambda = \lambda_M$ is "universal," i.e., independent of the particular family in question. In particular, $|c_n - c| \sim b\lambda^{-n}$.*

Given a finite family $\mathcal{L}$ of real Mandelbrot copies ("real family"), let $I_\mathcal{L} \subset [-2, 1/4]$ stand for the set of infinitely renormalizable parameter values of bounded type specified by $\mathcal{L}$, i.e.,

$$R^n P_c \in \bigcup_{M \in \mathcal{L}} \mathcal{T}_M, \; n = 0, 1, \dots .$$

One more application of the Hyperbolicity Theorem is the following:

HD THEOREM. *For any finite family $\mathcal{L}$ containing at least two elements, $I_\mathcal{L}$ is a Cantor set with Hausdorff dimension strictly in between $0$ and $1$.*

Along the lines of our work we also establish the following result:

QC THEOREM. *Any primitive Mandelbrot copy $M$ is quasi-conformally equivalent to the whole Mandelbrot set $M_0$.*

(Recall that a set $M$ is said to be nonprimitive, or satellite, if it is attached to some hyperbolic component of the Mandelbrot set.)

1.3. *Ingredients of the proofs.* The main idea of the proof of the Hyperbolicity Theorem is that for *complex analytic* transformations, lack of hyperbolicity on an invariant set $\mathcal{A}$ (satisfying certain assumptions) can be detected topologically. Namely, one can construct a point $f \in \mathcal{A}$ whose orbit is slowly shadowed by another orbit. For the renormalization operator, such a situation is ruled out by the Combinatorial Rigidity Theorem [L2].

Let us formulate here the above mentioned shadowing theorem in the simplest fixed-point case:

SMALL ORBITS THEOREM. *Let $\mathcal{B} \supset \mathcal{B}'$ be two complex Banach spaces such that the ball of $\mathcal{B}'$ is pre-compact in $\mathcal{B}$. Let $i : \mathcal{B}' \to \mathcal{B}$ stand for the natural embedding. Let $T : (\mathcal{U}, 0) \to (\mathcal{B}', 0)$ be a complex analytic map in a neighborhood $\mathcal{U} \subset \mathcal{B}$ of $0$, $R = i \circ T : (\mathcal{U}, 0) \to (\mathcal{B}, 0)$. Assume that the spectrum of $DR(0) : \mathcal{B} \to \mathcal{B}$ belongs to the closed unit disk and is not empty on the unit circle. Then $R$ has "slow small orbits"; that is, for any neighborhood $\mathcal{V} \ni 0$, there is an orbit $\{R^m f\}_{m=0}^\infty \subset \mathcal{V}$, such that*

$$\lim \frac{1}{m} \log \|R^m f\| = 0.$$

The idea of the proof of the Hairiness Conjecture is to pass to the unstable manifold of the renormalization fixed point where the Mandelbrot set becomes scaling invariant. Then we show that lack of hairiness would imply existence of a nontrivial automorphism of a tower with *a priori* bounds. But this situation is ruled out by the hairiness of the Feigenbaum Julia sets (McMullen [McM2]).

A substantial part of our work is to supply the space $\mathcal{QG}$ of quadratic-like germs and the space $\mathcal{E}$ of expanding circle maps with complex analytic structure and to demonstrate that the hybrid classes form a complex codimension-one foliation of the connectedness locus. Then a generalized version of the $\lambda$-lemma yields that this foliation is transversally quasi-conformal, which we exploit many times. In particular, this yields the QC Theorem.

At a Feigenbaum point $c$ we can do better, and show that the foliation is transversally $C^{1+\alpha}$-conformal along the hybrid class $\mathcal{H}_c$ (this is an expected regularity of a codimension-1 stable foliation). This yields the Self-Similarity and Universality Theorems. The HD Theorem follows from the hyperbolicity of the renormalization operator by a standard distortion argument.

1.4. *Structure of the paper.* This paper is organized as follows. In Section 2 we prove the Small Orbits Theorem.

In Section 3 we give a revisited account of the Douady-Hubbard theory of quadratic-like maps [DH2]. The main novelty of our approach is that the relation between quadratic-like and circle maps is given up to affine rather than conformal equivalence. This allows us, in particular, to extend the uniformization of $\mathbb{C} \setminus M_0$ to the "vertically holomorphic" uniformization of the complement $\mathcal{QG} \setminus \mathcal{C}$ of the connectedness locus.

In Section 4 we supply the space $\mathcal{QG}$ of quadratic-like germs (up to affine equivalence) and the space of expanding circle maps (up to rotation) with complex analytic structure (modeled on families of Banach spaces) and demonstrate that the Douady-Hubbard hybrid classes form a foliation of the connectedness locus $\mathcal{C}$ with complex codimension-one analytic leaves. Moreover, we show that $\mathcal{C}$ is the topological product of the space $\mathcal{E}$ of expanding circle maps by the Mandelbrot set $M_0$. We derive from this picture certain transversality results, and prove the QC Theorem.

In Section 5 we define the complex renormalization operator, analytically extend this operator beyond the renormalization strips, and show that it is transversally nonsingular. Then we state three crucial analytic results: *a priori bounds* ([MvS], [S2]), the Tower Rigidity Theorem [McM2], and the Combinatorial Rigidity Theorem [L2].

In Section 6 we prove the Hyperbolicity Theorem for stationary combinatorics. On the way to this result we give a new proof of the exponential convergence to the renormalization fixed point in its hybrid class based on the Schwarz Lemma in Banach spaces.

In Section 7 we prove in the stationary case the Hairiness, Self-Similarity and Universality Theorems, and discuss a relation of these results to the MLC Conjecture.

In the next two sections, 8 and 9, we extend the previous results from stationary to bounded combinatorics and prove the HD Theorem.

In Appendix 1 we collect for the reader's convenience some basic results and references on quasi-conformal maps.

In Appendix 2 we develop a theory of complex structures modeled on families of Banach spaces.

1.5. *Complex combinatorics and higher degrees.* Complex methods play so crucial a role for our proofs that one can wonder why we need the real line at all. Besides physical motivation, there is only one reason for that: a key technical result (complex *a priori* bounds) needed for the construction of the renormalization horseshoe $\mathcal{A}$ is not yet established for complex maps. Conjecturally, *a priori* bounds exist for all infinitely renormalizable maps with bounded combinatorics but they are established only for real maps ([MvS], [S2]) and for complex maps with sufficiently high combinatorics [L2]. Once the complex *a priori* bounds are established, our results become valid in the purely complex setting. In what follows we will state the results in this setting assuming *a priori* bounds.

Another natural extension of the Renormalization Conjecture is concerned with higher degrees of the maps under consideration at the critical point ("criticality"). All the results of this paper are extended in a straightforward way from the quadratic-like maps to polynomial-like maps with a single critical point of any even degree $d$. The only noteworthy point is that the Combinatorial Rigidity Theorem is still valid for infinitely renormalizable maps of this class with bounded combinatorics (see [L2, Remark at the end of §10.1]).

1.6. *Historical notes.* Feigenbaum made his first announcement of the Universality Phenomenon in 1976. The importance of this discovery was realized soon, and there has been a good effort to prove the conjectural renormalization picture. Prior to our work, the conjecture was proven (with the help of computers) in the period-doubling case with quadratic critical point; let us first summarize the development in this case.

The computer-assisted proof in the doubling case was given by Lanford [La1], with one missing ingredient (a transversality issue) filled by Eckmann and Wittwer [EW]. The unstable manifold at the corresponding fixed point was constructed numerically by Vul, Sinai and Khanin [VSK].

Later on, many ingredients of the picture were proven without computers. Existence of a renormalization fixed point (a solution of the "Cvitanović-Feigenbaum functional equation") was proven by Epstein [E1], [E2]. Existence of an unstable eigenvalue was proven by Epstein and Eckmann [EE]. The stable manifold was constructed by Sullivan and McMullen (see below). The ingredients which still required computers after all that (in the quadratic period doubling case) were the codimension and transversality issues. Thus, even in the period-doubling case our paper provides the first complete computer-free proof of the Renormalization Conjecture.

The importance of complex analytic machinery was realised early, particularly by Lanford and Epstein, who searched for renormalization fixed points in appropriate analytic functional spaces. In the mid 80's these ideas were greatly emphasized and expanded: ideas of holomorphic dynamics and geometric structures coming from the complex plane (and even 3D hyperbolic space) became the main tools in the field.

The renormalization operator was complexified by Douady and Hubbard [DH2]. A program of construction of the renormalization fixed point and its stable manifold by means of the Teichmüller theory was formulated by Sullivan in his address to ICM-86 in Berkeley [S1]. This program was carried out a few years later (see [MvS], [S2]). A different approach to the problem exploiting the idea of geometric limits was given by McMullen [McM2].

For a renormalization operator of bounded type, Sullivan and McMullen's theory provides us with the renormalization horseshoe $\mathcal{A}$ and proves uniform exponential contraction in the hybrid classes of $f \in \mathcal{A}$. However, existence of the unstable eigenvalue was not established even in the stationary period tripling case. Nor was it shown that the hybrid classes $\mathcal{H}(f)$, $f \in \mathcal{A}$, form a foliation with codimension-one complex analytic leaves. (In [McM2], [S2] the hybrid classes are treated intrinsically without embedding them into an ambient complex space.)

The above development based on complex methods treats the case of the quadratic critical point, or more generally, the "analytic" fixed point (i.e., having even criticality). However, the computer experiments suggest that the universality phenomenon is valid for any real criticality $\delta > 1$ as well. (By definition, a smooth unimodal map has criticality $\delta$ if it admits a representation $|\phi(x)|^\delta$ where $\phi$ is a diffeomorphism.) Important progress in this direction was made in the works of Collet, Eckmann, Epstein, Lanford and Martens [CEL], [E1], [E2], [EE], [Ma].

Let us also note that there is a parallel renormalization theory for circle dynamics which is also about to be completed (Lanford [La2], de Faria [dF], de Faria-de Melo [dFM]).

For background in "classical renormalization theory" (15 years old) see Collet and Eckmann [CE], Cvitanović [Cv] and Vul-Sinai-Khanin [VSK]. For more recent developments see de Melo - van Strien [MvS], McMullen [McM1] and the author [L3]. See also Tresser [T] for a lively historical retrospective.

Let us finally note that the Hairiness Theorem we prove here has a dynamical counterpart (hairiness of the Feigenbaum Julia sets) which was proven by McMullen [McM2].

1.7. *Further development.* We have recently proved the Renormalization Conjecture for all real combinatorial types [L5]. We conclude that the set of real infinitely renormalizable parameter values has zero linear measure, and that almost any real quadratic map $P_c : x \mapsto x^2 + c$, $c \in [-2, 1/4]$, has either an attracting cycle, or an absolutely continuous invariant measure.

1.8. *Notation and definitions.* As usual, $\mathbb{C}$ is the complex plane; $\mathbb{R}$ is the real line; $\mathbb{N} = \{0, 1, \ldots\}$ is the set of natural numbers; $\mathbb{Z}$ is the set of integers; $\mathbb{D}(a, r) = \{z : |z - a| < r\}$ is the open round disk of radius $r$; $\mathbb{D}_r \equiv \mathbb{D}(0, r)$, $\mathbb{D} \equiv D_1$;
$\mathbb{T}_r = \partial \mathbb{D}_r$ is the circle of radius $r$, $\mathbb{T} \equiv \mathbb{T}_1$;
$\mathbb{A}(r, R) = \{z : r < |z| < R\}$;
$\bar{X}$ denotes the closure of a set $X$;
$U \Subset V$ means that $U$ is *compactly contained* in $V$; that is, $\bar{U}$ is compact and is contained in $V$.
*A topological disk* means a simply connected domain in $\mathbb{C}$;
*A topological annulus* means a doubly connected domain in $\mathbb{C}$. The modulus of a topological annulus, mod$A$, is equal to $\log(R/r)$, provided $A$ is conformally equivalent to a round annulus $\mathbb{A}(r, R)$ (where $r = 0$ or $R = \infty$ are allowed).
The domain of a map $f$ is denoted by $\mathrm{Dom}(f)$;
Quasi-conformal maps will be abbreviated as "qc";
$\mathrm{Dil}(h)$ will stand for the dilatation of a qc map $h$;
$P_c(z) = z^2 + c$;
$M_0$ is the Mandelbrot set.
Given a map $f : X \to X$ and a point $x \in X$, let $\mathrm{orb}(x) = \{f^m x\}_{m=0}^{\infty}$, $\mathrm{orb}_n(x) = \{f^m x\}_{m=0}^{n}$;
$a \asymp b$ means that the ratio $a/b$ stays away from 0 and $\infty$;
$a \sim b$ means that $a/b \to 1$.

We assume the reader is familiar with the basic holomorphic dynamics (see e.g., [CG]) and the basic theory of hyperbolic dynamical systems (see e.g., [Sh]).

1.9. *Acknowledgement.* The initial part of this work was done during the Dynamics and Geometry Program at MSRI (spring 1995). My thought in the "unstable direction" was inspired by computer pictures concerning self-similarity of the Mandelbrot set shown to me by Dierk Schleicher at that time.
   It is also my pleasure to thank: Yu. S. Ilyashenko and Yu. I. Lyubich for helpful discussions of the central manifold theorem and operator theory in Banach spaces; J. Mather, C. McMullen and W. de Melo for many valuable comments on the manuscript; H. Epstein, M. Feigenbaum, O. Lanford, D. Sullivan, and C. Tresser for important historical comments; B. Hinkle for providing a reference for Lemma 5.1 and for making computer pictures.
   I am grateful to IHES (France) and UNAM in Cuernavaca (Mexico), where I completed writing this paper, for their hospitality. Thanks are also due to SUNY at Stony Brook for the generous support during the MSRI program and my sabbatical year. The work was partially supported by the Sloan Research Fellowship and NSF (grants DMS-8920768 (MSRI), DMS-9022140, DMS-9505833).

The main results of this paper were announced, among many other places, at the Paris-Orsay Conference (July 1995) and in the survey [L3].

## 2. Slow small orbits

2.1. *The one-dimensional case.* In the one-dimensional situation the Small Orbit Theorem says that any analytic map $R: z \mapsto e^{2\pi i\theta}z + bz^2 + \ldots$ near the origin has small orbits. This situation is well understood. There are three possible cases:

- *The parabolic case* when $\theta = q/p$ is rational. In this case $R$ is either of finite order, that is, $R^p = \mathrm{id}$, or there exist orbits converging to 0 (within the attracting petals).

- *The Siegel case* when $R$ is conformally equivalent to the rotation $z \mapsto e^{2\pi i\theta}z$. In this case all orbits which start sufficiently close to 0 do not escape a small neighborhood of 0.

- *The Cremer case* (neither of the above). In this case, for all sufficiently small $\varepsilon > 0$, the connected component $K_\varepsilon$ of the set
  $$\{z: |R^n z| \leq \varepsilon,\ n = 0, 1, \ldots\}$$
  is a continuum intersecting the circle $\mathbb{T}_\varepsilon$ (see Birkhoff [B, p. 95] and Perez-Marco [PM]).

Thus in all three cases there exist small orbits.

2.2. *The size of the basin of attraction.* We will consider a slightly more general setting than needed for the Small Orbits Theorem which will be suitable for further applications to bounded combinatorics. Let us consider a Banach space $\mathcal{B}$ split into two subspaces: $\mathcal{B} = E^s \oplus E^c$.

Let $D^s = D^s(\delta)$ and $D^c = D^c(\delta)$ stand for the open disks of radius $\delta$ centered at 0 in $E^s$ and $E^c$ respectively. Let us consider the bidisk $D = D(\delta) = D^s(\delta) \times D^c(\delta)$. Let $\partial^c D$ stand for $D^s \times \partial D^c$, and $\partial^s D$ be similar.

For $h \in \mathcal{B}$, let $h^s$ and $h^c$ denote the horizontal and vertical components of $h$, i.e, the projections of $h$ onto $E^s$ and $E^c$ respectively. Define the angle $\theta(h) \in [0, \pi/2]$ (between $h$ and $E^s$) by the condition:

$$\mathrm{tg}\theta(h) = \frac{\|h^c\|}{\|h^s\|}.$$

Let $C_f^\theta = \{h \in \mathrm{T}_f \mathcal{B} : h \neq 0, \theta(h) \geq \theta\}$ stand for the tangent cone with angle $\pi/2 - \theta$ about its axis $E^c$.

Let us also have another Banach space $\mathcal{B}'$ compactly embedded into $\mathcal{B}$; i.e., there is a linear injection $i : \mathcal{B}' \to \mathcal{B}$ such that the image of the unit ball of $\mathcal{B}'$ is relatively compact in $\mathcal{B}$.

Let us now consider a periodic point situation. Let us have $p$ pairs of complex Banach spaces $(\mathcal{B}_j, \mathcal{B}'_j)$ as above, and $p$ maps $T_j : (\mathcal{U}_j, 0) \to (\mathcal{B}'_{j+1}, 0)$ defined in some neighborhoods $\mathcal{U}_j \subset \mathcal{B}_j$ of the origin $0 = 0_j$ (where the index $j$ is considered $\bmod\, p$). We will naturally label all the above objects with the subscript $j$: $E_j^s$, $D_j(\delta)$, $i_j$ etc. Let $\mathcal{B} = \sqcup \mathcal{B}_j$, $\mathcal{U} = \sqcup \mathcal{U}_j$, $D = \sqcup D_j$ etc., and $T : \sqcup \mathcal{U} \to \mathcal{B}'$ be the operator acting as $T_j$ on $\mathcal{U}_j$. Let

$$R = i \circ T : \mathcal{U} \to \mathcal{B}, \quad R|\mathcal{U}_j = i_{j+1} \circ T_j.$$

LEMMA 2.1 (Basin of attraction). *Given the spaces and operators as above, assume that $T : \mathcal{B} \to \mathcal{B}'$ is complex analytic and $i : \mathcal{B}' \to \mathcal{B}$ is compact. Assume that the decomposition $\mathcal{B} = E^s \oplus E^c$ is invariant with respect to the differential $\mathrm{D}R(0)$, $\mathrm{D}R^p(0) | E^s \neq 0$, and moreover the following properties are satisfied*:

H0. *The origin is attracting*: $\mathrm{spec} \mathrm{D}R^p(0) \subset \mathbb{D}$.

H1. *Horizontal contraction*: $Rf \notin \partial^s D$ *for* $f \in D$ *and*

$$\|(\mathrm{D}R_f h)^s\| \leq q\|h\| \quad provided \quad f \in D,\ Rf \in D.$$

H2. *Invariant cone field*: *There exists a $\theta \in (0, \pi/2)$ such that the tangent cone field $C_f^\theta$ over $D$ is $R$-invariant*:

$$(\mathrm{D}R_f) C_f^\theta \subset C_{Rf}^\theta \quad provided \quad f \in D,\ Rf \in D.$$

*Then there is a point $f \in \partial^c D$ such that $\mathrm{orb}(f) \subset D$ and $\|R^m f\| \to 0$.*

*Remark.* Note that there are no assumptions relating the spectrum of $\mathrm{D}R^p(0)$ and the size of the bidisks $D_j$.

*Proof.* By assumption H0, the origins $0_j \in \mathcal{B}_j$ form an attracting cycle **O** of period $p$. Let us consider its basin of attraction in $\bar{D}$:

$$A = \{f : R^n f \in \bar{D},\ n = 0, 1, \ldots \quad \text{and} \quad \|R^n f\| \to 0\ as\ n \to \infty\}.$$

Clearly $A$ is forward invariant. We want to show that $A$ intersects the horizontal boundary $\partial^c D$. Assume this is not the case.

Let $A^o = A \cap D$. Then $A^o$ is forward invariant (indeed, if $f \in A^o$ then $Rf \notin \partial^c D$ by the assumption, and $Rf \notin \partial^s D$ by the assumption H1). It follows that $A^o$ is open.

Let $\partial^c A$ denote the part of the boundary of $A$ which does not belong to $\partial^s D$. Then

$$R(\partial^c A) \subset \partial^c A. \tag{2.1}$$

Otherwise there would be a point $f \in \partial^c A$ such that $Rf \in A^o$ and hence $f \in A$. As $A^o$ is open, $f$ can belong simultaneously to $A$ and $\partial A$ only if $f \in \partial D$ contradicting the assumptions.

We are going to show that the assumptions H0–H2 contradict (2.1). Since all these properties are inherited by the iterates, at this point we can replace $R$ with $R^p$ and assume without loss of generality that $p = 1$.

Note first that by the invariant cone field assumption H2, the linear operator $DR(0)|E^c$ does not have 0 in its spectrum. Since it is compact, $E^c$ is finite dimensional, $\dim E^c = d$.

Let us now consider a family $\mathcal{G}$ of immersed analytic manifolds $\psi : (\Omega, 0) \to (\Gamma, 0)$, where $\Omega = \Omega_\psi \subset \mathbb{C}^d$ and $\Gamma = \Gamma_\psi \subset A^o$, with the following properties:

A1. The tangent spaces $T_f \Gamma \equiv DT_z\Omega$, where $f = \psi(z)$, belong to the cones $C_f^\theta$.

A2. The manifolds are properly immersed into $A^o$. This means that if a curve $\gamma(t) \subset \Omega$, $0 < t < \infty$, tends to $\infty$ in $\Omega$ as $t \to \infty$ (i.e., it eventually escapes any compact subset of $\Omega$), then $\psi(\gamma(t))$ tends to $\partial^c A$.

*Remarks.* 1. The family $\mathcal{G}$ is nonempty: just let $\Gamma$ be the connected component of the vertical slice $D^c \cap A^o$ of the attracting basin.

2. By Property A1, the projection $P : \Gamma \to D^c$ of any $\Gamma \in \mathcal{G}$ onto the vertical subspace is nonsingular. Moreover, for any tangent vector $v \in T_f\Gamma$, $\|Pv\| \asymp \|v\|$.

3. Property A1 also implies that there is an $\varepsilon > 0$ such that any irreducible component $\Gamma(\varepsilon)$ of $(D^s \times D^c(\varepsilon)) \cap \Gamma$ containing 0 is a graph of an analytic map $\phi : D^c(\varepsilon) \to D^s$.

Let us supply the manifolds $\Gamma \in \mathcal{G}$ with the Kobayashi metrics. Recall that the Kobayashi norm of a tangent vector $v \in T_f\Gamma$ is defined as

$$\|v\|_\Gamma = \inf_\gamma \|w\|_P,$$

where $\|w\|_P$ stands for the Poincaré norm of a vector $w \in T_0\mathbb{D}$, and the infimum is taken over all holomorphic curves $\gamma : (\mathbb{D}, w) \to (\Gamma, v)$ (where by definition such a curve is factored via the parametrization $\psi : \Omega \to \Gamma$). The Kobayashi metric is invariant under holomorphic coverings and increases under shrinking the manifold.

*Remark.* A covering map between immersed manifolds is naturally defined using the parametrizations. A holomorphic covering of finite degree can be defined as a proper nonsingular holomorphic map.

It follows that for a tangent vector $v \in T_0\Gamma$, $\Gamma \in \mathcal{G}$, the Kobayashi norm is uniformly subordinate to the Banach one:

$$\|v\|_\Gamma \leq C\|v\|, \tag{2.2}$$

where the constant $C$ is independent of $v$. Indeed

$$\|v\|_\Gamma \leq \|v\|_{\Gamma(\varepsilon)} = \|Pv\|_{D^c(\varepsilon)}.$$

On the other hand, by Remark 2 above,

$$\|v\| \asymp \|Pv\| \asymp \|Pv\|_{D^c(\varepsilon)}.$$

Let us now consider a manifold transformation $R_* : \Gamma_\psi \mapsto \Gamma_{R\circ\psi}$. By the invariant cone field assumption H2 and (2.1), $R \circ \psi : \Omega_\psi \to A^o$ is an immersion satisfying properties A1, A2. Hence $R_*$ transforms $\mathcal{G}$ into itself. Moreover, the map $R : \Gamma \to R\Gamma$ is proper and nonsingular, and hence is a holomorphic covering of finite degree.

Since the Kobayashi metric is invariant under holomorphic coverings, for any tangent vector $v = D\psi_0(w) \in T_0\Gamma$ we have:

$$\|DR^n(v)\|_{R_*^n\Gamma} = \|v\|_\Gamma.$$

On the other hand, since 0 is an attracting point,

$$\|DR^n(v)\| \to 0 \quad \text{as} \quad n \to \infty.$$

These last two estimates contradict (2.2). $\square$

2.3. *Proof of the Small Orbits Theorem.* We are now ready to prove the Small Orbits Theorem stated in the introduction. The argument below exploits Perez-Marco's perturbation idea [PM].

Let $E^s$ stand for the spectral subspace of $R$ corresponding to the part of spec$R$ inside the unit disk $\mathbb{D}$, and let $E^c$ correspond to the part on the unit circle $\mathbb{T}$. After we replace $R$ by its iterate (or after adaptation of the Banach norm), $R$ becomes horizontally contracting and cone field preserving on a sufficiently small bi-disk $D = D(\delta)$.

For $\lambda \in (0, 1)$, let us consider the perturbation $R_\lambda = \lambda R$ which makes the origin attracting. This operator is more strongly horizontally contracting than $R$ and preserves the same cone field. Thus it satisfies assumptions H0–H2 of Lemma 2.1. Hence there is a point $f_\lambda \in \partial^c D \cap A_\lambda$, where $A_\lambda = \{f : R_\lambda^m f \in \bar{D},\ m = 0, 1, \ldots, \|R_\lambda^m f\| \to 0\}$ is the basin of 0.

Since the set $\{R_\lambda f_\lambda\}$ is pre-compact in $\mathcal{B}$, there is a convergent subsequence $R_{\lambda_n} f_{\lambda_n} \to g$ as $\lambda_n \to 1$. Clearly $g_m \equiv R^m g \in \bar{D}$, $m = 0, 1, \ldots$. Moreover, for $\delta$ sufficiently small, $\|g^c\| \geq \|g^s\|$, since $R_\lambda$ contracts more strongly in the $E^s$-direction than in the $E^c$-one. Since $E^c$ is the neutral direction,

$$\|g_{m+1}^c\| = \|g_m^c\|\,(1 + o(\|g_m^c\|)), \tag{2.3}$$

provided

(2.4) $$\|g_m^s\| \leq \|g_m^c\|.$$

But (2.4) is inductively satisfied for all $m$, so that (2.3) is satisfied for all $m$ as well. This implies that the $g_m$ may not exponentially converge to 0. □

*Remark.* The Small Orbits Theorem is still true if $R$ is allowed to have spectrum outside the unit disk.

## 3. External maps and hybrid classes

3.1 *Quadratic-like maps and germs.* The following fundamental notion was introduced by Douady and Hubbard [DH2]. A *quadratic-like map* $f : U \to U'$ is a holomorphic double branched covering (i.e., a proper map of degree 2) between topological disks $U, U'$ such that $U \Subset U'$. It has a single critical point which is assumed to be located at the origin 0, unless otherwise stated. Note that the restriction of a quadratic polynomial $P_c$ on a disk $\mathbb{D}_r$ of a sufficiently big radius $r$ is a quadratic-like map.

The filled Julia set of a quadratic-like map is defined as the set of nonescaping points:
$$K(f) = \{z : f^n z \in U, n = 0, 1 \ldots\}.$$

Its boundary is called the *Julia set*, $J(f) = \partial K(f)$. The Julia set $J(f)$ (and also $K(f)$) is connected if and only if the critical point itself is nonescaping: $0 \in K(f)$. Otherwise it is a Cantor set.

Any quadratic-like map has two fixed points, $\alpha$ and $\beta$ (counted with multiplicity). In the case of connected Julia sets these points can be distinguished. Namely, if these points are different then $\beta$ is a *nondividing* repelling point (i.e., its removal does not disconnect the Julia set), while $\alpha$ is either non-repelling, or *dividing*.

In what follows we will also assume (without loss of generality) that the domains $U$ and $U'$ of a quadratic-like map $f$ are bounded by *smooth Jordan curves*. The topological annulus $A = U' \setminus \bar{U}$ is called the *fundamental annulus* of $f$. Let us foliate the fundamental annulus by topological circles with winding number 1 around the origin (e.g., using the Riemann mapping onto a round annulus). This foliation can be pulled back by dynamics providing a foliation $\Phi$ with singularities on $U' \setminus K(f)$. Moreover, if $K(f)$ is connected, then $\Phi$ is nonsingular. Otherwise, the outermost singular leaf is a figure eight passing through the critical point 0. The leaves of this foliation will be called *equipotentials*.

We will consider quadratic-like maps up to *affine* conjugacy (rescaling), so that near the origin they can be normalized as $f : z \mapsto c + z^2 + \ldots$. Let $\mathcal{QM}$ stand for the set of normalized quadratic-like maps.

Two analytic maps $f$ and $g$ defined near 0 *represent the same germ at* 0 if they coincide in some neighborhood of 0. Let $\mathcal{G}_0$ stand for the space of germs at the origin. By taking all possible analytic continuations, a germ $f \in \mathcal{G}_0$ can be equivalently viewed as a *full analytic function* $f : S_f \to \mathbb{C}$ defined on a Riemann surface $S_f$ covering $\mathbb{C}$ (in general, nonevenly).

Let us now define quadratic-like germs. To this end consider the following relation on the space $\mathcal{QM}$: $f \sim \tilde{f}$ if the maps $f$ and $\tilde{f}$ have a common quadratic-like restriction. The following result gives a useful criterion for two maps to be in this relation:

PROPOSITION 3.1 ([McM1, §5.4]).  *Consider two quadratic-like maps $f : U \to U'$ and $\tilde{f} : \tilde{U} \to \tilde{U}'$ representing the same germ at 0. Let $W$ be the component of $U \cap \tilde{U}$ containing 0. Then $f \sim \tilde{f}$ if and only if $fW \ni 0$. Moreover, in this case the restriction $g = f|W$ is a quadratic-like map, and $K(f) = K(g) = K(\tilde{f})$.*

This lemma yields:

- If $f \sim g$ then $K(f) = K(g)$.

- For the maps with connected Julia sets, the above relation is an equivalence relation.

In general, let us consider the equivalence relation generated by $\sim$. Classes of this equivalence relation are called *quadratic-like germs*. Thus two quadratic-like maps $f$ and $\tilde{f}$ represent the same quadratic-like germ $[f] = [\tilde{f}]$ if there is a string of quadratic-like maps $f = f_0, f_2, \ldots, f_N = \tilde{f}$ such that $f_k$ and $\tilde{f}_k$ have a common quadratic-like restriction.

Let $\mathcal{QG}$ stand for the space of quadratic-like germs. There is a natural map $j : \mathcal{QG} \to \mathcal{G}_0$. If $g = j(f)$, one can say that $g$ is *marked*. Note that there exists at most one marking with connected Julia set [McM2, Lemma 7.1]. Also, any quadratic polynomial $P_c$ has a unique marking since any two quadratic-like restrictions $U \to U'$ and $V \to V'$ have a "common minorant" $\mathbb{D}_r \to P_c(\mathbb{D}_r)$ (such that $\mathbb{D}_r \supset U \cup V$). Thus the quadratic family $P_c$, $c \in \mathbb{C}$, is naturally embedded into $\mathcal{QG}$.

Furthermore, two markings with "close almost connected" Julia sets must coincide. More precisely, let us consider two quadratic-like maps $f : V \to V'$ and $\tilde{f} : \tilde{V} \to \tilde{V}'$ with the same germ at 0 whose Julia sets stay Hausdorff distance at most $\varepsilon$ apart. Let $\Omega_\varepsilon(J(f))$ denote the filled closure of the $\varepsilon$-neighborhood of $J(f)$ (where "filling" means adding bounded components of the complement). Assume that $\Omega_\varepsilon(J(f))$ is connected and

(3.1) $$\{0, f(0)\} \subset \Omega_\varepsilon(J(f)) \subset \Omega_{3\varepsilon}(J(f)) \subset V',$$

and the same holds for $\tilde{f}$. Then $f \sim \tilde{f}$. (Indeed, $\Omega_\varepsilon(J(f)) \cup \Omega_\varepsilon(J(\tilde{f}))$ is a connected set contained in $V \cap \tilde{V}$ and containing $\{0, f(0)\}$.)

Quadratic-like germs will be considered up to *affine* conjugacy. We will adopt the following notational and terminological conventions (except for special situations when they may lead to confusion). The germ of a quadratic-like map $f$ will still be denoted by $f$. However, sometimes the notation $f_V$ will be used for the quadratic-like representative $f_V : V \to V'$ of a germ $f$. The germ of a polynomial will still be called "polynomial."

Given a quadratic-like germ $f$, let

$$\mathrm{mod}(f) = \sup \mathrm{mod}(U' \setminus U),$$

where the supremum is taken over all quadratic-like representatives $f : U \to U'$ of the germ. Note that $\mathrm{mod}(f) = \infty$ if and only if $f$ is a polynomial.

Let $\mathcal{C}$ be the *connectedness locus* of $\mathcal{QG}$, that is, the subset of quadratic-like germs with connected Julia set. Let $\mathcal{QG}^\#$ stand for the set of quadratic-like germs which have representatives satisfying (3.1).

3.2. *External maps.* Let $g : \mathbb{T} \to \mathbb{T}$ be a degree two real analytic endomorphism of the unit circle $\mathbb{T}$. It can also be viewed as a complex analytic germ near the circle. We call $g$ *expanding* if it admits an analytic extension to a double covering $g : V \to V'$ between annular neighborhoods of $\mathbb{T}$ such that $V \Subset V'$. We consider such a map up to conjugacy by rotation, which is equivalent to normalizing it in such a way that $1 \in \mathbb{T}$ is a fixed point. Let $\mathcal{E}$ stand for the space of circle endomorphisms as above (up to rotation). Let

$$\mathrm{mod}(g) = \sup \mathrm{mod}(V' \setminus (V \cup \mathbb{D})),$$

where the supremum is taken over all extensions $g : V \to V'$ as above.

There is a projection $\pi : \mathcal{QG} \to \mathcal{E}$ which associates to $f \in \mathcal{QG}$ its *external map* $g \in \mathcal{E}$. In the case when $f \in \mathcal{C}$, the construction is easily provided by the Riemann Mapping Theorem. Namely, let $f : U \to U'$ be a quadratic-like representative of the germ. Let us conjugate $f : U \setminus K(f) \to U' \setminus K(f)$ by the Riemann mapping

$$\phi = \phi_f : \mathbb{C} \setminus K(f) \to \mathbb{C} \setminus \mathbb{D}$$

to a double covering $g : W \to W'$ between annuli with inner boundary $\mathbb{T}$. By the Reflection Principle, $g$ extends to a circle endomorphism of class $\mathcal{E}$. Since the Riemann mapping $\phi$ is well-defined up to post-composition with rotation, $g$ is well-defined up to conjugacy by rotation.

LEMMA 3.2. *Let $f \in \mathcal{C}$ and $g = \pi(f)$ be its external map. Then* $\mathrm{mod}(f) = \mathrm{mod}(g)$.

*Proof.* The Riemann mapping $\phi : \mathbb{C} \setminus K(f) \to \mathbb{C} \setminus \mathbb{D}$ obviously establishes one-to-one correspondence between the fundamental annuli of $f$ and $g$. □

In the case of a disconnected Julia set the construction is more subtle. Take a closed fundamental annulus $A = \overline{U' \setminus U}$ with real analytic boundary curves $E = \partial U'$ and $I = \partial U$. Then $f : I \to E$ is a real analytic double covering.

Let $\mu = \mathrm{mod} A$ and consider an abstract double covering $\xi_1 : A_1 \to A$ of an annulus $A_1$ of modulus $\mu/2$ over $A$. Let $I_1$ and $E_1$ be the "inner" and "outer" boundary of $A_1$; i.e., $\xi_1$ maps $I_1$ onto $I$ and $E_1$ onto $E$. Then there is a real analytic diffeomorphism $\theta_1 : E_1 \to I$ such that $\xi_1 = f \circ \theta_1$. This allows us to stick the annulus $A_1$ to the domain $\mathbb{C} \setminus U$ bounded by $I$. We obtain a Riemann surface $T_1 = (\mathbb{C} \setminus U) \cup_{\theta_1} A_1$. Moreover, the maps $f$ on $I$ and $\xi_1$ on $A_1$ match to form an analytic double covering $f_1 : A_1 \to A$.

This map $f_1$ restricts to a real analytic double covering of the inner boundary of $A_1$ onto its outer boundary. This allows us to repeat this procedure: we can attach to the boundary of $T_1$ an annulus $A_2$ of modulus $\frac{1}{4}\mu$, and extend $f_1$ to the new annulus $A_2$. Proceeding in this way, we will construct a Riemann surface

$$(3.2) \qquad T^A \equiv T^A(f) = \lim T_n = (\mathbb{C} \setminus U) \cup_{\theta_1} A_1 \cup_{\theta_2} A_2 \ldots$$

and an analytic double covering $F : \cup_{n \geq 1} A_n \to \cup_{n \geq 0} A_n$ extending $f$.

Since the trajectories of $F$ do not converge to the "inner" ideal boundary of $T^A$, it is a punctured disk and can be conformally mapped onto $\mathbb{C} \setminus \mathbb{D}$. Now by the reflection principle, this conformal representation of $F$ can be extended to an analytic expanding endomorphism $g_A : \mathbb{T} \to \mathbb{T}$.

For a given choice of the fundamental annulus $A$, the map $g_A : V \to V'$ (which comes together with the domains $(V, V')$) is well-defined up to rotation. Indeed, for two such maps $g_A$ and $\tilde{g}_A$, by construction there is a conformal isomorphism $h : \mathbb{C} \setminus \mathbb{D} \to \mathbb{C} \setminus \mathbb{D}$ conjugating them on an outer neighborhood of the circle. Reflecting $h$ to the unit disk, we conclude that $h$ is a rotation conjugating $g_A$ and $g'_A$ near the circle.

The endomorphism $g_A : \mathbb{T} \to \mathbb{T}$ does not actually depend on $A$, that is, on the choice of a representative of the quadratic-like germ. Indeed, let $A = U' \setminus U$ and $\tilde{A} = \tilde{U}' \setminus \tilde{U}$ be two fundamental annuli of $f$ such that $0$ and $f(0)$ are contained in the same connected component $W'$ of $U' \cap V'$. Then by Proposition 3.1, we have a quadratic-like map $f : W \to W'$ whose fundamental annulus $B = W' \setminus W$ is contained in both $U'$ and $V'$. Let us show that in this case $g_A = g_B = g_{\tilde{A}}$.

It is enough to check the first equality. Let us consider the tautological embedding $\mathbb{C} \setminus U \to \mathbb{C} \setminus W$. The latter Riemann surface is naturally embedded in $T^B$, so that we obtain an embedding $i : \mathbb{C} \setminus U \to T^B$ conjugating $f|\partial U$ to $g_B|i(\partial U)$. (Note that $g_B$ admits an analytic extension up to the curve $i(\partial U)$ provided by $f|U \setminus W$.)

Let us attach an abstract annulus $A_1$ of half modulus to the inner boundary of $A$. Then $i$ is naturally lifted to this annulus, so that we obtain a conformal embedding $i : (\mathbb{C} \setminus U) \cup A_1 \to T^B$. Now attach the next annulus $A_2$ and extend $i$ to the new Riemann surface, etc. At the end we obtain a conformal isomorphism $i : T^A \to T^B$ conjugating $g_A$ to $g_B$ near the ideal boundary. The conclusion is obvious.

Thus the external map $g \equiv \pi(f)$ is correctly defined up to rotation.

The above construction is due to Douady and Hubbard [DH2] except that we remember the embedding of the fundamental annulus $A$ into $\mathbb{C} \setminus U$. Due to this the external map is defined up to rotation rather than an analytic diffeomorphism.

We say that two quadratic-like germs $f$ and $\tilde{f}$ are *externally equivalent* if $\pi(f) = \pi(\tilde{f})$. Note that in the case of connected Julia sets this means that the appropriately normalized conformal map $h : \mathbb{C} \setminus K(f) \to \mathbb{C} \setminus K(\tilde{f})$ conjugates $f$ and $\tilde{f}$ near the Julia sets.

LEMMA 3.3. *The external map $\pi(f)$ is equal to $P_0 : z \mapsto z^2$ if and only if $f$ is a quadratic polynomial $P_c : z \mapsto z^2 + c$.*

*Proof.* The external map $g$ of a quadratic polynomial acts as a double covering on the whole Riemann surface $T \equiv \mathbb{C} \setminus \mathbb{D}$, and hence $g(z) = z^2$.

Vice versa, assume that $\pi(P_c) = z^2$. This means that there is a conformal map $\phi : \mathbb{C} \setminus U \to \mathbb{C} \setminus (V \cup \mathbb{D})$ such that $\phi(f(z)) = \phi(z)^2$ whenever this makes sense. But this functional equation allows us to extend $f$ analytically to the whole complex plane, so that $f$ is a restriction of a quadratic polynomial. □

The inverse map $\phi_f^{-1} : \mathbb{C} \setminus V \to \mathbb{C} \setminus U$ conjugating $g = \pi(f)$ on $\partial V$ to $f$ on $\partial U$ (and perhaps analytically extended elsewhere) will be called the *uniformization of $f$ at $\infty$*.

The fibers $\mathcal{Z}_g$ of the projection $\pi$ will be called *vertical curves*, or *vertical fibers*.

3.3. *The Riemann mapping.* Let us now construct a smooth map

$$\xi : \mathcal{QG} \setminus \mathcal{C} \to \mathbb{C} \setminus \bar{\mathbb{D}}$$

which conformally uniformizes the vertical fibers (its restriction to the quadratic family will coincide with the Riemann mapping $\mathbb{C} \setminus M_0 \to \mathbb{C} \setminus \bar{\mathbb{D}}$).

Let us say that a normalized map $g \in \mathcal{E}$, $g(1) = 1$, is *marked* if one has selected annuli neighborhoods $(V, V')$ of $\mathbb{T}$ according to the definition of an expanding map and a point $a \in V' \setminus \mathbb{D}$, up to the following equivalence relation. Two data $(g : V \to V', a)$ and $(g : \tilde{V} \to \tilde{V}', a)$ are considered to be equivalent if there is a string of representatives $g : V_k \to V'_k$, $k = 0, \ldots, N$, such that

- $(V_0, V'_0) = (V, V')$, $(V_N, V'_N) = (\tilde{V}, \tilde{V}')$;

- $\mathbb{T}$ and $a$ are contained in the same connected component of $V'_k \cap V'_{k+1}$, $k = 0, 1, \ldots, N - 1$.

In the case of a disconnected Julia set, the above construction of the external map $g = \pi(f)$ actually provides us with a marked map as follows. Given an $(f : U \to U') \in \mathcal{QG} \setminus \mathcal{C}$, let us consider its external map $(g : V \to V') = \pi(f) \in \mathcal{E}$, where $(V, V')$ naturally corresponds to $(U, U')$. By the construction of $g$, there is a conformal map $\phi = \mathbb{C} \setminus U \to \mathbb{C} \setminus (V \cup \mathbb{D})$ conjugating $f : \partial U \to \partial U'$ to $g : \partial V \to \partial V'$. Let $N$ be the maximal natural number such that $f^{-N}U \ni 0$. Then $\phi$ admits the analytic extension to the domain $\mathbb{C} \setminus f^{-N}U \ni f(0)$. Thus we can mark the point

$$a = \xi(f) = \phi(f(0))$$

in $V'$.

Selecting a different representative of $f$ does not change marking of $g$. Indeed, let $f : U \to U'$ and $f : \tilde{U} \to \tilde{U}'$ satisfy the property that $0$ and $f(0)$ belong to the the same component of $U' \cap \tilde{U}'$. Then considering the quadratic-like intersection $f : W \to W'$ we conclude that $\mathbb{T}$ and $a$ belong to the same component of $V' \cap \tilde{V}'$.

Let $\mathcal{E}^m$ denote the space of marked circle maps of class $\mathcal{E}$, and $\pi^m$ denote the projection $\mathcal{QG} \setminus \mathcal{C} \to \mathcal{E}^m$ just described. Let us also consider the following natural maps:

(3.3) $\qquad\qquad \zeta : \mathcal{E}^m \to \mathcal{E}$   (forgetting the marked point)

and

(3.4) $\qquad\qquad \eta : \mathcal{E}^m \to \mathbb{C} \setminus \bar{\mathbb{D}}$   (position of the marked point).

Let $\mathcal{S}_g$ denote the fibers of $\zeta$, $g \in \mathcal{E}$. Endow $\mathcal{S}_g$ with the following topology. Pick an $\varepsilon > 0$ and a representative $(g : V \to V', a)$ of a marked map $G \in \mathcal{S}_g$. Then a neighborhood of $G$ in $\mathcal{S}_g$ consists of marked maps represented by $(g : V \to V', b)$ with $b \in V'$, $|b - a| < \varepsilon$. Then $\eta : \mathcal{S}_g \to \mathbb{C} \setminus \bar{\mathbb{D}}$ is a local homeomorphism. Pulling the complex structure back via $\eta$, we make the fiber a Riemann surface covering (nonevenly) an outer neighborhood of the unit circle.

As in the case of $\mathcal{QG}$, there is a subset $\mathcal{E}^\# \subset \mathcal{E}^m$ with a preferred marking. Namely let $\mathcal{E}^\#$ consist of marked circle germs which have representatives $(g : V \to V', a)$ with the property that $\partial V'$ and $a$ are separated by some round circle $\mathbb{T}_r \subset V' \setminus \mathbb{D}$. Clearly this marking is uniquely determined by the germ $g$ and the point $a$. Thus the projection $\eta : \mathcal{S}_g^\# \to \mathbb{C} \setminus \bar{\mathbb{D}}$ is univalent (here $\mathcal{S}_g^\# = \mathcal{S}_g \cap \mathcal{E}^\#$) and can be identified with a domain $\Omega_g$ of the complex plane.

Note that for any $g \in \mathcal{E}$, there is a representative $g : V \to V'$ such that $V' \setminus \mathbb{D} \supset \mathbb{T}_r$ with $r = r(g) > 1$ depending only on $\mathrm{mod}(g)$. Hence any point $a \in \mathbb{A}(1, r)$ specifies a preferred marking of $g$, so that $\Omega_g \supset \mathbb{A}(1, r)$.

THEOREM 3.4 (External mating). *Let $g \in \mathcal{E}^m$ be a marked circle map. Then there is a unique quadratic-like germ $f = \theta(g) \in \mathcal{QG} \setminus \mathcal{C}$ such that $\pi^m(f) = g$.*

*Proof.* Let us consider the marked point $a = \eta(g)$. Let $g^N a \in A \equiv V' \setminus (V \cup \mathbb{D})$ (without loss of generality we can assume that $g^N a \notin \partial A$). On the other hand, take any quadratic polynomial $P = P_c$ with Cantor Julia set. Let $G$ be the Green function on $\mathbb{C} \setminus K(P)$ with pole at $\infty$, and let $\Omega_r = \{z : G(z) < \log r\}$ ($r > 1$). Select a fundamental annulus $\Omega_{R^2} \setminus \Omega_R$ containing $P^N c$. Let $r = R^{1/2^N}$ and $U = g^{-N} V$.

There exists a diffeomorphism $\phi : (\mathbb{C} \setminus \Omega_r, c) \to (\mathbb{C} \setminus U, a)$ conjugating $P|(\Omega_R \setminus \Omega_r)$ to $g|(V \setminus U)$. Let us consider the conformal structure $\mu = \phi^* \sigma$ on $\mathbb{C} \setminus \Omega_r$, and pull it back by $P$ to the complement of $J(P)$. (Here $\sigma$ is the standard structure on $\mathbb{C} \setminus U$; see Appendix 1.) Straightening this conformal structure by the Measurable Riemann Mapping Theorem we obtain a desired quadratic-like map $f$.

Let us have two quadratic-like maps $f : U \to U'$ and $\tilde{f} : \tilde{U} \to \tilde{U}'$ in $\mathcal{QG} \setminus \mathcal{C}$ with $\pi^m(f) = \pi^m(\tilde{f}) = g$, in particular, $\xi(f) = \xi(\tilde{f}) = a$. Let us show that these maps represent the same quadratic-like germ.

First assume that the domains $(V, V')$ of $g$ corresponding to $U$ and $\tilde{U}$ coincide. Then $f$ and $\tilde{f}$ have fundamental annuli $U' \setminus U$ and $\tilde{U}' \setminus \tilde{U}$ whose inner boundaries are "figures eight" passing through 0, and such that there is a conformal map $\psi : (\mathbb{C} \setminus U, 0) \to (\mathbb{C} \setminus \tilde{U}, 0)$ conjugating $f$ to $\tilde{f}$ on the boundaries of these domains. This map admits a dynamical analytic extension to the complements of the Julia sets. Since the Julia sets are removable, $\psi$ is affine.

Assume now that we have two representatives of the marked germ, $g : V \to V'$ and $g : \tilde{V} \to \tilde{V}'$, such that $\mathbb{T}$ and $a$ belong to the same component $\Omega'$ of $V' \cap \tilde{V}'$. Let $\Omega = g^{-1} \Omega'$. Then $\Omega \Subset \Omega'$ and $g : \Omega \to \Omega'$ is a double covering.

Assume also that the map $g : V \to V'$ corresponds to a quadratic-like map $f : U \to U'$. Then there is a restriction $f : W \to W'$ corresponding to $g : \Omega \to \Omega'$. Moreover, by means of the functional equation $\phi(fz) = g(\phi(z))$, $f$ analytically extends to a domain $\tilde{U} \to \tilde{U}'$ corresponding to $g : \tilde{V} \to \tilde{V}'$.

Let us now consider a string of circle maps $g : V_k \to V_k'$, $k = 0, \ldots, N$, such that $\mathbb{T}$ and $a$ are contained in the same component of $V_k' \cap V_{k+1}'$, and $(V_0, V_0') = (V, V')$, $(V_N, V_N') = (\tilde{V}, \tilde{V}')$. Then by the above construction, we have the corresponding string of quadratic-like maps $f : U_k \to U_k'$ such that the $U_k' \cap U_{k+1}'$ contain both 0 and $f(0)$, and $(U_0, U_0') = (U, U')$. Then $f : U \to U'$ and $f : U_N \to U_N'$ represent the same germ. On the other hand, $f_N : U_N \to U_N'$ and $\tilde{f} : \tilde{U} \to \tilde{U}'$ have the same marked external map $g : \tilde{V} \to \tilde{V}'$ considered with its domain. As we have shown above, $f_N$ is affinely equivalent to $\tilde{f}$. □

The above operation $\theta$ will be called the *external mating*. The reason is that one can think of it as the mating of a circle map $g \in \mathcal{E}$ with a point $a \in \mathbb{C} \setminus \bar{\mathbb{D}}$ which produces a quadratic-like germ $f$ outside the connectedness locus $\mathcal{C}$.

Note that the restriction $\xi : \mathbb{C} \setminus M \to \mathbb{C} \setminus \mathbb{D}$ of $\xi$ to the quadratic family coincides with the uniformization of $\mathbb{C} \setminus M$ tangent to id at $\infty$ (see [DH1]). The preimage of a round circle $\mathbb{T}_r$ under this uniformization is called the (*parameter*) *equipotential* of level $r$.

We will see that the map $\pi^m : \mathcal{Z}_g \setminus \mathcal{C} \to \mathcal{S}_g$ is a conformal isomorphism (see Lemma 4.14). The preimages of the round circles under $\xi = \eta \circ \pi^m$ will be called (external) equipotentials on $\mathcal{Z}_g$.

The external equipotentials are the traces of *equipotential hypersurfaces* in $\mathcal{QG}$, the preimages of the round circles under $\xi : \mathcal{QG} \setminus \mathcal{C} \to \mathbb{C} \setminus \bar{\mathbb{D}}$. We will see that they are codimension 1 smooth submanifolds in $\mathcal{QG}$ (see Lemma 4.14).

3.4. *Conjugacies.* Two quadratic-like maps $f : U \to U'$ and $\tilde{f} : \tilde{U} \to \tilde{U}'$ are called *topologically conjugate* if there exists a homeomorphism $h : (U', U) \to (\tilde{U}', \tilde{U})$ such that $h(fz) = \tilde{f}(hz)$, $z \in U$. Two quadratic-like germs $f$ and $\tilde{f}$ are called topologically conjugate if there is a choice of topologically conjugate quadratic-like representatives. A self-conjugacy $h$ of a map/germ is called its *automorphism*.

LEMMA 3.5 (see e.g. [L2, Lemma 10.4]). *Let $f$ be a quadratic-like germ with connected Julia set and $h$ be its automorphism. Then $h|J(f) = \mathrm{id}$.*

Two maps/germs are called *quasi-conformally/smoothly* etc. conjugate (or equivalent) if they admit a conjugacy $h$ with the corresponding regularity. If two maps/germs are qc conjugate with $\bar{\partial} h = 0$ almost everywhere on the filled Julia set, then $f$ and $\tilde{f}$ are called *hybrid equivalent*.

LEMMA 3.6. *Let $f \in \mathcal{C}$ and $g \in \mathcal{C}$ be two quadratic-like germs with $\mathrm{mod}(f) \geq \nu > 0$ and $\mathrm{mod}(g) \geq \nu > 0$. Assume that they are qc conjugate by a map $h$ with*

$$\mathrm{ess} - \sup_{z \in K(f)} \frac{|\bar{\partial} h(z)|}{|\partial h(z)|} \leq \kappa.$$

*Then there exist quadratic-like representatives* $f : U \to U'$ *and* $g : V \to V'$ *such that*:

(i) $\mathrm{mod}(U' \setminus U) \geq \mu(\nu) > 0$ *and* $\mathrm{mod}(V' \setminus V) \geq \mu(\nu) > 0$;

(ii) *These representatives are $K$-qc conjugate with* $K = K(\kappa, \nu)$.

*Proof.* A quadratic-like germ $f \in \mathcal{C}$ with $\mathrm{mod}(f) \geq \nu > 0$ admits a representative $f : U \to U'$ satisfying (i) and such that the boundaries $\partial U, \partial U'$ are smooth $\gamma(\nu)$-quasi-circle (see e.g. [McM2, Prop. 4.10]). If we have two maps $f : U \to U'$ and $g : V \to V'$ with this property then there is a $Q$-qc map $H : U' \setminus U \to V' \setminus V$ respecting the boundary dynamics, with $Q = Q(\mu, \gamma)$. This map can be pulled back to a $K$-qc conjugacy between between $f$ and $g$ on the complements of the filled Julia sets. This conjugacy glues with $h|K(f)$ to a single $K$-qc conjugacy, where $K = \max(Q, (\kappa+1)/(\kappa-1))$ (see e.g., [L2, Cor. 10.5] and Lemma 10.4). □

COROLLARY 3.7. *Let $f \in \mathcal{C}$ and $g \in \mathcal{C}$ be two quadratic-like germs with $\mathrm{mod}(f) \geq \nu > 0$ and $\mathrm{mod}(g) \geq \nu > 0$. If these maps are hybrid equivalent then there exists a $K(\nu)$-qc hybrid conjugacy between them.*

3.5. *Hybrid classes.* Let $\mathcal{H}(f)$ stand for the hybrid class of $f \in \mathcal{QG}$. If $f \in \mathcal{C}$ then the hybrid class $\mathcal{H}(f)$ can be endowed with the following *Teichmüller-Sullivan metric*:

$$\mathrm{dist}_T(f, g) = \inf \log \mathrm{Dil}(h),$$

where $h$ runs over all hybrid conjugacies between $f$ and $g$ (see [S1]).

LEMMA 3.8. *For any germ $f_0 \in \mathcal{C}$ with connected Julia set and any circle map $g \in \mathcal{E}$, there exists a unique (up to affine conjugacy) germ $f \in \mathcal{H}(f_0)$ whose external map is equal to $g$. Moreover, if $\mathrm{mod}(f_0) \geq \mu > 0$ and $\mathrm{mod}(g) \geq \mu > 0$ then $\mathrm{dist}_T(f_0, f) \leq K(\mu)$.*

*Proof.* Let $g_0 : V_0 \to V_0'$ be the external map of $f_0$. Any two expanding circle maps, in particular $g : V \to V'$ and $g_0$, are quasi-conformally conjugate. Indeed, let $A = V' \setminus (V \cup \mathbb{D})$ and $A_0 = V_0' \setminus (V_0 \cup \mathbb{D})$ be outer fundamental annuli of $g$ and $g_0$ respectively. Their boundaries can be selected as smooth quasi-circles with dilatation controlled by $\mu$. Then there exists a $K = K(\mu)$-qc diffeomorphism $\psi : (\mathbb{C} \setminus V, A) \to (\mathbb{C} \setminus V_0, A_0)$ conjugating $g : \partial V \to \partial V'$ to $g_0 : \partial V_0 \to \partial V_0'$. It admits a unique extension to a $K$-qc map $\psi : \mathbb{C} \setminus \mathbb{D} \to \mathbb{C} \setminus \mathbb{D}$ conjugating $g : V \setminus \mathbb{D} \to V' \setminus \mathbb{D}$ to $g_0 : V_0 \setminus \mathbb{D} \to V_0' \setminus \mathbb{D}$. By the Reflection Principle, $\psi$ admits an extension to a $\mathbb{T}$-symmetric $K$-qc map $(\mathbb{C}, V', V) \to (\mathbb{C}, V_0', V_0)$ conjugating $g$ and $g_0$ on their domains.

Let us consider a $\psi$-push-forward of the standard conformal structure $\sigma$ from $\mathbb{C} \setminus \mathbb{D}$ to $\mathbb{C} \setminus \mathbb{D}$, $\nu = \psi_* \sigma$. It is preserved by $g_0 : V_0 \to V_0'$. Recall now that in the case of a connected Julia set, $g_0 : V_0 \setminus \mathbb{D} \to V_0' \setminus \mathbb{D}$ is conformally conjugate to $f_0 : U_0 \setminus K(f) \to U_0' \setminus K(f)$ by means of the Riemann mapping $\phi : \mathbb{C} \setminus K(f) \to \mathbb{C} \setminus \mathbb{D}$. Hence the structure $\mu = \phi^* \nu$ on $\mathbb{C} \setminus K(f_0)$ is preserved by $f_0$ near the Julia set. Let us extend it to $K(f)$ in a standard way: $\mu | K(f) = \sigma$. Straightening $\mu$ by the Measurable Riemann Mapping Theorem, we obtain a desired map $f$.

Let us now have two maps like this, $f$ and $\tilde{f}$. Since they are hybrid equivalent, there is a qc conjugacy $h : (U, K(f)) \to (\tilde{U}, K(\tilde{f}))$ near the filled Julia sets, such that $\bar{\partial} h = 0$ a.e. on $K(f)$. On the other hand, $f$ and $\tilde{f}$ are externally equivalent, so that there is a conformal map $H : \mathbb{C} \setminus K(f) \to \mathbb{C} \setminus K(\tilde{f})$ conjugating $f$ and $\tilde{f}$ near the Julia sets. These two maps match on the Julia sets (see [DH2, Prop. 6]), and hence glue together into a single conformal, and thus affine, map. □

THEOREM 3.9 (Straightening [DH2]). *If $f$ is a quadratic-like germ with connected Julia set then its hybrid class $\mathcal{H}(f)$ contains a unique quadratic polynomial $P : z \mapsto z^2 + \chi(f)$, where $c = \chi(f)$ is a point of the Mandelbrot set $M_0$. Moreover, $\mathrm{dist}_T(f, P) \leq K(\mathrm{mod}(f))$.*

*Proof.* By the previous lemma, there is a unique map $P \in \mathcal{QG}$ which is hybrid equivalent to $f$ and externally equivalent to $P_0 : z \mapsto z^2$. By Lemma 3.3, $P$ is the unique quadratic polynomial in $\mathcal{H}(f)$. □

Let us summarize Lemma 3.8 and Theorem 3.9:

THEOREM 3.10 (Internal mating). *Any parameter $c \in M$ can be mated with any circle map $g \in \mathcal{E}$ to obtain a unique (up to affine conjugacy) germ $f \equiv i_c(g) \in \mathcal{QG}$ such that $\chi(f) = c$ and $\pi(f) = g$.*

The hybrid class passing through a point $c \in M_0$ will also be denoted as $\mathcal{H}_c$. Thus we have a partition of the connectedness locus $\mathcal{C}$ into the hybrid classes labeled by the points of the Mandelbrot set and parametrized by the space $\mathcal{E}$ of circle maps.

Note that all quadratic-like germs with disconnected Julia set are hybrid equivalent, so that $\mathcal{QG} \setminus \mathcal{C}$ is a single hybrid class.

Let us finish with the following important remark: Any two germs $f_0$ and $f_1$ in the same hybrid class $\mathcal{H}$ can be included in a certain complex one-parameter family of maps called the *Beltrami disk*. Let $h$ be a hybrid conjugacy between $f_0$ and $f_1$, and $\mu = \bar{\partial} h / \partial h$ be the corresponding Beltrami differential. Let us consider a complex one-parameter family of Beltrami differentials $\mu_\lambda = \lambda \mu$, $\lambda \in \mathbb{D}_{1+\varepsilon}$. Let $h_\lambda$ be the solution of the corresponding Beltrami equation. Then by definition, the family $\{f_\lambda = h_\lambda \circ f_0 \circ h_\lambda^{-1}\}$ is the Beltrami disk via $f_0$ and $f_1$.

The real one-parameter family $\{f_\lambda : |\lambda| < 1 + \varepsilon\}$, is called the *Beltrami path* joining $f_0$ and $f_1$.

## 4. The space of quadratic-like germs

4.1. *Topology and analytic structure.* McMullen has supplied the space $\mathcal{QM}$ of quadratic-like maps with the Carathéodory *convergence structure* by declaring that a sequence $f_n : U_n \to U'_n$ converges to $f : U \to U'$ if the pointed domains $(U_n, 0)$ Carathéodory converge to $(U, 0)$, and $f_n \to f$ uniformly on compact subsets of $U$ (see [McM1, §4]). This structure can be pushed down to the space $\mathcal{QG}$ of quadratic-like germs by declaring $[f_n] \to [f]$ if the sequence $[f_n]$ can be split into finitely many subsequences $[f_m^i]$ which admit representatives $f_m^i$ Carathéodory converging to certain representatives $f_i$ of $[f]$ (the splitting of the sequence is actually not needed in the case when $f \in \mathcal{C}$ ). Below we will show that this convergence structure on $\mathcal{QG}$ is consistent with a certain topology, which in turn can be refined to a natural complex analytic structure modeled on a family of Banach spaces. For the background for this section the reader should consult Appendix 2.

As in subsection 11.3,

- $\mathbb{V}$ will stand for the directed set of topological discs $V \ni 0$ with piecewise smooth boundary, and $U \succ V$ if $U \Subset V$;

- $\mathcal{B}_V$ will denote the space of normalized analytic functions $f(z) = c + z^2 + \ldots$ on $V \in \mathbb{V}$ continuous up to the boundary supplied with sup-distance;

- $\mathcal{B} = \lim \mathcal{B}_V$ will stand for the space of normalized analytic germs at 0.

Let us now supply the space $\mathcal{QG}$ of normalized quadratic-like germs with topology and complex analytic structure modeled on a family of Banach spaces $\mathcal{B}_V$. Given $f \in \mathcal{QG}$, let $\mathbb{V}_f$ stand for the family of topological disks with piecewise smooth boundary such that $f : V \to fV$ is quadratic-like. If $g \in \mathcal{B}_V(f, \varepsilon)$ is sufficiently close to $f$ in the Banach space $\mathcal{B}_V$ then it is quadratic-like on a domain $U$ slightly smaller than $V$. Hence $g$ represents a quadratic-like germ. Thus we have an injection $j_{f,V} : \mathcal{B}_V(f, \varepsilon) \to \mathcal{QG}$. This family of injections obviously satisfies properties P1–P3 stated in Appendix 2 (with linear transition maps), and hence endows $\mathcal{QG}$ with topology and complex analytic structure.

Note that by Lemma 11.6, convergence in this topology coincides with the quotient of the Carathéodory convergence.

If $f \in \mathcal{C}$ then the domains $V$ on which $f$ is quadratic-like form a directed set $\mathbb{V}_f$, so that $f$ is a regular point of $\mathcal{QG}$ (see Appendix 2). Since the transition maps are linear, the tangent space $\mathrm{T}_f \mathcal{QG}$ is naturally identified with the inductive limit

(4.1) $$\mathcal{B}(f) = \lim_{V \in \mathbb{V}_f} \mathcal{B}_V,$$

which is in turn identified with the space of germs of analytic vector fields near the filled Julia set $K(f)$ normalized at the origin as $v(z) = \delta + az^3 + \ldots$.

Since any finitely dimensional submanifold locally sits in some space $\mathcal{B}_V$, it can be locally identified with an analytic finitely parameter family of functions $f_\lambda(z)$ on $V$ (so that $f_\lambda(z)$ is analytic in two variables).

Given $\mu > 0$, $R > 0$, let $\mathcal{QG}(\mu, R)$ stand for the space of quadratic-like germs which have normalized representatives $f : V \to V'$, $z \mapsto c + z^2 + \ldots$ with $\operatorname{mod}(f) \geq \mu$, $|c| \leq R$, and $\operatorname{dist}_{\mathrm{hyp}}(0, c) \leq R$, where the hyperbolic distance is measured in $V'$. Let

$$\mathcal{QG}(\mu) = \{f \in \mathcal{QG} : \operatorname{mod}(f) \geq \mu\}.$$

Note that

(4.2) $$\mathcal{C}(\mu) \equiv \mathcal{QG}(\mu) \cap \mathcal{C} \subset \mathcal{QG}(\mu, R(\mu))$$

(indeed, for a quadratic map $z \mapsto z^2 + c$ with connected Julia set, we have: $c \leq 2$, and the statement for quadratic-like maps follows from the Straightening Theorem). Similarly, let $\mathcal{E}(\mu) = \{g \in \mathcal{E} : \operatorname{mod}(g) \geq \mu\}$.

LEMMA 4.1 (Compactness).  *A subset $\mathcal{K}$ of $\mathcal{QG}$ (respectively: of $\mathcal{C}$ or $\mathcal{E}$) is pre-compact if and only if it is contained in some $\mathcal{QG}(\mu, R)$ (respectively: in $\mathcal{C}(\mu)$ or $\mathcal{E}(\mu)$). Any compact set $\mathcal{K}$ sits in a union of finitely many Banach slices and bears a Montel metric $\operatorname{dist}_M$ (see Appendix 2) well-defined up to Hölder equivalence.*

*Proof.* Sequential pre-compactness of $\mathcal{QG}(\mu, R)$ follows from [McM1, Th. 5.6]), and by Lemma 11.6 it implies pre-compactness.

Vice versa, let $\mathcal{K} \subset \mathcal{QG}$ be pre-compact, thus sequentially pre-compact. Since $c = f(0)$ continuously depends on $f$, it is bounded on $\mathcal{K}$. Moreover, $\operatorname{mod}(f)$ and $\operatorname{dist}_{\mathrm{hyp}}(0, f(0))$ are sequentially continuous on the level of maps: If $f_n \to f$ then

$$\lim \operatorname{mod}(f_n) = \operatorname{mod}(f) \quad \text{and} \quad \lim \operatorname{dist}_{\mathrm{hyp}}(0, f_n(0)) = \operatorname{dist}_{\mathrm{hyp}}(0, f(0)).$$

Since convergence of germs is described in terms of representatives,

$$\operatorname{mod}(f) \geq \mu > 0 \text{ and } \operatorname{dist}_{\mathrm{hyp}}(0, f(0)) \leq R < \infty \quad \text{for} \quad f \in \mathcal{K}.$$

The criterion for $\mathcal{C}$ follows by (4.2), and the criterion for $\mathcal{E}$ is completely analogous.

The last statement follows from Lemma 11.5. □

Let us consider a holomorphic family $f_\lambda$, $\lambda \in \Lambda$, of quadratic-like germs over a Banach domain $(\Lambda, 0)$, $f_0 \equiv f_{\lambda_0}$. Then this family locally sits in some Banach slice $\mathcal{B}_U$ and is represented there by a holomorphic family of quadratic-like maps $f_\lambda : V_\lambda \to V'_\lambda$. Take a thickened fundamental annulus $\tilde{A}_0$ (i.e., a little neighborhood of the fundamental annulus $A_0$) of $f_0$ with piecewise smooth boundary compactly sitting in $V_0'$. The following useful statement shows that this annulus can be included in a holomorphic family:

LEMMA 4.2. *There is a smooth holomorphic motion $h_\lambda$ of $\tilde{A}_0$ over a neighborhood $\Lambda' \subset \Lambda$ of 0, such that $\tilde{A}_\lambda \equiv h_\lambda \tilde{A}_0$ is a thickened fundamental annulus of $f_\lambda$, which respects the dynamical relation near the boundary of $A_0$.*

*Proof.* Fix a collar neighborhood of the outer boundary of $A_0$, and let the corresponding collar neighborhood of the inner boundary move as prescribed by dynamics. By Lemma 11.2, this motion can be smoothly interpolated through the whole annulus. □

Douady and Hubbard [DH, Prop. 9] formulated this statement (for one-parameter families) as existence of horizontally analytic smooth *tubing*, i.e., a smooth map $\Psi_\lambda(z)$,

(4.3) $$\Psi : \Lambda' \times \mathbb{A}(2 - \varepsilon, 4 + \varepsilon) \to \bigcup_{\lambda \in \Lambda'} \tilde{A}_\lambda,$$

analytic in $\lambda \in \Lambda'$ for any given $z \in \mathbb{A}(2 - \varepsilon, 4 + \varepsilon)$, fibered over $\mathrm{id}|\Lambda'$ and such that $\Psi_\lambda$ conjugates $P_0 : z \mapsto z^2$ near $\mathbb{T}_2$ to $f_\lambda$ near the inner boundary of $A_\lambda$.

4.2. *Complex structure on the space of circle maps.* In a similar way we can supply the space $\mathcal{E}$ of expanding circle maps with the inductive limit topology and real analytic structure. Namely, let us represent $\mathbb{T}$ as $\mathbb{R}/(\gamma : x \mapsto x + 1)$ so that $1 \in \mathbb{T}$ corresponds to $0 \in \mathbb{R}$. Let $V$ be a $\gamma$-invariant $\mathbb{R}$-symmetric neighborhood of $\mathbb{R}$, and let $\mathcal{B}_V$ stand for the Banach space of functions $f$ analytic on $V$, real on $\mathbb{R}$, normalized as $f(0) = 0$, and satisfying the following equation: $f(z + 1) = f(z) + 2$. (In other words, this is the space of degree-two circle maps analytic in a given neighborhood of $\mathbb{T}$ and fixing 1.)

Let $\mathcal{E}_V$ be the set of expanding circle maps $f \in \mathcal{E}$ which belong to $\mathcal{B}_V$. It is clearly an open subset of $\mathcal{B}_V$. Thus we have a natural representation of $\mathcal{E}$ as the inductive limit of real Banach manifolds $\mathcal{E}_V$.

It is not obvious that $\mathcal{E}$ can also be endowed with complex analytic structure. To see this let us consider the hybrid class of $z \mapsto z^2$, $\mathcal{H}^0 \equiv \mathcal{H}(z^2) = \{f \in \mathcal{QG} : f(0) = 0\}$. Since the condition $f(0) = 0$ specifies in every Banach space $\mathcal{B}_V$ a codimension 1 linear subspace, $\mathcal{H}^0$ is naturally endowed with topology and complex analytic structure.

LEMMA 4.3. *The space $\mathcal{E}$ of circle maps and the hybrid class $\mathcal{H}_0$ are homeomorphic.*

*Proof.* The homeomorphism $i \equiv i_0 : \mathcal{E} \to \mathcal{H}_0$ is constructed as the mating of $c = 0 \in M$ with $g \in \mathcal{E}$. By the Mating Theorem 3.10, $i$ is one-to-one.

Take a map $g_0 \in \mathcal{E}_V$ with a fundamental annulus $A_0$. Then for all nearby $g \in \mathcal{E}_V$, we can select a continuously moving fundamental annulus $A = A_g$ (as in Lemma 4.2). It follows that the map $\psi_g : \mathbb{A}(r, r^2) \to A_g$ conjugating $P_0 : \mathbb{T}_r \to \mathbb{T}_{r^2}$ to the boundary restriction of $g$ can be selected continuously in $g$.

Let us consider the conformal structure $\psi_g^* \sigma$ on $\mathbb{A}(r, r^2)$. Pull it back by the iterates of $P_0$ and extend it as $\sigma$ beyond $\mathbb{A}(1, r^2)$. We obtain a continuous family of conformal structures $\nu_g$ on $\mathbb{C}$. By the Measurable Riemann Mapping Theorem, the straightening map

$$h_g : (\mathbb{C}, 0) \to (\mathbb{C}, 0), \quad (h_g)_* \nu_g = \sigma, \quad (h_g)'(0) = 1,$$

depends continuously on $g$. Hence $f = i(g) = h_g \circ P_0 \circ h_g^{-1}$ depends continuously on $g$ as well.

Vice versa, it is easy to see that the filled Julia set $K(f)$ depends continuously on $f \in \mathcal{H}_0$. Hence the normalized Riemann mapping $\phi \equiv \phi_f : \mathbb{C} \setminus K(f) \to \mathbb{C} \setminus \bar{\mathbb{D}}$ depends continuously on $f$ in compact-open topology. Let $\gamma \subset \mathbb{C} \setminus \bar{\mathbb{D}}$ be a closed curve homotopic to $\mathbb{T}$ which belongs to the domain $V_0$ of $g_0 \in \mathcal{E}$. Then for $f$ near $f_0$, the external map $g = \phi \circ f \circ \phi^{-1}$ restricted on $\gamma$ depends continuously on $f$. By the reflection and maximum principles, $g$ depends continuously on $f$ on the whole annulus enclosed by $\gamma$ and the symmetric curve. Thus the map $i^{-1}$ is continuous. $\square$

*Remark.* In the above proof we implicitly use the following fact. Let $f_\lambda = h_\lambda \circ f_0 \circ h_\lambda^{-1}$, $\lambda \in D \in \mathbb{C}$, where $h_\lambda$ is a holomorphic family of qc maps and $f_\lambda$ is a family of holomorphic maps. Then $f_\lambda$ is holomorphic in $\lambda$ (it is not trivial because $h_\lambda^{-1}$ need not depend holomorphically on $\lambda$). Indeed, taking $\partial/\partial \bar{\lambda}$ of the expression $h_\lambda \circ f_0 = f_\lambda \circ h_\lambda$, we obtain:

$$0 = \frac{\partial h_\lambda}{\partial \bar{\lambda}} \circ f_0 = f'_\lambda \circ \frac{\partial h_\lambda}{\partial \bar{\lambda}} + \frac{\partial f_\lambda}{\partial \bar{\lambda}} \circ h_\lambda = \frac{\partial f_\lambda}{\partial \bar{\lambda}} \circ h_\lambda.$$

Now the natural complex analytic structure on $\mathcal{H}_0$ can be transferred to $\mathcal{E}$ via the above homeomorphism. (This construction is inspired by the construction of the complex structure on the Teichmüller spaces via the Bers embedding [Be].)

Let

(4.4) $\quad \Pi = i_0 \circ \Pi : \mathcal{QG} \to \mathcal{H}_0, \quad I_c = i_c \circ \pi : \mathcal{H}_0 \to \mathcal{H}_c.$

4.3. *Analyticity of $\pi$ and $i_c$.*

LEMMA 4.4.

- *The projection $\pi : \mathcal{C} \to \mathcal{E}$ is proper;*

- *The projection $\pi : \mathcal{QG} \to \mathcal{E}$ is complex analytic.*

*Proof.* • By definition, "proper" means that preimages of compact sets are compact. Let $\mathcal{K} \subset \mathcal{E}$ be compact. Then by Lemma 4.1, $\mod(g) \geq \mu > 0$, $g \in \mathcal{K}$.

Let $f \in \mathcal{C} \cap \pi^{-1}\mathcal{K}$. Then by Lemma 3.2, $\mod(f) \geq \mu$ as well, so that by Lemma 4.1 $\pi^{-1}\mathcal{K}$ is compact.

• Let $f_0 \in \mathcal{QG}_V$, $A_0$ be its fundamental annulus. For $f \in \mathcal{QG}_V$ near $f_0$, select a holomorphically moving fundamental annulus $A_f$ (see Lemma 4.2). This defines a holomorphic family $\mu_f$ of conformal structures on $A_0$. Pulling these structures back to the Riemann surface $T_0 \equiv T(f_0)$ constructed above, see (3.2), we obtain a holomorphic family of conformal structures $\nu_f$ on $T_0$. Realize $T_0$ as $\mathbb{C} \setminus \bar{D}$, put the standard structure on $\mathbb{D}$, and solve the Beltrami equation: $(h_f)_*\nu_f = \sigma$. The analytic dependence in the Measurable Riemann Mapping Theorem ensures that the maps $\Pi(f) = h_f \circ f_0 \circ h_f^{-1} \in \mathcal{H}_0$ analytically depend on $f$. As the complex structure on $\mathcal{E}$ is by definition transferred from $\mathcal{H}_0$, we conclude that $g_f = \pi(f)$ analytically depends on $f$ as well. □

LEMMA 4.5. *For any $c \in M$, $f = i_c(g)$ depends analytically on $g \in \mathcal{E}$.*

*Proof.* The proof is similar to that for the previous lemma. Let $A_0$ be a fundamental annulus for $f_0 \in \mathcal{C}$. Select a fundamental annulus $B_g$ which moves holomorphically with $g$, and a family of diffeomorphisms $h_g : A_0 \to B_g$ respecting the dynamics on the boundaries and depending analytically on $g$. Then the conformal structure $\mu_g = (h_g)^*\sigma$ depends analytically on $g$. Pulling it back by iterates of $f_0$ and straightening, we complete the proof (using analytic dependence in the Measurable Riemann Mapping Theorem). □

Thus one can say that the mating $f = i_c(g)$ is *horizontally analytic*.

4.4. *Infinitesimal deformations.* Let us introduce spaces needed for the description of the tangent spaces to the hybrid classes. For $f \in \mathcal{C}$, let $\mathbb{B}(f)$ be the space of $f$-invariant Beltrami differentials $\mu$ of class $L^\infty$ near $K(f)$ such that $\mu = 0$ a.e. on $K(f)$. Consider the Beltrami path $h_t$ in the direction of $\mu \in \mathbb{B}(f)$, i.e., the family of normalized solutions of the the Beltrami equations $\bar\partial h_t / \partial h_t = t\mu$ for small $|t|$ (see §3.5). The velocity of this path at $f$,

$$w = \frac{dh_t}{dt}\Big|_{t=0},$$

is a vector field near $K(f)$ which has locally square integrable distributional derivatives (i.e., of class $H$; see Appendix 1) and satisfies the $\bar\partial$-equation $\bar\partial w = \mu$. Let $\mathbb{F}(f)$ stand for the space of such vector fields (corresponding to all possible $\mu \in \mathbb{B}(f)$).

For $f \in \mathcal{C}$, let us consider the tangent space to the hybrid class of $f$,

(4.5) $$E^h(f) = \mathrm{T}_f\mathcal{H}(f).$$

By definition, it consists of the velocities at $f$ of all smooth curves $\gamma(t) \in \mathcal{H}(f)$ through $f$. The space $E^h(f)$ and its vectors will be called *horizontal*.

*Remarks.* 1. Vector fields $v(z)/dz \in \mathrm{T}_f \mathcal{QG}$ are normalized so that $v'(0) = 0$ and considered up to adding a vector field $az^2/dz$. By this we can make $v''(0) = 0$ but sometimes we will prefer a different normalization.

2. On the plane $\mathbb{C}$ we will freely identify functions $v(z)$ with the corresponding vector fields $v(z)/dz$.

LEMMA 4.6. *The horizontal space $E^h(f)$ consists of (normalized) holomorphic vector fields $v(z)/dz$ near $K(f)$ which admit a representation $v(z) = w(fz) - f'(z)w(z)$ with some $w/dz \in \mathbb{F}(f)$.*

*Proof.* Let us consider a smooth path $f_t \in \mathcal{H}(f)$ tangent to $v(z)/dz \in E^h(f)$. Then Lemma 4.2 and the Measurable Riemann Mapping Theorem imply that there is a smooth path of qc maps $h_t$ conjugating $f$ to $f_t$, $h_0 = \mathrm{id}$. Let $w/dz \in \mathbb{F}(f)$ be the velocity of this path at id. Linearizing the curve $t \mapsto h_t \circ f \circ h_t^{-1}$ at $t = 0$, we conclude that $\{f_t\}$ is tangent to $(w \circ f - f'w)/dz$ at $f$, so that $v$ admits a desired representation.

Vice versa, let $v = w \circ f - f'w$ with $w/dz \in \mathbb{F}(f)$. Then the Beltrami differential $\mu = \bar{\partial}w$ belongs to $\mathbb{B}(f)$. The corresponding Beltrami path $f_t = h_t \circ f \circ h_t^{-1}$ (i.e., $\bar{\partial}h_t/\partial h_t = tw$) is a smooth curve in $\mathcal{H}(f)$ tangent to $v$ at $f$. □

Let $f \in \mathcal{C}$. A vector field $v(z)/dz \in \mathrm{T}_f \mathcal{QG}$ is called *vertical* if there is a holomorphic vector field $\alpha(z)/dz$ on $\bar{\mathbb{C}} \setminus K(f)$ vanishing at $\infty$ such that

(4.6) $$v(z) = \alpha(fz) - f'(z)\alpha(z)$$

near the Julia set. (Note that the above condition is equivalent to saying that $\alpha$ is a holomorphic function on $\mathbb{C} \setminus K(f)$ with at most simple pole at $\infty$.) The space of vertical tangent vectors at $f$ will be denoted by $E^v(f)$. (We will eventually show that it is the tangent space to the vertical fiber $\mathcal{Z}_f$.)

LEMMA 4.7. *For $f \in \mathcal{C}$, $\mathrm{T}_f \mathcal{QG} = E^h(f) \oplus E^v(f)$.*

*Proof. Existence.* Let us consider a holomorphic vector field $v(z)/dz \in \mathrm{T}_f \mathcal{QG}$ in a neighborhood of $K(f)$. Select a quadratic-like representative $f : U \to U'$ such that $v$ is well-defined in a neighborhood of $\bar{U}'$. Then there exists a smooth vector field $w(z)/dz$ in a neighborhood of $\bar{\mathbb{C}} \setminus U$ vanishing near $\infty$ and such that

(4.7) $$v(z) = w(f(z)) - f'(z)\, w(z), \quad z \in U \cap \mathrm{Dom}(w).$$

By means of this equation $w(z)/dz$ can be extended to a smooth vector field in $\mathbb{C} \setminus K(f)$ satisfying (4.7) in $U$.

Let us consider the Beltrami differential $\mu = \bar{\partial} w$ in $\bar{\mathbb{C}} \setminus K(f)$ extended by 0 to the filled Julia set $K(f)$. Since $v$ is holomorphic on $U$, (4.7) implies that $\mu$ is $f$-invariant over there. Hence it has a bounded $L^\infty$-norm on the whole sphere (equal to its $L^\infty$-norm on $\bar{\mathbb{C}} \setminus U$).

Let us solve the $\bar{\partial}$-problem $\bar{\partial} u = \mu$, where $u(z)/dz$ is a vector field on $\mathbb{C}$ of class $H$ vanishing at $\infty$. Then:

a) The vector field $v^h = u \circ f - f'u$ on $U$ is holomorphic since $\mu$ is $f$-invariant. By adding a linear function $az + b$ to $u$ we can normalize $v^h$ so that $(v^h)'(0) = (v^h)''(0) = 0$. Since $\mu = 0$ on $K(f)$, $v^h$ is horizontal.

b) Let $\alpha = w - u$. The corresponding vector field $\alpha(z)/dz$ is holomorphic on $\bar{\mathbb{C}} \setminus K(f)$ (since $\bar{\partial}\alpha = 0$) and vanishes at $\infty$. Moreover, $v - v^h = \alpha \circ f - f'\alpha$ on $U$. Hence $(v - v^h)/dz$ is a vertical vector field.

*Uniqueness.* Assume that there exists a vector field $v/dz \in E^h(f) \cap E^v(f)$, $v \neq 0$. Then $v = w \circ f - f'w$ with $w/dz \in \mathbb{F}(f)$ and $v = \alpha \circ f - f'\alpha$ where $\alpha(z)/dz$ is holomorphic on $\bar{\mathbb{C}} \setminus K(f)$ and vanishes at $\infty$.

Let us consider a vector field $u/dz = (w - \alpha)/dz$ on $U' \setminus K(f)$. Since it is $f$-invariant, it is bounded with respect to the hyperbolic norm on $U$. Hence $|u(z)| \to 0$ as $z \to J(f)$, $z \in U$, so that $u$ admits a continuous extension to $\bar{U}$ vanishing on the Julia set $J(f)$.

Thus the vector fields $w(z)/dz$ and $\alpha(z)/dz$ match on the Julia set $J(f)$, i.e., the vector field $\beta(z)/dz$ which is equal to $w(z)/dz$ on $K(f)$ and equal to $\alpha(z)/dz$ on $\bar{\mathbb{C}} \setminus K(f)$ is continuous on the whole sphere. Moreover, this vector field has distributional derivatives of class $L^2$ and $\bar{\partial}\beta = 0$ (see e.g., Lemma 10.4). By Weil's lemma, $\beta(z)/dz$ is holomorphic on the whole sphere. Since it vanishes at $\infty$, $\beta(z) = az + b$ is linear.

Thus $v(z) = af(z) + b - f'(z)(az + b)$, where $f(z) = c + z^2 + \ldots$. If $v$ is normalized by $v'(0) = v''(0) = 0$, then a straightforward calculation yields $a = b = 0$. Hence $v = 0$ as well. □

4.5. *Horizontal foliation of $\mathcal{C}$.* In this section we will show that the hybrid classes $\mathcal{H}_c$, $c \in M$, are codimension-1 complex analytic submanifolds in $\mathcal{QG}$.

As usual, a foliation is called *analytic* (or smooth etc.) if it locally admits an analytic (smooth etc.) straightening.

LEMMA 4.8. *The partition of $\mathrm{int}\,\mathcal{C}$ into the hybrid classes is a complex analytic foliation, with an analytic straightening given by the mating.*

*Proof.* On the hyperbolic components of $\mathrm{int}\,\mathcal{C}$, the analytic map $f \mapsto (\pi(f), \lambda(f))$ straightens the leaves, where $\lambda(f)$ is the multiplier of the attracting periodic point. To obtain an analytic straightening over a "queer" component $U$ of $\mathrm{int}\,M_0$ (see §10.2), select a reference point $c_0 \in U$, and a nontrivial measurable line field $\gamma: z \mapsto e^{2i\theta(z)}$ on the Julia set $J(P_{c_0})$. Then the unit disk of Beltrami differentials $\lambda\gamma$, $\lambda \in \mathbb{D}$, analytically parametrizes $U$. Its product with $\mathcal{E}$ analytically parametrizes the corresponding component of int$\mathcal{C}$ (by analytic dependence on parameters in the Measurable Riemann Mapping Theorem). □

LEMMA 4.9. *For $f \in \text{int}\mathcal{C}$, the tangent space to the vertical fiber $\mathcal{Z}_f$ coincides with $E^v(f)$. Moreover, if $f : V \to V'$ is a quadratic-like representative of $f$ then the vertical vector fields extend to holomorphic fields on $V$.*

*Proof.* For $g \in \mathcal{Z}_f$ near $f$, let us consider a normalized qc map $h_g : \bar{\mathbb{C}} \to \bar{\mathbb{C}}$ conjugating $f$ to $g$ and conformal on $\bar{\mathbb{C}} \setminus K(f)$. By the Measurable Riemann Mapping Theorem, it depends holomorphically on $g$. Hence the map $(g, z) \mapsto h_g(z)$ is holomorphic in two variables near $\{f\} \times (\bar{\mathbb{C}} \setminus K(f))$.

Let us consider a smooth curve $g(t)$ in the vertical fiber $\mathcal{Z}_f$ tangent to a vector field $v(z)/dz$ at $f$ (for $t = 0$). Since $h_{g(t)}$ smoothly depends on $t$, we can consider the vector field $\alpha(z)/dz$ tangent to this curve at $t = 0$:

$$\alpha = \frac{dh_{g(t)}}{dt}\Big|_{t=0}. \tag{4.8}$$

Then $v = \alpha \circ f - f'\alpha$ on the domain $V$ of $f$. But by the above discussion, $\alpha$ is holomorphic on $\bar{\mathbb{C}} \setminus K(f)$. Hence $v \in E^v(f)$.

Thus $\mathrm{T}_f \mathcal{Z}_f \subset E^v(f)$. Since both spaces are one-dimensional, they coincide.

The second statement of the lemma also follows from the above discussion. □

*Remark.* Formula (4.8) is valid for all $f \in \mathcal{C}$ (not only for $f \in \text{int}\mathcal{C}$). In the general formula, $h_g$ should be understood as $\phi_g \circ \phi_f^{-1}$, where the $\phi_g^{-1}$ are the uniformizations at $\infty$ (see §3.2). Justification for this formula requires proof that the vertical fibers are submanifolds (see §4.12) and that $\phi_g$ smoothly depends on $g$.

We will also need the following technical lemma:

LEMMA 4.10. *Let $f : V' \to V''$ be a quadratic-like map with connected Julia set, and let $V = f^{-1}V'$. Consider vector fields $v(z)/dz$ and $\alpha(z)/dz$ satisfying (4.6) on $V'$, where $v$ is holomorphic in $V''$ and $\alpha$ is holomorphic in $\bar{\mathbb{C}} \setminus K(f)$. Then*

$$\|\alpha\|_{\bar{\mathbb{C}}\setminus V} \leq C\|v\|_{V'} \quad \text{and} \quad \|v\|_V \geq C^{-1}\|v\|_{V'}$$

*with a constant $C$ depending only on $\mu = \text{mod}(V' \setminus V)$.*

*Proof.* Let $\gamma = \partial V$. By the standard normality or hyperbolic metric arguments, the inverse branches of $f^{-N}|\gamma$ are uniformly contracting on $\gamma$ with some constant $\lambda^{-1} < 1$ (where $N$ and $\lambda$ depend on $\mu$ only).

Let $\|v\|_V = \varepsilon$ and $\|\alpha\|_{\bar{\mathbb{C}} \setminus V} = M$. Incorporating the above contracting property and the Maximum Principle into (4.6) we obtain:

$$M \leq \|\alpha\|_{f^{-N}\gamma} \leq (M + A\varepsilon)/\lambda,$$

where $A$ depends only on $\mu$. Hence $M \leq A\varepsilon/(\lambda - 1)$, which proves the first desired estimate. Using (4.6) once again we conclude that $\|v\|_{V'} \leq A\varepsilon(K + 1)/(\lambda - 1)$, where $K = \|f'\|_{V'}$, which proves the second estimate. □

THEOREM 4.11 (Leaves). *Hybrid classes $\mathcal{H}_c$, $c \in M$, are connected codimension-one complex analytic submanifolds of $\mathcal{QG}$. The quadratic family $\mathcal{Q} \equiv \{P_c(z) = z^2 + c\}$ is a transversal to these submanifolds.*

*Proof.* We have: $\pi \circ i_c = \text{id}$. Hence by the definitions (see Appendix 2) and the last two lemmas, $\pi : \mathcal{QG} \to \mathcal{E}$ is a complex analytic submersion and $i_c : \mathcal{E} \to \mathcal{QG}$ is a complex analytic embedding, so that its image $\mathcal{H}_c$ is a complex analytic submanifold in $\mathcal{QG}$. Since any two points in a hybrid class can be joined by a Beltrami path, the hybrid classes are connected. The quadratic family is transverse to all the leaves since by Lemma 3.3, it is a fiber of $\pi$.

Let us now show that $\text{codim}\,\mathcal{H}_c = 1$ at any point $f \in \mathcal{H}_c$ (note that by the last statement, it is true at $f = P_c$). By Lemma 4.8 it is true for $c \in \text{int}\,M_0$. Let $c \in \partial M$. Select a sequence of maps $f_n \in \text{int}\,\mathcal{C}$ converging to $f$. Let us consider the tangent projection $P = \text{D}(i_c \circ \pi)(f) : \text{T}_f\mathcal{QG} \to E^h(f)$ parallel to the vertical space $E^v(f)$, and the analogous projections $P_n$ at $f_n$. By Lemma 4.9, the latter projections have corank 1.

Let us first show that $\text{corank}\,P \leq 1$. Otherwise there would be a two-dimensional tangent space $F \subset E^v(f)$ sitting in some Banach slice $\mathcal{B}_V \subset \text{T}_f\mathcal{QG}$. For $n$ sufficiently big, this slice is naturally contained in the tangent spaces $\text{T}_{f_n}\mathcal{QG}$ as well. Let us consider the slices $L_n = \text{Im}\,P_n \cap \mathcal{B}_V$ of the horizontal spaces. By Lemma 11.3, they have codimension 1 in $\mathcal{B}_V$. Hence there exist vectors $v_n \in F \cap L_n$, $\|v_n\|_V = 1$.

Since $P_n \to P$, there exists a $U \succ V$ such that eventually $P_n(\mathcal{B}_V) \subset \mathcal{B}_U$ and the Banach operators $P_n : \mathcal{B}_V \to \mathcal{B}_U$ converge to $P : \mathcal{B}_V \to \mathcal{B}_U$. But since $F$ is finite-dimensional, $\|v_n\|_U \geq c > 0$. Thus $\|P_n v_n\|_U \geq c$ while $P_0 v_n = 0$ contradicting the operator convergence.

Let us finally prove the opposite inequality: $\text{corank}\,P \geq 1$. To this end let us consider vertical vector fields $v_n$ at $f_n$. Let us select quadratic-like representatives $f_n : V_n \to V'_n \to V''_n$ Carathéodory converging to $f : V \to V' \to V''$. By Lemma 4.9, the $v_n$ holomorphically extend to $V''_n$. Normalize the vector fields so that $\|v_n\|_{V'_n} = 1$. Then Lemma 4.10 implies that

(4.9) $$\|v_n\|_V \geq c > 0.$$

Select a subsequence of these vector fields converging uniformly on compact subsets of $V'$ to a vector field $v$ on $V'$. By (4.9), $v \neq 0$. Since $P_n \to P$, we conclude that $Pv = 0$, and we are done. □

*Remark.* Under more usual circumstances the codimension-1 statement would immediately follow from the facts that $\mathcal{H}_c$ is connected and has codimension one at the point $P_c$. However, a justification of this argument in our setting would be more involved than the above proof.

One can extract extra useful information from the proof above:

LEMMA 4.12. *For $f \in \mathcal{C}$,*

$$E^v(f) = \mathrm{Ker} D\pi(f). \tag{4.10}$$

*Proof.* In the proof of Theorem 4.11 we constructed a sequence of vector fields $v_n$ and $\alpha_n$ satisfying (4.6) and passed to a limit $v = \lim v_n \in \mathrm{Ker} D\pi(f)$ (along a subsequence). By Lemma 4.10 we can also pass to a limit $\alpha$ on $\bar{\mathbb{C}} \setminus \bar{V}$ along a subsequence of the $\alpha_n$'s. Then $v$ and $\alpha$ are related by (4.6) on $V' \setminus \bar{V}$. By means of this equation we can now extend $\alpha$ to $\bar{\mathbb{C}} \setminus K(f)$. Hence $v \in E^v(f)$. Since $\dim \mathrm{Ker} D\pi(f) = 1$ (by Theorem 4.11), $\mathrm{Ker} D\pi(f) \subset E^v(f)$.

Furthermore, by Lemma 4.7, $E^v(f)$ complements $E^h(f)$, and by Theorem 4.11, the latter space has codimension 1. Hence $\dim E^v(f) = 1$, and the conclusion follows. □

Denote by $\mathcal{F}$ the foliation of $\mathcal{C}$ into the hybrid classes. Accordingly the hybrid classes in the connectedness locus will also be called the *leaves* of $\mathcal{F}$.

Let us summarize the above information:

THEOREM 4.13 (Product structure). *The connectedness locus $\mathcal{C}$ is homeomorphic to the product $\mathcal{E} \times M_0$. The homeomorphism is provided by mating $f = i_c(g)$. It is horizontally analytic everywhere, and analytic in both variables for $c \in \mathrm{int} M_0$.*

*Proof.* By the Mating Theorem, the mating $f = i_c(g)$ provides one-to-one correspondence between $\mathcal{E} \times M$ and $\mathcal{C}$. Moreover, the inverse map $i^{-1} : f \mapsto (\pi(f), \chi(f))$ is continuous, as $\pi$ is continuous by Lemma 4.4 and $\chi$ is continuous by [DH2], [McM2, Prop. 4.7]. Since by Lemma 4.4, $i^{-1}$ is proper, $i$ is sequentially continuous. By Lemma 11.6, it is continuous.

The last statement is the content of Lemmas 4.5 and 4.8. □

4.6. *Regularity of the external mating.* Let us now study analytic properties of the external mating defined in subsection 3.3.

LEMMA 4.14. *The external mating $\theta : \mathcal{E}^m \to \mathcal{QG} \setminus \mathcal{C}$ is a diffeomorphism fibered over $\mathcal{E}$. Moreover, it is vertically conformal, i.e., it is conformal on the fibers $\mathcal{S}_g \to \mathcal{Z}_g \setminus \mathcal{C}$.*

*Proof.* Let us consider a pair $\lambda = (G, b)$ where $G : W \to W' \in \mathcal{H}_0$, $b \in W'$, varying near some $\lambda_0 = (G_0, b_0)$. Then by Lemma 4.2 and Lemma 11.2, there is a local smooth holomorphic motion $H_\lambda : (B_0, b_0) \to (B_G, b)$, where $B_G$ is a fundamental annulus of $G$. Consider the corresponding holomorphic family $\mu_\lambda = H_\lambda^* \sigma$ of conformal structures on $B_0$.

Let $R = R_G : \mathbb{C} \setminus K(G) \to \mathbb{C} \setminus \bar{\mathbb{D}}$ stand for the normalized Riemann mapping. Let us consider the corresponding marked external map

$$(g_\lambda : V_\lambda \to V'_\lambda, a_\lambda) = (\pi(G), R_G(b)) \in \mathcal{E}^m.$$

Transfer the structures $\mu_\lambda$ to this external model: $\nu_\lambda = (R_0)_* \mu_\lambda$ on $A_0 = V'_0 \setminus V_0$ (with 0 we label the objects corresponding to $\lambda_0$). Pulling $\nu_\lambda$ back by iterates of $g_0$ we obtain a holomorphic family of conformal structures on $V'_0 \setminus \bar{\mathbb{D}}$ which we will still denote by $\nu_\lambda$. Let $\hat{\nu}_\lambda$ be the symmetrization of $\nu_\lambda$ with respect to the unit circle. Then $\hat{\nu}_\lambda$ depends smoothly on $\lambda$.

Let us solve the Beltrami equation, $(h_\lambda)_* : \hat{\nu}_\lambda \to \sigma$, $h(0, 1, \infty) = (0, 1, \infty)$. The solution depends smoothly on $\lambda$ and hence $a_\lambda = h_\lambda(a_0)$ also depends smoothly on $\lambda$. Moreover, $g_\lambda$ depends on $\lambda$ holomorphically by definition of the complex structure on $\mathcal{E}$.

Reversing this construction we see that $(G, b)$ depends smoothly on $(g, a)$ as well.

Moreover, if $G \equiv G_0$ does not vary, then $g_\lambda \equiv g_0$ does not vary either, and hence $h_\lambda$ commutes with $g_0$. It follows that the identical map in the interior of the unit disk glues with $h_\lambda$ outside to a holomorphic motion (see e.g., [L2, Lemma 10.3]). Hence in this case $a_\lambda$ depends holomorphically on $\lambda$ (i.e., on $b$ only as $G$ stays fixed).

In a similar way one sees that $f_\lambda = \theta(g, a) \in \mathcal{QG}$ depends holomorphically on $\lambda$ (compare Lemma 4.5). Hence it depends smoothly on $(g, a)$ and, moreover, depends holomorphically on $a$ once $g$ is fixed. □

4.7. *Uniformization at $\infty$.* We need for further reference a statement on continuous dependence of the uniformization at $\infty$ on a quadratic-like map. Let $f : U \to U'$ be a quadratic-like map, $g = \pi(f) : V \to V'$ be its external, and $\phi_f : \mathbb{C} \setminus U \to \mathbb{C} \setminus V$ be the conformal map respecting the boundary dynamics (see §3.2). Let us normalize $\phi_f$ so that it has a positive derivative at $\infty$. If $N = N(f)$ is the maximal natural number such that $f^{-N} U \ni 0$ then $\phi_f$ admits an analytic extension

$$\phi_f : \Omega(f) \equiv \mathbb{C} \setminus f^{-N} U \to \mathbb{C} \setminus g^{-N} V \equiv \Delta(f).$$

In the case of a connected Julia set, $\Omega(f) = \mathbb{C} \setminus K(f)$ and $\Delta(f) = \mathbb{C} \setminus \bar{\mathbb{D}}$.

LEMMA 4.15. *Let a sequence of quadratic-like maps $f_n : U_n \to U'_n$ Carathéodory converge to a map $f : U \to U'$ with connected Julia set. Then the maps $\phi_{f_n}^{-1}$ converge to $\phi_f^{-1}$ uniformly on compact sets $K \subset \mathbb{C} \setminus \bar{\mathbb{D}}$.*

*Proof.* Let $g_n : V_n \to V'_n$ be the external maps of $f_n$, and $g : V \to V'$ be the external map of $f$. These maps can be selected in such a way that

- The $g_n$ Carathéodory converge to $g$ (by Lemma 4.4);

- The $\partial V_n$ are smooth quasi-circles with bounded dilatation (see [McM2, Prop. 4.10];

- The curves $\partial V_n$ uniformly converge to $\partial V$.

Hence for any $r > 1$ there exists an $m \in \mathbb{N}$ such that eventually $g_n^{-m} V \subset \mathbb{D}_r$. Moreover, since the critical point of $f$ is nonescaping, $N(f_n) \to \infty$. Hence $\Delta(f_n) \supset \mathbb{C} \setminus g_n^{-m} V$ for $n$ big enough. Altogether it follows that for any $r > 1$, $\Delta(f_n) \supset \mathbb{C} \setminus \mathbb{D}_r$ for $n$ big enough. In other words, all functions $\phi_n^{-1}$ are eventually defined on any $\mathbb{C} \setminus \mathbb{D}_r$, $r > 1$.

Moreover, these functions form a normal family on each domain $\mathbb{C} \setminus \mathbb{D}_r$. Indeed, the $\beta$-fixed points $\beta_n$ of the $f_n$ converge to the $\beta$-fixed point of $f$, and hence stay away from 0 and $\infty$. Since $\text{Im}\phi_n^{-1}$ does not assume the values $\{0, \beta_n\}$, the statement follows from Montel's theorem.

Take any limit function $\psi$ of the family $\{\phi_n^{-1}\}$. It is defined on the whole complement $\mathbb{C} \setminus \bar{\mathbb{D}}$ of the unit disk. We need to show that $\psi = \phi_f^{-1}$.

If $f^m z \in U' \setminus \bar{U}$ then for all sufficiently big $n$, $f_n^m z \in U' \setminus \bar{U}$ as well. As $N(f_n) \to \infty$, we conclude that $\Omega(f_n)$ is eventually contained in any neighborhood of $K(f)$. Hence the $\text{Im}\psi$ contains $\mathbb{C} \setminus K(f)$.

On the other hand, if $\text{Im}\psi \cap K(f) \neq \emptyset$, then $\text{Im}\psi \cap J(f) \neq \emptyset$. But then

$$\text{Im}\phi_n^{-1} \cap J(f_n) \neq \emptyset \tag{4.11}$$

since the Julia set depends lower semi-continuously on a quadratic-like map (as the repelling periodic points are stable under perturbations). But (4.11) is absurd.

Thus $\text{Im}\psi = \mathbb{C} \setminus K(f)$. By the normalization at $\infty$, $\psi$ must coincide with $\phi_f^{-1}$. □

4.8. *Banach slices of the foliation $\mathcal{F}$.* Let us say that the leaves $\mathcal{H}_c$ depend $C^1$-*continuously* on $c \in M_0$ if for any $c_0 \in M_0$, $g_0 \in \mathcal{E}$, and a Banach slice $\mathcal{E}_V \ni g_0$, there exist a Banach neighborhood $\mathcal{V} = \mathcal{E}_V(g_0, \varepsilon)$ and a Banach slice $\mathcal{B}_U \ni f_0$ such that $i_c \mathcal{V} \subset \mathcal{B}_U$ and

$$\mathrm{D}i_c(g) \to \mathrm{D}i_{c_0}(g) \quad \text{as} \quad c \to c_0,\ c \in M_0, \tag{4.12}$$

where the convergence is understood to be in the Banach operator norm, and where it is uniform over $g \in \mathcal{V}$.

LEMMA 4.16. *The leaves $\mathcal{H}_c$ depend $C^1$-continuously on $c \in M_0$.*

*Proof.* By the product structure of the foliation $\mathcal{F}$ (Theorem 4.13), the leaves depend $C^0$-continuously on $c \in M_0$: For any $c_0 \in M$, $g_0 \in \mathcal{E}$, and a Banach slice $\mathcal{E}_V \ni g_0$, there exist a Banach neighborhood $\mathcal{V} = \mathcal{E}_V(g_0, \varepsilon)$ and a Banach slice $\mathcal{B}_U \ni f_0$ such that $i_c \mathcal{V} \subset \mathcal{B}_U$ and $\|i_c(g) - i_{c_0}(g)\|_U \to 0$ as $c \to c_0$, $c \in M_0$, uniformly over $g \in \mathcal{V}$. Now the Cauchy inequality (see Appendix 2) yields (4.12). □

A Banach slice $\mathcal{F}_V$ of the foliation $\mathcal{F}$ is the restriction of $\mathcal{F}$ to the Banach space $\mathcal{B}_V$, so that the leaves of $\mathcal{F}_V$ are $\mathcal{H}_V(f) = \mathcal{H}(f) \cap \mathcal{B}_V$, $f \in \mathcal{C}$. We will show that the sufficiently deep Banach slices of $\mathcal{F}$ are still foliations with complex codimension-one analytic leaves (in the corresponding Banach space).

Let $E_V^h(f) = E^h(f) \cap \mathcal{B}_V(f)$ denote the slices of the horizontal spaces (4.5). Since the codimension of a subspace does not drop after restriction to a Banach slice (by the density property C1 from Appendix 2), the $E_V^h(f)$ are codimension-one subspaces in $\mathcal{B}_V$.

LEMMA 4.17. *For any $f_0 \in \mathcal{C}$ there exists a domain $V_0 \in \mathbb{V}_{f_0}$ such that for any $V \subset V_0$, $V \in \mathbb{V}_{f_0}$, the slice $\mathcal{F}_V$ near $f_0$ is a foliation in $\mathcal{B}_V$ with complex codimension-one analytic leaves.*

*Proof.* Let us first assume that $f_0 \in \text{int}\,\mathcal{C}$, so that $c_0 = \chi(f_0) \in \text{int}\,M$. By Lemma 4.8, $\mathcal{F}$ is analytic near $f_0$. Hence it has a local analytic transversal $\mathcal{S}$ parametrized near $c_0$ by $f = \phi(c)$, where $\chi \circ \phi = \text{id}$. Being a one-dimensional submanifold in $\mathcal{QG}$, $\mathcal{S}$ locally sits in some Banach slice $\mathcal{B}_V$.

Let us consider the Banach slice $\mathcal{F}_V$. The leaves of this foliation-to-be are the fibers of the straightening $\chi_V : \mathcal{QG}_V \to \mathbb{C}$, which is analytic near $f_0$. Since $\chi_V|\mathcal{S} = \chi|\mathcal{S}$ is nonsingular, it is a submersion near $f_0$. By the Implicit Function Theorem in Banach spaces (see [D1], [Lang]), the fibers of $\chi_V$ form a codimension-one analytic foliation in $\mathcal{QG}_V$ near $f_0$. (Note that the only requirement on the slice $\mathcal{B}_V$ is that some local transversal $\mathcal{S}$ should sit in it).

Let now $f_0 \in \partial \mathcal{C}$, $c_0 = \chi(f_0)$, $G_0 = \Pi(f_0) \in \mathcal{H}_0$, where the projection $\Pi : \mathcal{QG} \to \mathcal{H}_0$ is defined by (4.4). The vertical line $E^v(f_0)$ is naturally embedded into the space $\mathcal{B}(f_0)$ and hence into some Banach slice $\mathcal{B}_U \ni f_0$. Let us take a neighborhood $\mathcal{U} \subset \mathcal{B}_U$ which is analytically projected by $\Pi$ into some Banach neighborhood

$$\mathcal{V} \subset \mathcal{B}_V^0 \equiv \mathcal{H}_0 \cap \mathcal{B}_V = \{f \in \mathcal{B}_V : f(0) = 0\}$$

of $G_0$. Then

(4.13) $$D\Pi(E_U^h(f)) \subset T_{\Pi(f)}\mathcal{V} \approx \mathcal{B}_V^0, \quad f \in \mathcal{U}.$$

For $c \in M_0$ near $c_0$, let us consider immersed submanifolds $\mathcal{X}_c \subset \mathcal{QG}$ parametrized by $f = I_c(G)$ over $\mathcal{V}$ (recall that $I_c = i_c \circ \pi : \mathcal{H}_0 \to \mathcal{H}_c$).

By Lemma 4.16, these manifolds sit in some Banach slice $\mathcal{B}_W \supset \mathcal{B}_U$ and $C^1$ converge there to $\mathcal{X}_{c_0}$ as $c \to c_0$, $c \in M_0$. Hence the tangent planes $T_W(f) \subset \mathcal{B}_W$ to $\mathcal{X}_c$ at $f = I_c(G)$ are almost parallel to the tangent plane $T_W(f_0) \subset \mathcal{B}_W$ to $\mathcal{X}_{c_0}$ at $f_0$ (once $f$ is sufficiently close to $f_0$ in $\mathcal{B}_U$). Let us shrink $\mathcal{U}$ so that this takes place for all $f \in \mathcal{U}$.

Now consider the following decomposition: $\mathcal{B}_W = E_W^h(f_0) \oplus E^v(f_0)$. For a vector $u \in \mathcal{B}_W$, let $u^h$ and $u^v$ stand for the horizontal and vertical coordinates of $u$ with respect to this decomposition. Define the angle $\alpha = \text{ang}_W(u, E^v(f_0))$ by

$$\text{tg}\alpha = \frac{\|u^h\|_W}{\|u^v\|_W}, \quad \alpha \in [0, \pi/2].$$

Since the spaces $T_W(f)$ are almost parallel to $T_W(f_0)$ for $f \in \mathcal{U}$, the angle between any $u \in T_W(f)$ and $E^v(f_0)$ in $\mathcal{B}_W$ stays away from 0. If this $u$ actually belongs to $T_U(f) \equiv T_W(f) \cap \mathcal{B}_U$, then

$$\text{ang}_U(u, E^v(f_0)) \geq q\, \text{ang}_W(u, E^v(f_0)),$$

where $q = q(U, W) > 0$. Indeed, $\|u^h\|_U \geq \|u^h\|_W$, while $\|u^v\|_U \asymp \|u\|_W$, since all norms on a one-dimensional space are equivalent. Thus for $f \in \mathcal{U}$, the angle between $T_U(f)$ and $E^v(f_0)$ in $\mathcal{B}_U$ stays away from 0 as well.

But by (4.13), $E_U^h(f) \subset T_U(f)$. Hence, the horizontal spaces $E_U^h(f)$ also have a definite angle in $\mathcal{B}_U$ with the vertical line $E^v(f_0)$.

Let us now consider a local coordinate system in $\mathcal{U}$ provided by the decomposition $\mathcal{B}_U = E_U^h(f_0) \oplus E^v(f_0)$. Without loss of generality we can assume that $\mathcal{U}$ has a local product structure with respect to this decomposition, $\mathcal{U} = \mathcal{U}^h \times \mathcal{U}^v$. Let $L(f) \subset \mathcal{B}_U$ stand for the vertical lines through $f \in \mathcal{U}$ parallel to $E^v(f_0)$. As we have shown, these lines have a definite angle with the horizontal subspaces. In particular, for $f \in \mathcal{U} \cap \text{int}\mathcal{C}$, $L(f)$ provides a local transversal to $\mathcal{F}$ in the slice $\mathcal{B}_U$. As we have shown above, this implies that $\mathcal{F}_V$ is a codimension-one analytic foliation near $f$. Hence the leaf $\mathcal{H}_U(f)$ is the graph of an analytic function $\psi_f$ with a bounded slope over $\mathcal{U}^h$.

Take now a point $f \in \mathcal{U}$ and approximate it with a sequence $f_n \to f$, $f_n \in \mathcal{U} \cap \text{int}\mathcal{C}$. Then the slices $\mathcal{H}_U(f_n)$ uniformly converge in $\mathcal{B}_W$ to $\mathcal{H}_U(f)$. Hence the functions $\psi_{f_n}$ uniformly converge to a function $\psi_f$ parametrizing the slice $\mathcal{H}_U(f)$ (note that this statement does not depend on the choice of topology on $\mathcal{U}^h$). Since the uniform limit of analytic functions is analytic, we conclude that $\mathcal{H}_U(f)$ is a codimension-one analytic submanifold in $\mathcal{B}_U$.

Finally, the map $f \mapsto (p^h(f), \psi_f(f_0)) \in \mathcal{U}^h \times L(f_0)$ (where $p^h : \mathcal{B}_U \to E_U^h(f_0)$ stands for the horizontal projection) provides a local topological straightening of $\mathcal{F}_U$. $\square$

4.9. *Extension of the foliation.* Let us show that the foliation $\mathcal{F}_V$ admits a local smooth extension beyond $\mathcal{C}$. The leaves of this foliation are given by the position of the critical point in an appropriate local chart.

THEOREM 4.18. *For any $f_0 \in \mathcal{C}$ and any Banach slice $\mathcal{B}_U \ni f_0$, $U \in \mathbb{V}_f$, as in Lemma 4.17, there exists a Banach neighborhood $\mathcal{U} \subset \mathcal{B}_U$ of $f_0$ such that the foliation $\mathcal{F}_U$ admits an extension to $\mathcal{U}$ (with codimension-one complex analytic leaves) which is smooth on $\mathcal{U} \setminus \mathcal{C}$.*

*Proof.* Let us take a Banach slice $\mathcal{B}_U \ni f_0$ as in Lemma 4.17, so that $\mathcal{F}_U$ is a Banach foliation. Let $G_0 = \Pi(f_0)$ where $\Pi : \mathcal{QG} \to \mathcal{H}_0$ is projection (4.4) associating to a map $f \in \mathcal{QG}$ the map $G \in \mathcal{H}_0$ in the same vertical fiber. Then there exist Banach neighborhoods $\mathcal{U} \subset \mathcal{B}_U$ and $\mathcal{W} \subset \mathcal{H}_{0,V} \subset \mathcal{H}_0$ such that $\Pi(\mathcal{U}) \subset \mathcal{W}$. In what follows the neighborhoods $\mathcal{U}$ and $\mathcal{W}$ will be shrunk several times without change of the notation, keeping the above inclusions.

Let us select a fundamental annulus $A_0 = V_0' \setminus V_0$ for $G_0$ with a piecewise smooth boundary supplied with an invariant real analytic foliation by equipotentials. Then by Lemma 4.2, for maps $G \in \mathcal{H}_{0,V}$ sufficiently close to $G_0$, there is a choice of the fundamental annulus $A_G$ holomorphically moving with $G$ so that this motion respects the boundary dynamics and is smooth in both variables. This holomorphic motion admits an extension to a hybrid conjugacy $h_G : \mathbb{C} \to \mathbb{C}$ between $G_0$ and $G$ holomorphically depending on $G$. Let $\nu_G = h_G^* \sigma$ be the corresponding holomorphic family of Beltrami differentials.

Recall that for $f \notin \mathcal{C}$, $\xi(f)$ denotes the position of the critical value in the external model (see §3.3). Let us transfer this point to the $G$-plane:

$$(4.14) \qquad a_f = R_G^{-1}(\xi(f)),$$

where $G = \Pi(f)$ and $R_G : \mathbb{C} \setminus K(G) \to \mathbb{C} \setminus \bar{\mathbb{D}}$ is the Riemann mapping with positive derivative at $\infty$.

Take some point $a_0 \in V_0' \setminus K(G_0)$. By Theorem 3.4, there exists a quadratic-like map $f_* : U_* \to U_*'$ with disconnected Julia set such that $\Pi(f_*) = G_0$ and $a_{f_*} = a_0$. Consider the equipotential $\Gamma_0 \subset V_0$ through $G_0^{-1}(a_0)$ and the corresponding figure eight curve $\gamma_* \subset U_*$ centered at 0. Transfer the conformal structures $\nu_G$ to the exterior of $\gamma_*$ and pull it back by the iterates of $f_*$. We obtain a holomorphic family $\mu_G$ of $f_*$-invariant conformal structures, and the corresponding qc deformation $f_G$ of $f_*$ analytically depending on $G$. We define the leaves of the desired extension of $\mathcal{F}$ as the Banach slices of the holomorphic families $\{f_G\}_{g \in \mathcal{W}}$.

Let us show that for two different maps $\tilde{f}_* \neq f_*$ in the same vertical fiber (i.e., $\Pi(f_*) = \Pi(\tilde{f}_*) = G_0$), the corresponding families $f_G$ and $\tilde{f}_G$ are disjoint. Since $\Pi(f_G) = G$, they can intersect only at a point with the same $G$. But $a_{f_G} = h_G(a_0)$, since the deformation $f_G$ was obtained by lifting the Beltrami differential of $h_G$. Hence $a_{f_G} \neq a_{\tilde{f}_G}$, and the statement follows from Theorem 3.4.

Let us now consider the horizontally analytic tubing over $\mathcal{W}$:

$$(4.15) \qquad \Psi : \mathcal{W} \times \mathbb{A}(1,4) \to \mathbb{V} \equiv \bigcup_{G \in \mathcal{W}} V_G, \quad \Psi_G : \mathbb{A}(1,4) \to V_G,$$

such that the graphs $G \mapsto \Psi_G(z)$ coincide with the orbits if the holomorphic motion $\zeta \mapsto h_G(\zeta)$ (see Lemma 4.2 and (4.3)). Moreover, $(G, z) \mapsto \Psi_G(z)$ is smooth in two variables. Then the leaves of the extension of $\mathcal{F}$ to $\mathcal{U} \setminus \mathcal{C}$ are the fibers of the smooth map

$$(4.16) \qquad f \mapsto \Psi^{-1}_{\Pi(f)}(a_f).$$

This map is a submersion, since its restriction to any vertical fiber is a diffeomorphism. Hence its fibers form a smooth foliation.

Finally, the leaves of this foliation are complex analytic since they are the Banach slices of complex analytic families, and hence have complex tangent spaces. □

4.10. *$\mathcal{F}$ is transversally quasi-conformal.* Let us say that the foliation $\mathcal{F}$ is transversally *quasi-conformal* if the holonomy between two transversals $\mathcal{S}$ and $\mathcal{T}$ is locally a restriction of a qc map. Dilatation of the holonomy at $p \in \mathcal{S}$ is defined as the infimum of the dilatations of the local qc extensions. Let

$$\mathrm{mod}(\mathcal{S}) = \inf_{f \in \mathcal{S}} \mathrm{mod}(f).$$

THEOREM 4.19. *The foliation $\mathcal{F}$ is transversally quasi-conformal. The dilatation of the holonomy between two transversals $\mathcal{S}$ and $\mathcal{T}$ depends only on $\mu = \min(\mathrm{mod}(\mathcal{S}), \mathrm{mod}(\mathcal{T}))$.*

*Proof.* Let us take two transversals $\mathcal{S}$ and $\mathcal{T}$ to a leaf $\mathcal{H}$ of the foliation and a Beltrami path $\gamma$ in $\mathcal{H}$ joining two intersection points. Being compact, this path is contained in finitely many Banach slices (Lemma 11.6), whose number $N$ depends only on $\mu$. Hence by Lemma 4.18 this path can be covered by finitely many Banach balls $\mathcal{B}_{V_i}(f_i, \varepsilon_i)$ such that $\mathcal{F}$ extends to the twice bigger balls $\mathcal{B}_{V_i}(f_i, 2\varepsilon_i)$. Thus the holonomy between $\mathcal{S}$ and $\mathcal{T}$ can be decomposed into $N$ Banach holomorphic motions, which extend to the twice bigger domains. By the $\lambda$-lemma (see Appendix 2), each of the Banach motions is locally transversally quasi-conformal with uniform dilatation. □

Taking the quadratic family $\mathcal{Q}$ as one of the transversals, we obtain:

COROLLARY 4.20. *Let us consider a complex one-dimensional transversal $\mathcal{S} = \{f_\lambda\}$ to $\mathcal{F}$ in $\mathcal{QG}$. Then the straightening $\chi : \mathcal{S} \to \mathcal{Q}$ is locally $K$-quasi-conformal, with $K$ depending only on $\mathrm{mod}(\mathcal{S})$. Moreover, $K \to 1$ as $\mathrm{mod}(\mathcal{S}) \to \infty$.*

*Remarks.* 1. From the above point of view, the "miracle" of continuity of straightening in the quadratic-like case (see Douady [D2]) is directly related to the miracle of the $\lambda$-lemma. Note also that the source of discontinuity of the straightening for higher degrees is the failure of the $\lambda$-lemma for foliations of codimension bigger than 1.

2. The foliation $\mathcal{F}$ is *not transversally smooth*. For example, take the Ulam-Neumann quadratic $P = P_{-2} : z \mapsto z^2 - 2$ with a postcritical fixed point $\beta = 2$. Approximate it with superattracting parameter values $c_n \to -2$, for which $P_n^k(0) > 0$, $k = 2, \ldots, n-1$, while $P_n(0) = 0$, where $P_n = P_{c_n}$. Then $c_n - 2 \asymp 4^{-n}$, where 4 is the multiplier of $\beta$.

Let us now take another quadratic-like map $f \in \mathcal{H}_{-2}$ in the same hybrid class. Then in the vertical fiber $\mathcal{Z}_f$ there is a sequence of superattracting maps $f_n \in \mathcal{H}_{c_n}$ converging to $f$. Now the rate of convergence will be $\lambda^{-n}$, where $\lambda$ is the multiplier of the $\beta$-fixed point of $f$. Since the latter can be made different from 4, the holonomy $c_n \mapsto f_n$ is not smooth at $-2$. (To see that the multiplier can be efficiently changed in the hybrid class of $P$, consider, e.g., a quadratic-like deformation $P + \varepsilon Q$, where $Q$ is a polynomial with roots at $0, -2, 2$, and $Q'(2) \neq 0$.)

For the same reason the foliation is not smooth at other Misiurewicz points. Quasi-conformality seems to be the best transverse regularity of $\mathcal{F}$ which is satisfied everywhere. However, we will prove in subsections 7.2 and 9.3 that $\mathcal{F}$ is transversally smooth at Feigenbaum points.

4.11. *Full families.* Let $D$ be a Riemann surface. Consider a quadratic-like family $\mathcal{S} = \{f_\lambda\}, \lambda \in D$, over $D$, i.e., a complex analytic one-dimensional submanifold of $\mathcal{QG}$ parametrized by $D$. Such a family is called *full* if its Mandelbrot set $M_\mathcal{S} = \mathcal{S} \cap \mathcal{C}$ is compactly contained in $\mathcal{S}$. It is called *unfolded* if the straightening $\chi : M_\mathcal{S} \to M_0$ is injective.

If $D$ is a topological disk then by definition, the *winding number* of a full family over $D$ is the winding number of the critical value $\lambda \mapsto f_\lambda(0)$ around the critical point 0 as $\lambda$ runs once around a Jordan curve in $D$ which is close to $\partial D$ in the Hausdorff metric.

THEOREM 4.21 (Douady and Hubbard [DH2]). *Let $\mathcal{S}$ be a full quadratic-like family over $D$.*

- *If $\mathcal{S}$ is unfolded then the straightening $\chi : \mathcal{S} \to M_0$ is a homeomorphism;*

- *If $D$ is a topological disk then $\mathcal{S}$ is unfolded if and only if it has winding number 1.*

4.12. *Vertical fibers.* Recall that the external fibers $\mathcal{Z}_g = \pi^{-1}g$, $g \in \mathcal{E}$, are the fibers of the projection $\pi : \mathcal{QG} \to \mathcal{E}$. Let $M_g = \mathcal{Z}_g \cap \mathcal{C}$. In this section we will show that the external fibers are complex analytic curves. This is natural to expect since by Lemmas 4.4 and 4.5, $\pi$ is a submersion. However, since the Implicit Function Theorem is not available on the manifolds under consideration, we will give a special argument. First let us show that the fibers are topological curves:

LEMMA 4.22. *Let $g \in \mathcal{E}$. Then there is a neighborhood $N$ of the Mandelbrot set $M_0$ and a continuous injection $\gamma : N \to \mathcal{Z}_g$ such that for $c \in M_0$, $\gamma(c) = i_c(g)$.*

*Proof.* By Theorem 4.13, the map $M_0 \to \mathcal{Z}_g$, $c \mapsto i_c(g)$, is a homeomorphism onto $M_g$. Let us extend it beyond $M_0$.

We select a representative $g : V \to V'$ and conjugate it to $P_0 : z \mapsto z^2$ by a qc diffeomorphism

$$\psi : (V' \setminus \mathbb{D},\, V \setminus \mathbb{D}) \to (\mathbb{A}(1,4),\, \mathbb{A}(1,2)).$$

Let $\Omega \subset \mathbb{C} \setminus M_0$ be a collar around $M_0$ bounded by $\partial M_0$ and the parameter equipotential of level 4. We parametrize a collar around $M_g$ in $\mathcal{Z}_g$ by $\Omega$ as follows:

$$\gamma = \xi^{-1} \circ \psi^{-1} \circ \xi. \tag{4.17}$$

Since the middle map is smooth and the two others are conformal (by Lemma 4.14), $\gamma|\Omega$ is smooth as well. We will now show that it continuously matches $c \mapsto i_c(g)$ on $M_0$ (adopting the Douady and Hubbard argument for continuity of the straightening [DH2]). Let $c(k) \in \Omega$ converge to $c_* \in \partial M_0$, $f_k = \gamma(P_{c_k})$. For $c \in \Omega$, $P_c$ is $K$-qc conjugate to $f_c \equiv \gamma(P_c)$ with $K = \operatorname{Dil}\psi$. Hence the sequence $\{f_k\}$ is pre-compact and any limit $f$ of this sequence is qc conjugate to $P_{c_*}$. But since $c \in \partial M_0$, $P_{c_*}$ does not have invariant line fields on the Julia set (see §10.2). Hence $f$ is hybrid equivalent to $P_{c_*}$. Thus $\chi(f) = c_*$ and $\pi(f) = g$, and $f$ is uniquely defined. It follows that $f_{c(k)} \to f$. $\square$

THEOREM 4.23. *The vertical fibers $\mathcal{Z}_g$, $g \in \mathcal{E}$, are complex analytic curves.*

*Proof.* Let $\mathcal{Z} \equiv \mathcal{Z}_g$. First of all, by Lemma 4.8, in the interior of $\mathcal{C}$ the mating $c \mapsto i_c(g)$ provides a complex analytic parametrization of $\mathcal{Z}$.

Second, by Lemma 4.14, $\mathcal{Z} \setminus \mathcal{C}$ admits a complex analytic parametrization by the Riemann surface $\mathcal{S}_g$, so that $\mathcal{Z} \setminus \mathcal{C}$ is a holomorphic curve.

Let us prove that $\mathcal{Z}$ is analytic near any $f_0 \in \mathcal{Z} \cap \mathcal{C}$. To this end, consider the decomposition

$$\mathcal{B}(f_0) \approx \mathrm{T}_{f_0}\mathcal{QG} = E^h \oplus E^v$$

into the horizontal and vertical subspaces at $f_0$; see Lemma 4.7. Let $p^h : \mathcal{B}(f_0) \to E^h$ and $p^v : \mathcal{B}(f_0) \to E^v$ stand for the corresponding projections.

Let us first show that $p^v : \mathcal{Z} \to E^v$ is injective near $f_0$. Consider the map

$$A : \mathcal{U} \to E^h, \quad A = C \circ \pi,$$

where $\mathcal{U}$ is a neighborhood in $\mathcal{B}(f_0)$, and $C = \mathrm{D}i_a(g)$ with $a = \chi(f_0)$. Then $\mathrm{D}A(f_0) = p^h$.

For $u \in E^v$, let $H(u) = \{h \in \mathcal{B}(f_0) : p^v(h) = u\}$ stand for the hyperplane via $u$ parallel to $E^h \equiv H(0)$. Consider the restrictions $A_u$ of $A$ to these hyperplanes (defined on appropriate neighborhoods). Since all of these hyperplanes are naturally isomorphic, $A_u$ can be viewed as acting on the same space $E^h$. Then $\mathrm{D}A_0 = \mathrm{id}$ and by smoothness of $A$,

(4.18) $$\mathrm{D}A_u(f) \to \mathrm{id} \quad \text{as} \quad u \to 0, \ f \to f_0.$$

Take now a Banach slice $\mathcal{B}_V \ni f_0$ containing a neighborhood of $f_0$ in $\mathcal{Z}$. Let $h \in E_V^h$ and $\|h\|_V \leq 1$. By Lemma 11.5, for any $\varepsilon > 0$ there exist domains $W \Subset \Omega \Subset V$ in the family $\mathbb{V}_f$ such that

(4.19) $$\|h\|_W \geq \|h\|_\Omega^{1+\varepsilon}.$$

By (4.18),

$$\mathrm{D}A_u(f)h \in E_W^h \quad \text{and} \quad \|\mathrm{D}A_u(f)h\|_W \geq \frac{1}{2}\|h\|_W,$$

provided $(u, f)$ is close to $(0, f_0)$. Together with (4.19) it yields:

(4.20) $$\|\mathrm{D}A_u(f)h\|_W \geq q\|h\|_\Omega^{1+\varepsilon},$$

with some $q > 0$.

Let us take some $f \in \mathcal{Z}$ near $f_0$, so that $Af = f_0$. Then

$$A_u(f + h) = f_0 + \mathrm{D}A_u(f)h + O(\|h\|_\Omega^2),$$

where the big $O$ is uniform for $(u, f)$ near $(0, f_0)$. Incorporating (4.20), we obtain:

$$\|A_u(f + h) - f_0\|_W \geq \frac{q}{2}\|h\|_\Omega^{1+\varepsilon}.$$

Hence $A_u(f + h) \neq f_0$ for small enough $h \neq 0$. It follows that $\pi(f + h) \neq \pi(f) = g$ either, so that $f + h \notin \mathcal{Z}$.

Now, the projection $p^v : \mathcal{Z} \to E^v$ is injective near $f_0$. But by Lemma 4.22, $\mathcal{Z}$ is a topological curve. By the Open Mapping Theorem, the image $p^v \mathcal{Z}$ covers a neighborhood of 0. Thus, $\mathcal{Z}$ near $f_0$ is the graph of a continuous map $\psi : E^v \to E_V^h$.

Let us show that $\psi$ is differentiable. Select points $u$ and $u + \Delta u$ on $E^v$, and the corresponding points $f = \psi(u)$ and $f + \Delta f = \psi(u + \Delta u)$ on $\mathcal{Z}$. Let $B = \mathrm{D}A(f)$. Then

$$0 = A(f + \Delta f) - A(f) = B \cdot \Delta f + O(\|\Delta f\|_\Omega^2),$$

so that

$$\|B \cdot \Delta f\|_W = O(\|\Delta f\|_\Omega^2). \quad (4.21)$$

Note that the projection $p^v : E^v(f) \to E^v$ is nonsingular since $E^v(f)$ is close to $E^v$. Hence there exists a linear map $L : (E^v, 0) \to (\mathcal{B}(f_0), f)$ parametrizing the line $E^v(f)$ such that $L \circ p^v | E^v(f) = \text{id}$. Then $p^v f = u = p^v(Lu)$ and $p^v(f + \Delta f) = u + \Delta u = p^v(L(u + \Delta u))$. Hence

$$\Delta f = L \cdot \Delta u + \omega, \quad (4.22)$$

where $\omega \in E_V^h$ is a horizontal vector. Applying $B$ (taking into account that $B \cdot L = 0$), we conclude that $B \cdot \Delta f = B \cdot \omega$. Together with (4.21) this yields:

$$\|B\omega\|_W = O(\|\Delta f\|_\Omega^2).$$

But according to (4.20), $q\|\omega\|_\Omega^{1+\varepsilon} \leq \|B\omega\|_W$. Hence

$$\|\omega\|_\Omega = o(\|\Delta f\|_\Omega). \quad (4.23)$$

Together with (4.22) this yields:

$$\|L \cdot \Delta u\|_\Omega \geq \frac{1}{2}\|\Delta f\|_\Omega.$$

But since all one-dimensional norms are equivalent, $\|L \cdot \Delta u\|_\Omega \asymp \|\Delta u\|$ (with any choice of the latter norm). Hence $\|\Delta f\|_\Omega \leq \text{const} \cdot \|\Delta u\|$. By (4.23), $\|\omega\|_\Omega = o(\|\Delta u\|)$. By (4.22), $\psi$ is differentiable as a curve in $\mathcal{B}_\Omega$.

Thus, the external fibers are differentiable smooth curves. Since they have complex tangent lines, they are analytic. □

COROLLARY 4.24. *The vertical fibers are full unfolded quadratic-like families.*

4.13. *The transversality criterion.* The following lemma will give us an efficient way to check transversality of one-parameter families to the foliation $\mathcal{F}$:

LEMMA 4.25. *Let us consider an analytic one-dimensional submanifold $\mathcal{S} = \{f_\lambda\}$ in $\mathcal{QG}$, $f_0 \equiv f_{\lambda_0} \in \mathcal{C}$. If the straightening $\chi$ is locally injective on $\mathcal{M} \equiv \mathcal{S} \cap \mathcal{C}$ near $f_0$, then $\mathcal{S}$ is transverse to the foliation $\mathcal{F}$ at $f_0$.*

*Proof.* Injectivity of the straightening means that $\mathcal{S}$ intersects the leaves of $\mathcal{F}$ at single points. We should show that this yields transversality. Taking a Banach slice locally containing $\mathcal{S}$ and using Lemma 4.17 we reduce the situation to the Banach setting. By the Hurwitz Theorem (see Appendix 2), $\mathcal{S}$ is either transverse to $\mathcal{F}$ near $f_0$ or persistently tangent. But by the Intersection Lemma from the same appendix, the latter is impossible in $\text{int}\mathcal{C}$ where by Lemma 4.8 the foliation is transversally analytic. As $\text{int}\mathcal{C}$ is dense in $\mathcal{C}$, the conclusion follows. □

Thus full unfolded quadratic-like families (in particular the external fibers $\mathcal{Z}_g$) are transverse to the foliation $\mathcal{F}$ and we conclude:

THEOREM 4.26.    *If $\mathcal{S}$ is a full unfolded quadratic-like family then the straightening $\chi : M_\mathcal{S} \to M_0$ is a* qc *homeomorphism. Moreover, the dilatation of the straightening depends only on* $\mathrm{mod}(\mathcal{S})$ *and tends to* 0 *as* $\mathrm{mod}(\mathcal{S}) \to \infty$.

## 5. Renormalization, bounds and rigidity

5.1. *Little Mandelbrot copies.* Let us consider a hyperbolic component $H$ of the Mandelbrot set $M_0$ centered at the superattracting parameter value $c \in H$. Douady and Hubbard [DH2], [D2] proved that $H$ originates a "(little) Mandelbrot copy" $M = M_c$ canonically homeomorphic to the whole set. Let $\sigma = \sigma_M : M \to M_0$ be the corresponding homeomorphism. It transforms the component $H$ to the domain $H_0$ bounded by the main cardioid of $M_0$, so that $\sigma(c) = 0$. The inverse homeomorphism $M_0 \to M_c$ is called *tuning* and is denoted as $z \mapsto c * z$ (see Milnor [M]).

A basic combinatorial parameter of the Mandelbrot copy $M = M_c$ is its *period* $p_M \equiv p_c$ defined as the period of 0 under $P_c$. Except for the period doubling, there are several Mandelbrot copies with the same period. They are distinguished by their *combinatorics*, i.e., the Thurston type of the superattracting map $P_c$ (see [DH3]). Note that $M_c$ determines $c$ as the superattracting parameter value in $M_c$ with the smallest period. Thus we can use the copy itself to label the combinatorics.

The *root* $r_M$ of $M$ is the point corresponding to the cusp $1/4 \in M$ under the homeomorphism $\sigma$. A little Mandelbrot copy is called *primitive* if it is not attached at its root to any other hyperbolic component (geometrically it is recognized by the cusp of the originating component $H$ at its root). Otherwise it is called *satellite* (for such components, $\partial H$ is smooth at the root).

Let $\hat{M} = M$ in the primitive case, and $\hat{M} = M \setminus \{r_M\}$ ("unrooted" $M$) otherwise.

A Mandelbrot copy $M_c$ is called *real* if $c \in \mathbb{R}$, or equivalently $M_c$ is symmetric with respect to $\mathbb{R}$. The combinatorics of a real copy $M_c$ with period $p = p_c$ is determined by the order of the points $0, P_c(0), \ldots, P_c^{p-1}(0)$ on the real line.

A Mandelbrot copy is called *maximal* if it does not belong to any other copy except $M_0$ itself. These copies are pairwise disjoint, and any other copy, except $M_0$ itself, belongs to a unique maximal one (compare with the discussion of maximal renormalizations in the next section). All maximal copies are primitive except for the ones attached to the main cardioid. In particular, all real maximal Mandelbrot copies are primitive except for the one corresponding to the doubling bifurcation (i.e., with period 2).

All the copies $M \neq M_0$ are obtained from $M_0$ by iterated tunings

$$(5.1) \qquad M = c_l * \ldots c_1 * M_0 \equiv \sigma_{M_l}^{-1} \circ \ldots \circ \sigma_{M_1}^{-1} M_0,$$

where $M_k = M_{c_k}$ are maximal Mandelbrot copies. Thus any two Mandelbrot copies are either disjoint or nested.

Let $\mathcal{N}$ stand for the full family of the little Mandelbrot copies (not including $M_0$ itself), which is naturally identified with the set of all superattracting parameter values except 0. We will use $\mathcal{L}$ to denote a subfamily of pairwise disjoint copies of $\mathcal{N}$. Let $\mathcal{N}_{\max}$ stand for the family of maximal Mandelbrot copies.

5.2. *Renormalization.* For $M \in \mathcal{N}$, let $\mathcal{T}_M = \chi^{-1} M \subset \mathcal{C}$ (resp. $\mathcal{T}_{\hat{M}} = \chi^{-1} \hat{M} \subset \mathcal{T}_M$) stand for the union of the hybrid classes via $M$ (resp. $\hat{M}$). These sets will be called the (*horizontal*) *renormalization strips*. The strips $\mathcal{T}_M$ are closed (for instance, by the Product Structure Theorem 4.13). The renormalization strip is called *maximal* if it corresponds to a maximal Mandelbrot copy. Note that the maximal renormalization strips are pairwise disjoint.

There is a canonical renormalization operator $R_M : \mathcal{T}_{\hat{M}} \to \mathcal{C}$ defined as the $p = p_M$-fold iterate of $f$ restricted to an appropriate neighborhood $U$ of the critical point, up to rescaling. This neighborhood is selected in such a way that $g = f^p | U$ is a quadratic-like map with connected Julia set, and moreover the "little Julia sets" $f^k K(g)$, $k = 0, 1, \ldots, p-1$, are pairwise disjoint except, perhaps, touching at their $\beta$-fixed points (see [D2], [DH2], [L2], [McM1] for an extensive discussion of this notion). The maps $f \in \mathcal{T}_{\hat{M}}$ are called *renormalizable* with combinatorics $M$.

Among all renormalizations of a map $f$ there is the *maximal* one, with the smallest possible period (see [L2, §3.4], [McM1, §7.3]). It corresponds to the maximal renormalization strip containing $f$. Decomposition (5.1) can be rewritten as follows:

$$R_M = R_{M_1} \circ \ldots \circ R_{M_l},$$

where the $R_{M_k}$ are maximal renormalizations. In this sense any renormalization is induced by maximal ones.

Let us consider a family $\mathcal{L} \subset \mathcal{N}$ of pairwise disjoint Mandelbrot copies, e.g., $\mathcal{L} = \mathcal{N}_{\max}$. The operators $R_M$, $M \in \mathcal{L}$, can be unified into a single operator

$$R_{\mathcal{L}} : \bigcup_{M \in \mathcal{L}} \mathcal{T}_{\hat{M}} \to \mathcal{C}$$

whose restriction to a strip $\mathcal{T}_{\hat{M}}$ coincides with $R_M$. All operators $R_{\mathcal{L}}$ are induced by the maximal renormalization operator $R_{\max}$ corresponding to the family $\mathcal{L} = \mathcal{N}_{\max}$.

Similarly, the homeomorphisms $\sigma_M : M \to M_0$, $M \in \mathcal{L}$, can be unified into a single map
$$\sigma_\mathcal{L} : \bigcup_{M \in \mathcal{L}} M \to M_0.$$

The latter homeomorphism is the quotient of the renormalization by the straightening $\chi : \mathcal{C} \to M_0$:

(5.2) $$\sigma|\hat{M} = \chi \circ R_\mathcal{L}|\hat{M}.$$

(Note that $R_\mathcal{L}$ is not defined at the roots of satellite components, while $\sigma$ is extended over there.)

If $\mathcal{L}$ is a finite family then $R_\mathcal{L}$ is called a renormalization operator of *bounded type*. If $\mathcal{L} = \{M\}$ consists of a single Mandelbrot copy then $R_\mathcal{L} \equiv R_M$ is called the renormalization operator of *stationary type*.

We finish with the following useful fact:

LEMMA 5.1 (de Melo-van Strien [MvS, p. 440]). *The renormalization operator $R_{\max}$ is injective on the real slice of its domain.*

5.3. *Analytic extension.* For any $f_0 \in \mathcal{T}_{\hat{M}}$ and any Banach slice $\mathcal{B}_V \supset f_0$, the renormalization $R_M$ admits an analytic extension to a Banach neighborhood $\mathcal{B}_V(f_0, \varepsilon)$. That is, we take a quadratic-like representative $g_0 = f_0^p : U \to U'$ of the renormalization $R_M f_0$ with $U \Subset V$. Then for $f$ sufficiently close to $f_0$ in the Banach space $\mathcal{B}_V$, the restricted iterate $g = f^p|U$ represents a quadratic-like germ. Set by definition $R_M(f) = g$, up to rescaling.

As $R_M f$ is the normalized restricted iterate of $f$, $R_M : \mathcal{B}_V(f_0, \varepsilon) \to \mathcal{QG}$ is complex analytic. For instance, let us consider the doubling case when $R_M f$ corresponds to $f^2$. It is a composition of the second iterate operator $L : f \mapsto f^2|U$, $f : z \mapsto c + z^2 + \ldots$ and the normalization operator $N$. The former operator is obviously complex analytic with the differential

$$DL(f)v = (f' \circ f) \cdot v + v \circ f|U,$$

where $v : z \mapsto \delta + az^3 + \cdots$ is a vector field on $V$ with vanishing first and second order terms. The normalization operator $f \mapsto \lambda f(\lambda^{-1} z)$, where $\lambda = \lambda(f) = f''(0)/2$, is certainly analytic as well. (Note that $N$ transforms a small Banach neighborhood $\mathcal{B}_V(f, \varepsilon)$ to a Banach slice $\mathcal{B}_W$, where $W$ is a slightly shrunken domain $\lambda V$).

LEMMA 5.2. *There exists a $\rho > 0$ and a neighborhood $\mathcal{U}_M = \mathcal{U}_M(\mu, \rho)$ of the renormalization strip $\mathcal{T}_{\hat{M}}(\mu)$ in some slice $\mathcal{QG}(\mu, \rho)$ such that the renormalization $R_M$ admits an analytic extension to $\mathcal{U}_M$.*

*Proof.* Take a $\rho = \rho(\mu) > 0$ such that $\mathcal{C}(\mu) \subset \mathcal{QG}(\mu, \rho/2)$. Select a Montel distance $\text{dist}_M$ on $\mathcal{QG}(\mu, \rho)$, and take a small $\gamma > 0$. For any renormalizable quadratic-like map $g : U \to U'$ with $\mod(U' \setminus U) \geq \mu/2 > 0$, we have constructed local analytic extensions $R_g$ of the renormalization $R_M$ to the Banach neighborhood $\mathcal{B}_U(g, \gamma)$. We need to show that these local extensions glue together to form a single operator.

If $\gamma$ is sufficiently small then for all $g$ as above and all $f \in \mathcal{B}_U(g, \gamma)$:

(i) $\mod(R_g f) \geq \nu = \nu(\mu) > 0$;

(ii) $R_g f \in \mathcal{QG}^\#$, and the connected filled neighborhood $\Omega_\varepsilon(J(R_g f))$ from (3.1) can be selected with an $\varepsilon$ independent of $f$.

These properties follow from compactness of $\mathcal{C}(\mu)$ and the fact that the same renormalization domain can be used for a perturbed map.

Let us take a germ $f \in \mathcal{QG}(\mu, \rho) \setminus \mathcal{C}$. Assume there are two representatives $R_1 f \equiv Rg_1 f : V_1 \to V_1'$ and $R_2 f \equiv Rg_2 f : V_2 \to V_2'$ satisfying (i) and (ii) (where possibly $g_1$ and $g_2$ represent the same germ). Then $\|g_1 - g_2\|_U < 2\gamma$, and hence $\text{dist}_M(Rg_1, Rg_2) < \delta = \delta(\gamma)$, where $\delta \to 0$ as $\gamma \to 0$. But $\text{dist}_M(R_i f, Rg_i) < \delta$ as well. Thus $\text{dist}_M(R_1 f, R_2 f) < 3\delta$.

Since the Julia set $J(f)$ depends semi-continuously on the map (see [D3]), both points 0 and $R_1 f(0) = R_2 f(0)$ are contained in the same connected component of the intersection $\Omega_\varepsilon(J(R_1 f)) \cap \Omega_\varepsilon(J(R_2 f))$ (provided $\delta$ is sufficiently small). As $\Omega_{3\varepsilon}(J(R_1 f)) \cup \Omega_{3\varepsilon}(J(R_2 f))$ is contained in a connected component of $V_1' \cap V_2'$, the maps $R_1 f : V_1 \to V_1'$ and $R_2 f : V_2 \to V_2'$ represent the same germ. □

In what follows, referring to the analytic extension of $R_M$ beyond $\mathcal{T}_M$, we will mean the above extension to $\mathcal{U}_M(\mu, \rho)$ for some $\mu > 0$, $\rho > 0$.

5.4. *QC Theorem.* Let us pick a map $f \in \mathcal{T}_{\hat{M}}$ and a complex tangent line $E^t \subset \mathcal{B}(f)$ transverse to the leaf $\mathcal{H}(f)$, so that $E^t \oplus \mathrm{T}_f \mathcal{H}(f) = \mathcal{B}(f)$. We say that $R$ is *transversally nonsingular* at $f$ if the restriction of the differential $DR_f$ to $E^t$ is nonsingular. (Since the foliation $\mathcal{F}$ is $R$-invariant, this definition is independent of the choice of $E^t$.)

LEMMA 5.3 (Transversal nonsingularity). *The renormalization is transversally nonsingular at any $f \in \mathcal{T}_{\hat{M}}$.*

*Proof.* Let $\mathcal{S}$ be a one-dimensional local transversal to $\mathcal{F}$ through $f$, and $\mathcal{M} = \mathcal{S} \cap \mathcal{T}_{\hat{M}}$. Then the straightening $\chi : \mathcal{M} \to M$ is injective. Since $\sigma : M \to M_0$ is also injective, (5.2) yields the injectivity of $\chi : R(\mathcal{M}) \to M_0$. By Lemma 4.25, $R(\mathcal{S})$ is transverse to $\mathcal{F}$ at $Rf$. □

Douady and Hubbard gave a sufficient condition for the canonical homeomorphism $\sigma : M \to M_0$ to be qc ([DH2, Prop. 22]). We will show now that it is always so in the primitive case.

LEMMA 5.4. *A primitive copy $M$ of the Mandelbrot set $M_0$ is locally* qc *equivalent to the whole set $M_0$.*

*Proof.* Let us consider the analytic extension of $R_M$ to a neighborhood $D \subset \mathbb{C}$ of $M$.

By (5.2), the regularity of the homeomorphism $\sigma : M \to M_0$ is ruled by the regularity of the straightening $\chi$. By Theorem 4.19, the latter is locally qc on analytic transversals to $\mathcal{F}$. As the analytic family $\mathcal{S} = RD$ is transverse to $\mathcal{F}$ (by Lemma 5.3), the conclusion follows. $\square$

THEOREM 5.5 (The QC Theorem). *A primitive copy $M$ of the Mandelbrot set $M_0$ is* qc *equivalent to the whole set $M_0$.*

*Proof.* We will use the notations of the previous lemma. The straightening $\chi : RM \to M_0$ admits a continuous extension to a neighborhood $N$ of $RM$ in the transversal $\mathcal{S}$ which is qc on $N \setminus RM$ (see [L4, Lemma 3.1]). By the Gluing Lemma from Appendix 1, this extension glues with the local qc extensions provided by Lemma 5.4 into a single qc homeomorphism (see Lemma 3.2 of [L4] for details). $\square$

*Remark.* By the same argument, a satellite Mandelbrot set $M$ is qc equivalent to $M_0$ after removal of neighborhoods of the roots. Presumably the whole satellite set $M$ is qc equivalent to the "one half" of the Mandelbrot set $\tilde{M}_0$ of the family $z \mapsto \lambda z + z^2$. (Note that the latter is a holomorphic double branched covering of $M_0$ by the map $c = \lambda/2 - \lambda^2/4$ branched at $\lambda = 1$ over the cusp $c = 1/4$.)

5.5. *Combinatorial type.* From now until the end of Section 5 we fix a family $\mathcal{L} \subset \mathcal{N}$ of disjoint Mandelbrot copies, and let $R = R_{\mathcal{L}}$. A map $f \in \mathcal{C}$ is called $N$ times renormalizable by $R$ ($0 \leq N \leq \infty$), if

$$R^n f \in \bigcup_{M \in \mathcal{L}} \mathcal{T}_{\hat{M}}, \quad n = 0, 1, \ldots, N-1, \quad \text{and } R^N f \in \mathcal{C}.$$

The *itinerary* $\tau_N(f)$ of such a map $f$ is the sequence $\tau(f) = \{M_0, M_1, \ldots M_{N-1}\}$ of copies $M_n \in \mathcal{L}$ such that $R^n f \in \mathcal{T}_{M_n}$. One says that the combinatorics of such an $f$ is bounded by $\bar{p}$ if $p(M_n) \leq \bar{p}$, $n = 0, 1, \ldots, N$.

The itinerary $\tau(f) \equiv \tau_\infty(f)$ of an infinitely renormalizable map is also called its *combinatorial type*. Two infinitely renormalizable maps are called *combinatorially equivalent* if they have the same combinatorial type.

Let us now consider an orbit

(5.3) $$\{f_n = R_{M_{n-1}} \cdot \ldots \cdot R_{M_0} f \equiv R^n f\},$$

where the maps $f_n$ can have a disconnected Julia set and the $R_{M_k}$ are understood as the analytic continuation of the renormalization. To keep the notation simple, we will still denote $f_n$ as $R^n f$ keeping in mind its meaning.

The combinatorial type of the orbit (5.3) is naturally defined as the string $\{M_0, M_1, \ldots\}$. Somewhat loosely, it will also be called the combinatorial type of $f$. Accordingly two orbits as above (or the corresponding germs) are called combinatorially equivalent if they have the same combinatorial type.

5.6. *A priori bounds.* A real quadratic-like map $f$ is close to the cusp if it has an attracting fixed point with multiplier at least $1/2$ (one can use any $1 - \varepsilon$ in place of $1/2$ but then the bounds below will depend on $\varepsilon > 0$).

THEOREM 5.6 (A priori bounds). *Let $f$ be a real $N$ times renormalizable quadratic-like germ with itinerary $\tau_N(f) = \{M_0, M_1, \ldots, M_{N-1}\}$. Assume that $p(M_k) \leq \bar{p}$ and $\mathrm{mod}(f) \geq \nu > 0$. Then there exist $\mu = \mu(\bar{p}) > 0$ and $l = l(\nu)$ such that*

$$\mathrm{mod}(R^n f) \geq \mu > 0, \ n = l, \ldots, N - 1.$$

*Moreover, $\mathrm{mod}(R^N f) \geq \mu$ as well, unless the last renormalization is of doubling type and $R^N f$ is close to the cusp.*

An infinitely renormalizable germ is said to have *a priori* bounds if $\mathrm{mod}(R^n f) \geq \mu > 0, \ n = 0, 1, \ldots$.

COROLLARY 5.7. *Any real infinitely renormalizable quadratic-like germ $f$ with bounded combinatorics has a priori bounds. More precisely, if the combinatorics of $f$ is bounded by $\bar{p}$ and $\mathrm{mod}(f) \geq \nu > 0$, then there exist $\mu = \mu(\bar{p}) > 0$ and $l = l(\nu)$ such that*

$$\mathrm{mod}(R^n f) \geq \mu, \quad n = l, l + 1, \ldots.$$

*Remark.* The existence of *a priori* bounds for maps with bounded combinatorics was proved in [MvS], [S2]. The refined finitely renormalizable version appeared in [LS], [LY]. The latter works actually prove that the above bounds are independent of $\bar{p}$.

For an $n$ times renormalizable *map* $f_V : V \to V'$, let us say that a quadratic-like representative $R^n f_V : V_n \to V'_n$ is *subdued* to $f$ if it is a restricted iterate of the map $f_V$ itself (so that no analytic continuation of $f$ is allowed). The family of subdued quadratic-like maps represents the *subdued renormalization germ* which will also be denoted as $R^n f_V$. For a subdued germ $g$, $\mathrm{mod}(g)$ means the supremum of moduli of the subdued representatives.

It is easy to see that if the renormalizations of a germ $f$ have bounds $R^n f \geq \mu, \ n = 0, 1, \ldots, N$, then the subdued renormalizations have delayed bounds:

(5.4) $$R^n f_V \geq \mu/2, \quad n = l(\mu, \mathrm{mod}(V' \setminus V)), \ldots, N.$$

5.7. *Combinatorial rigidity.* An orbit $\{R^n f = R_{M_{n-1}} \cdot \ldots \cdot R_{M_0} f\}_{n \in \mathbb{N}}$ is called *nonescaping* if there exist $\mu > 0$, $\rho > 0$ such that $R^n f \in \mathcal{U}_{M_{n+1}}(\mu, \rho)$ for all $n \in \mathbb{N}$, where the $\mathcal{U}_{M_n}(\mu, \rho)$ are the domains of analyticity of the $R_{M_n}$ constructed in subsection 5.3. We will also say that $f$ is nonescaping (keeping in mind that this notion depends on the choice of branches $R_{M_n}$). In particular, $f$ has *a priori* bounds: $\mathrm{mod}(R^n f) \geq \mu > 0$, $n = 0, 1, \ldots$.

LEMMA 5.8.   *If an orbit $\{R^n f\}_{n=0}^{\infty}$ is nonescaping then $f$ is infinitely renormalizable.*

*Proof.* Let $f_n$ stand for non-rescaled germs representing subdued renormalizations $R^n f_V$. Since $f = f_V$ is nonescaping, there exist quadratic-like representatives $f_n = f^{r_n} : U_n \to V_n$ such that $U_n$ and $V_n$ have bounded geometry.

Let us show $\mathrm{diam} U_n \to 0$. Otherwise, there would be a disk $\mathbb{D}_r$ contained in all $U_n$. Since $\deg(f^{r_n}|\mathbb{D}_r)) \leq 2$, $\mathbb{D}_r$ does not intersect the Julia set $J(f)$. But then $\mathbb{D}_r$ is escaping under iterates of $f$, so that the big iterates $f^{r_n}$ are not well-defined on $\mathbb{D}_r$.

It follows that the Julia set $J(f)$ is connected. Indeed, since $J(f) \cap U_n \supset J(f_n) \neq \emptyset$, we have: $\mathrm{dist}(0, J(f)) \leq \mathrm{diam} U_n \to 0$. Hence $0 \in J(f)$. For the same reason, all Julia sets $J(f_n)$ are connected.

By definition, $f_{n+1} = R_{M_n} f_n$ where $R_{M_n}$ means the analytic extension of the renormalization. But once $J(f_{n+1})$ is connected, $f_n$ is renormalizable and $f_{n+1}$ is its canonical renormalization. Hence $f$ is infinitely renormalizable. □

THEOREM 5.9 (Combinatorial rigidity).   *Consider two nonescaping germs $f_1$ and $f_2$ in $\mathcal{QG}$ with bounded combinatorics. If $f_1$ and $f_2$ are combinatorially equivalent then they are hybrid equivalent.*

*Proof.* By Lemma 5.8, $f_1$ and $f_2$ are infinitely renormalizable quadratic-like germs with *a priori* bounds. By the Rigidity Theorem of [L2], $f_1$ and $f_2$ are hybrid equivalent. □

5.8. *McMullen towers.* Let $-\infty \leq l \leq 0 \leq n \leq \infty$. By definition, an $(l, n)$-*tower* $\mathbf{f}$ (related to the renormalization operator $R = R_\mathcal{L}$) is a sequence of quadratic-like germs $f_k : V^k \to U^k$ with connected Julia set, $k = l, \ldots, n$, such that $f_k = f_{k-1}^{r_k} | V_k$, where $f_{k-1}^{r_k} | V_k$ represents the renormalization $R f_{k-1}$. The germs can be *simultaneously* rescaled, so that $f_0$ can be normalized as $z \mapsto c + z^2 + \ldots$. A tower is called infinite if $-l = \infty$, and it is called bi-infinite if $-l = n = \infty$.

*The combinatorics* $\tau(\mathbf{f})$ of the tower $\mathbf{f}$ is the string $\{M_l, \ldots, M_{n-1}\}$ of little Mandelbrot copies such that $\chi(f_k) \in M_k$. We say that the tower has a *bounded combinatorics* if there are only finitely many different copies in this

string. Let
$$p(\mathbf{f}) = \sup_{l \le k \le n} p(f_k).$$
We say that combinatorics of the tower is bounded by $\bar{p}$ if $p(\mathbf{f}) \le \bar{p}$.

Let
$$\mathrm{mod}(\mathbf{f}) = \inf_{l \le k \le n} \mathrm{mod}(f_k).$$
We say that a tower has *a priori bounds* if $\mathrm{mod}(\mathbf{f}) > 0$.

The space of towers is supplied with the weak topology: $\{g_{m,k}\}_k \equiv \mathbf{g}_m \to \mathbf{f}$ if:

- Given a $k$, the coordinates $g_{m,k}$ are eventually well-defined if and only if the coordinate $f_k$ is well-defined as well;

- For each $k$ with a well-defined $f_k$, $g_{m,k} \to f_k$ as $m \to \infty$.

Note that, in particular, finite towers can converge to an infinite one.

By means of the diagonal procedure, Lemma 4.1 yields:

LEMMA 5.10. *Take $\bar{p}$ and $\mu > 0$. The set of towers with $p(\mathbf{f}) \le \bar{p}$ and $\mathrm{mod}(\mathbf{f}) \ge \mu > 0$ is sequentially compact.*

The filled Julia set $K(\mathbf{f})$ of a tower is defined as the union $\cup K(f_k)$ (without the closure).

THEOREM 5.11 (Hairiness of the Julia set [McM2]). *If $\mathbf{f}$ is an infinite tower with bounded combinatorics and a priori bounds then the filled Julia set $K(\mathbf{f})$ is dense in $\mathbb{C}$.*

Two $(l,n)$-towers $\mathbf{f}$ and $\mathbf{g}$ are called topologically conjugate if there is a homeomorphism $h$ defined in a neighborhood of $K(\mathbf{f})$ which simultaneously conjugates each $f_k$ to $g_k$. A self-conjugacy of some tower with itself is called its *automorphism*. The last theorem together with Lemma 3.5 yield:

COROLLARY 5.12 (No automorphisms). *An infinite tower $\mathbf{f} = \{f_k\}$ with a priori bounds and empty $\mathrm{int} K(f_0)$ does not admit nontrivial automorphisms. In particular, bi-infinite towers with a priori bounds do not admit nontrivial automorphisms.*

If a conjugacy $h$ between two towers can be selected to be quasi-conformal then the towers are called qc conjugate. If additionally $\bar{\partial} h = 0$ a.e. on $K(\mathbf{f})$ then the towers are called *hybrid equivalent*.

LEMMA 5.13. *Two towers $\mathbf{f} = \{f_k\}_{k=l}^n$ and $\mathbf{g} = \{g_k\}_{k=l}^n$ with bounded combinatorics and a priori bounds are qc equivalent if and only if all pairs $f_k$ and $g_k$ are $K$-qc equivalent with uniform $K$.*

*Proof.* Since a qc conjugacy between $f_k$ and $g_k$ serves as a qc conjugacy between $f_s$ and $g_s$ for all $s \geq k$, the statement is not totally obvious only when $l = -\infty$.

First note that by Lemma 3.6, $f_k$ and $g_k$ are $L$-qc conjugate by a map $h: V_k \to U_k$ such that $\mod(V_k \setminus K(f_k)) \geq \nu > 0$ and $\mod(U_k \setminus K(g_k)) \geq \nu > 0$, where $\nu = \nu(\mu) > 0$, and $L = L(K, \mu)$.

Second, the diameters of $J(f_k)$ and $J(g_k)$ grow exponentially as $k \to -\infty$ (see [McM2, Prop. 8.1]).

It follows that the domains $U_k$ and $V_k$ exhaust the plane as $k \to -\infty$. Since the space of normalized $K$-qc maps is compact, we can select a subsequence converging to a conjugacy between the towers. □

THEOREM 5.14 (The Tower Rigidity Theorem [McM2]).

(i) *If two bi-infinite towers with bounded combinatorics and a priori bounds are quasi-conformally equivalent then they are affinely equivalent.*

(ii) *If two infinite towers $\{f_k\}_{k=0}^{-\infty}$ and $\{\tilde{f}_k\}_{k=0}^{-\infty}$ with bounded combinatorics and a priori bounds are quasi-conformally equivalent and $\chi(f_0) = \chi(\tilde{f}_0)$, then they are affinely equivalent.*

## 6. Hyperbolicity of the renormalization (the stationary case)

6.1. *A renormalization fixed point and its stable manifold.* Throughout this section $R \equiv R_M$ will stand for a renormalization operator with stationary combinatorics $M \in \mathcal{N}$.

Let us consider a renormalization fixed point $f_*$, $Rf_* = f_*$.

*Definition* 6.1. Given an invariant set $\mathcal{W} \subset \mathcal{QG}$, let us say that the orbits of $\mathcal{W}$ *uniformly exponentially converge to* $f_*$ if for any quadratic-like germ $f \in \mathcal{W}$, the orbit $\{R^n f\}_{n \geq N(\mod(f))}$, belongs to some $\mathcal{B}_V \ni f_*$ and uniformly exponentially converges to $f_*$ in this Banach space:

$$\|R^n f - f_*\|_V \leq Cq^{n-N},$$

where $C > 0$ and $q \in (0, 1)$ are independent of $f$.

*Remark.* Since the property of exponential convergence is Hölder invariant, it can be understood in the sense of the natural Hölder structure on the precompact sets $\mathcal{W}(\mu)$ (see §11.3).

Let us define the *stable manifold* of $f_*$ as

$$\mathcal{W}^s(f_*) \equiv \mathcal{W}^s_* = \{f \in \mathcal{QG} : R^n f \to f_*\}.$$

The following theorem summarizes results of Sullivan [S2] and McMullen [McM2], and the author ([L2] and this work).

THEOREM 6.1 (Stable Manifold). *Assume that there exists an infinitely renormalizable (under the operator R) map F with a priori bounds. Then the operator R has a unique fixed point $f_*$. The stable manifold $\mathcal{W}^s_* \equiv \mathcal{W}^s(f_*)$ of this point is a complex analytic submanifold in $\mathcal{QG}$ of codimension* 1 *coinciding with the hybrid class $\mathcal{H}_* = \mathcal{H}(f_*)$. The orbits in $\mathcal{W}^s_*$ uniformly exponentially converge to $f_*$.*

*Proof.* The maps $F$ and $RF$ are combinatorially equivalent and have *a priori* bounds. By the Combinatorial Rigidity Theorem, they are hybrid equivalent, so that the hybrid class $\mathcal{H}_* = \mathcal{H}(F)$ is $R$-invariant.

By Lemma 4.1, the orbit $\{R^n F\}_{n=0}^\infty$ is pre-compact, so that its $\omega$-limit set $\Omega$ is compact. Since $R|\Omega$ is obviously surjective and $\mathrm{mod}(g) \geq \mu > 0$ for all $g \in \Omega$, any $f \equiv f_0 \in \Omega$ can be included into the two-sided tower $\mathbf{f} = \{f_k \in \Omega\}_{k=-\infty}^\infty$ with stationary combinatorics and *a priori* bounds.

Take two such towers $\mathbf{f}$ and $\mathbf{g}$. Since $f_k$ and $g_k$ belong to the same hybrid class $\mathcal{H}_*$, by Lemma 5.13, these towers are quasi-conformally equivalent. By the Towers Rigidity Theorem, they are affinely equivalent. In particular, $f_0 = g_0$, so that $\Omega$ consists of a single fixed point $f_* \equiv f_0$.

It follows that $R^n f \to f_*$ for any $f \in \mathcal{H}_*$. Moreover, this convergence is uniform in the following sense:

*Statement.* There exists a quadratic-like representative $f_* : V \to V'$ with the following property. For any $\nu > 0$ and $\varepsilon > 0$, there exists an $N = N(\nu, \varepsilon)$ such that: If $\mathrm{mod}(f) \geq \nu > 0$ then for $n \geq N$, $R^n f \in \mathcal{B}_V(f_*, \varepsilon)$.

Note first that Lemma 3.6 and the fact that a conjugacy between $f_*$ and $f$ restricts to a conjugacy between their renormalizations imply that $R^n f \in \mathcal{H}_*(\eta)$, $n = 0, 1, \ldots$, where $\eta = \eta(\nu)$. Let us show that $\mathrm{dist}_\mathrm{M}(R^n f, f_*) \leq \varepsilon$ for $n \geq N(\nu, \varepsilon)$, where $\mathrm{dist}_\mathrm{M}$ is the Montel distance on $\mathcal{H}_*(\eta)$.

Otherwise we would find a sequence of maps $f_m \in \mathcal{H}_*(\nu)$ and moments $n_m \to \infty$ such that $\mathrm{dist}_\mathrm{M}(R^{n_m} f_m, f_*) \geq \delta$. Let $h_m = R^{n_m} f_m$. Consider towers $\mathbf{h}_m = \{R^k h_m\}_{k=-n_m}^\infty$. As these towers have *a priori* bounds, by the Compactness Lemma 5.10 we can select diagonally a sequence converging to a two-sided tower with *a priori* bounds. Since $\mathrm{dist}_\mathrm{M}(h_m, f_*) \geq \delta$, this tower is different from the stationary tower $(\ldots, f_*, f_*, f_*, \ldots)$, which contradicts the rigidity of towers.

Moreover, there exists a quadratic-like representative $f_* : V \to V'$ (*a priori* depending on $\nu$) such that $R^n f \in \mathcal{B}_V(f_*, \varepsilon)$ for $n \geq N(\nu, \varepsilon)$. Otherwise let us consider a nested sequence $V_1 \supset V_2 \supset \ldots$ of domains shrinking to $K(f_*)$, and find a sequence of germs $f_m \in \mathcal{H}_*(\nu)$ and moments $n_m \to \infty$ such that $R^{n_m} f_m \notin \mathcal{B}_{V_m}$. But as we have just shown, $R^{n_m} f_m \to f_*$ in $\mathcal{H}_*(\eta)$. This means that there exists a quadratic-like representative $f_* : V \to V'$ such that $R^{n_m} f_m \in \mathcal{B}_V$ for all sufficiently big $m$. As this $V$ contains some $V_m$, we arrive at a contradiction.

To complete the proof of the statement we need to show that the same is true for a domain $V$ independent of $\nu$. Take a representative $R^l f_W : W_l \to W'_l$ $\varepsilon$-close to $f_{*,V}$ in $\mathcal{B}_V$. If $\varepsilon$ is small enough, they can be $(1+\delta)$-qc conjugate in slightly smaller domains. This conjugacy provides a $(1+\delta)$-qc conjugacy between the further renormalizations $R^{l+m} f_W : W_{l+m} \to W'_{l+m}$ and $R^m f_{*,V} : V_m \to V'_m$, $m \geq 0$, subdued to the above representatives. But since $\text{mod}(R^m f_*) = \text{mod}(f_*) > 0$, the subdued renormalizations $R^m f_{*,V}$ are eventually (for $m \geq N = N(\eta)$) well defined on the same domain $U$; see (5.4). Hence for $m \geq N$, $R^{l+m} f_W$ is well defined on a slightly smaller domain and is close to $f_*$ there, and the statement follows.

Let us now consider the analytic diffeomorphism $\Pi : \mathcal{H}_* \to \mathcal{H}_0$ (4.4) and the inverse map $I_* : \mathcal{H}_0 \to \mathcal{H}_*$. Let $G_* = \Pi(f_*)$, and

$$R_0 = \Pi \circ R \circ I_* : \mathcal{H}_0 \to \mathcal{H}_0.$$

Then $\Pi(\mathcal{B}_V(f_*, r)) \subset \mathcal{H}_{0,W}$ for some $r > 0$ and some Banach slice $\mathcal{H}_{0,W} \ni G_*$, and this Banach restriction is continuous. It follows that the orbits of $R_0$ uniformly converge to $G_*$: For any sufficiently small $\varepsilon > 0$ and $\delta > 0$ there is an $N$ such that

$$R_0^N \mathcal{H}_{0,W}(G_*, \varepsilon) \subset \mathcal{H}_{0,W}(G_*, \delta).$$

By the Schwarz Lemma (see Appendix 2), $R_0^N$ is uniformly contracting if $\delta < \varepsilon/2$. Thus the orbits of $R_0$ converge to $G_*$ exponentially fast in the $\|\cdot\|_W$-norm.

Finally, there exists a Banach slice $\mathcal{B}_U \ni f_*$ such that for sufficiently small $\varepsilon$, $I_* \mathcal{H}_{0,W}(G_*, \varepsilon) \subset \mathcal{B}_U$, and this Banach restriction is continuous. It follows that the orbits of $f \in \mathcal{B}_V(f_*, r)$ converge to $f_*$ exponentially fast in the $\|\cdot\|_U$-norm, hence in the Montel metric on $\mathcal{B}_V$.

So, we have proved that $\mathcal{H}_* \subset W_*^s$ and that the orbits in $\mathcal{H}_*$ uniformly exponentially converge to $f_*$. The opposite inclusion, $W_*^s \subset \mathcal{H}_*$, follows from Theorem 5.9. Thus $W_*^s = \mathcal{H}_*$, and by Theorem 4.11 this is a codimension-one complex analytic submanifold in $\mathcal{QG}$. □

*Remarks.* 1. In [McM2], [S2] the following extra assumption was needed: $RF$ is *hybrid* equivalent to $F$ (which was proved by Sullivan for real $F$). The Combinatorial Rigidity Theorem 5.9 allows us to eliminate this assumption. The above proof of existence of the renormalization fixed point which attracts all the hybrid class is due to McMullen [McM2]. However, the proof of the exponential convergence based on the Schwarz Lemma is new. The inclusion $\mathcal{H}_* \subset W^s$ is due to Sullivan and McMullen but the opposite inclusion $W_*^s \subset \mathcal{H}_*$ is new. The statement that $W_*^s$ is a codimension-one analytic submanifold is also new.

2. Note that the above argument does not use *uniform a priori* bounds, i.e., the bounds which are eventually independent of the map in question. Vice versa, it shows how the uniform bounds follow from the relative ones.

A fixed point $f_*$ is called *attracting* if it has a neighborhood $\mathcal{U} \subset \mathcal{QG}$ contained in the stable manifold $W^s(f_*)$.

COROLLARY 6.2. *Fixed points of the renormalization operator are not attracting.*

*Proof.* Otherwise the stable manifold $W^s(f_*)$ would have codimension 0 rather than 1. □

6.2. *Hyperbolicity.* We are now ready to prove hyperbolicity of the renormalization transformation at its fixed point $f_*$. Let $R_* \equiv DR(f_*)$ : $T\mathcal{H}_* \to T\mathcal{H}_*$ stand for the differential of $R$ at $f_*$. Note that the tangent space $T\mathcal{H}_*$ is naturally identified with the space of germs of analytic vector fields $z \mapsto a + bz^3 + \ldots$ near $K(G_*)$, where $G_* = \Pi(f_*) \in \mathcal{H}_0$. Thus it has a natural structure of the inductive limit of Banach spaces. We say that $R_*$ is uniformly exponentially contracting in this space if its iterates uniformly exponentially converge to 0 (which is defined in the same way as in the nonlinear situation; see the previous section).

THEOREM 6.3 (Hyperbolicity). *The tangent space $\mathcal{B}_* \equiv T_{f_*}\mathcal{QG}$ admits an $R_*$-invariant splitting $\mathcal{B}_* = E^s \oplus E^u$, where $E^s = T\mathcal{H}_*$ and $\dim E^u = 1$. Moreover, $R_*|E^s$ is uniformly exponentially contracting, while the absolute value of the eigenvalue $\lambda_*$ of $R_*|E^u$ is greater than 1.*

*Proof.* By Theorem 6.1, the map $R|\mathcal{H}_*$ is uniformly exponentially contracting. By the Schwarz lemma, its differential $R_*|E^s$ is uniformly exponential contracting as well.

Let us make a selection of Banach spaces. Select a quadratic-like representative $f_{*,V} : V \to V'$ of the renormalization fixed point $f_*$ so that the requirements of Definition 6.1 are satisfied on the stable manifold $\mathcal{W}^s(f_*)$. In particular, for any representative $f_{*,W} : W \to W'$ there exists an $N = N(W)$ such that the subdued renormalization $R^N f_{*,W}$ has a representative $U \to U'$ with $V \Subset U$ (and of course $R^N f_{*,W}|V = f_{*,V}$). Then for any $\delta > 0$ there is an $\varepsilon > 0$ such that for all maps $f \in \mathcal{B}_W(f_*, \varepsilon)$, the renormalization $R^N f_W$ belongs to $\mathcal{B}_V(f_*, \delta)$. Thus for $W \Subset V$, $R^N$ gives rise to a Banach operator $A : \mathcal{B}_W(f_*, \varepsilon) \to \mathcal{B}_W$ which is the composition of the following two operators:

$$\mathcal{B}_W(f_*, \varepsilon) \underset{R^N}{\to} \mathcal{B}_V \underset{i}{\hookrightarrow} \mathcal{B}_W,$$

where $i \equiv i_{WV}$ is the natural inclusion. Since the latter embedding is compact, $A$ is compact as well. Since $R$ is complex analytic, $A$ is complex analytic as well. Let $A_*$ stand for the differential of $A$ at $f_*$. This linear operator is compact as well since it is also factored via the inclusion $i : \mathcal{B}_V \hookrightarrow \mathcal{B}_W$.

Let us consider the slice $\mathcal{W}_W^s = \mathcal{W}^s(f_*) \cap \mathcal{B}_W$ of the stable manifold. If $W$ is sufficiently small then it is a codimension-one complex analytic submanifold in $\mathcal{B}_W$ (Lemma 4.17). By the Stable Manifold Theorem, the orbits of $A|\mathcal{W}_W^s$ uniformly exponentially converge to $f_*$. Hence the spectrum of the restriction $A_*$ to the tangent space $E_W^s = T_{f_*}\mathcal{W}_W^s$ is a discrete set in the open unit disk $\mathbb{D}$ accumulating on 0.

Let us consider the quotient linear operator $\bar{A}_* : \mathcal{B}_W/E_W^s \to \mathcal{B}_W/E_W^s$. Being one-dimensional, it is a multiplication operator, $v \mapsto \rho v$. Let us show that $|\rho| > 1$.

By Corollary 6.2, $|\rho| \geq 1$.

If $|\rho| = 1$ then by the Small Orbits Theorem, for any $\gamma > 0$, there is an $f \in \mathcal{B}_W$ such that $A^n f \in \mathcal{B}(f_*, \gamma)$, $n = 0, 1, \ldots$ but the $A$-orbit of $f$ does not exponentially converge to $f_*$. By the Stable Manifold Theorem, $f \notin \mathcal{H}(f_*)$.

But if $\gamma$ is sufficiently small then $A$ is the analytic continuation of $R^N$ to the Banach slice $\mathcal{B}_W$, and the orbit $\{A^n f\}_{n=0}^\infty$ is nonescaping. By the Combinatorial Rigidity Theorem, $f \in \mathcal{H}(f_*)$ – a contradiction.

Thus $|\rho| > 1$. As the rest of the spectrum of $A$ belongs to the unit disk, $A$ has an eigenvector $h \in \mathcal{B}_W \setminus E_W^s$ corresponding to $\rho$. Let us show that this is also an eigenvector for $R_* : \mathcal{B}_* \to \mathcal{B}_*$ (corresponding to a root $\lambda_* = \rho^{1/N}$). Indeed, let us consider two vectors $h$ and $R_* h = \lambda_* h + w$, where $w \in E^s$. Since both of them are eigenvectors of $R_*^N$ corresponding to the same eigenvalue $\rho$, $R_*^N w = \rho w$ as well. Let us select a Banach slice $E_U^s$ containing $w$. If $w \neq 0$ then

$$\|R_*^{Nm} w\|_U = |\rho|^m \|w\|_U \to \infty$$

contradicting the fact that the orbits of $R_*|E^s$ converge to 0. The theorem is proved. $\square$

6.3. *Unstable manifold.* We keep considering a renormalization operator $R = R_M$ of stationary type near its fixed point $f_*$. As above, $\mathcal{B}_* = E^s \oplus E^u$ means the hyperbolic splitting constructed in the previous section, and $\lambda_*$ stands for the unstable eigenvalue.

THEOREM 6.4 (Local unstable manifold). *There exists a complex analytic one-dimensional manifold $\mathcal{W}_{\mathrm{loc}}^u(f_*) \equiv \mathcal{W}_*^u \subset \mathcal{QG}$ via $f_*$ satisfying the following properties*:

(i) $\mathcal{W}_*^u$ *belongs to some Banach slice $\mathcal{B}_W \ni f_*$, is tangent to $E^u$, and transverse to $\mathcal{W}_*^s$;*

(ii) $\mathcal{W}_*^u \Subset R\mathcal{W}_*^u$ *and the inverse map* $R^{-1} : \mathcal{W}_*^u \to \mathcal{W}_*^u$ *is well-defined;*

(iii) *For any* $f \in \mathcal{W}_*^u$, $\|R^{-n}f - f_*\|_W \sim C_f \lambda_*^{-n}$.

*Proof.* In the proof of the Hyperbolicity Theorem we have constructed a Banach slice $\mathcal{B}_W \ni f_*$ locally invariant under some renormalization iterate $A = R^N$. Moreover, we have proved that this Banach operator is hyperbolic. By the standard hyperbolicity theory, $A$ has a local unstable manifold $\mathcal{W}_*^u$ satisfying properties (i)–(iii) (with $R$ replaced by $A$).

Let us consider the image $\mathcal{S} = R\mathcal{W}_*^u$. Since $R_* : E^u \to E^u$ is a nonsingular operator, $\mathcal{S}$ is a complex one-dimensional manifold tangent to $E^u$ (provided $\mathcal{W}_*^u$ was taken small enough). Moreover, $\mathcal{S}$ sits in some Banach slice $\mathcal{B}_U$, $U \subset W$, and $R^N \mathcal{S} \supset \mathcal{S}$. But $\mathcal{B}_U$ is locally invariant under some iterate $R^{lN}$ whose restriction to $\mathcal{B}_U$ is compact (by the same argument as was used for construction of the operator $A$).

Thus in $\mathcal{B}_U$ we have two analytic submanifolds, $\mathcal{W}_*^s$ and $\mathcal{S}$, tangent to $E^u$ and expanded by the compact hyperbolic analytic map $R^{Nl}$. By the standard hyperbolic theory, these submanifolds must represent the same germ at $f_*$. In other words, the germ of $\mathcal{W}_*^s$ is $R$-invariant. Since $f_*$ is a repelling fixed point for the local restriction $R_*|\mathcal{W}_*^u$, a small disk around $f_*$ in this manifold is strictly expanded under $R$. □

Let us now globalize the unstable manifold. Let us define the *unstable Mandelbrot set* $\mathcal{M}^u$ as the set of infinitely anti-renormalizable points $f \in \mathcal{C}$ such that $R^{-n}f \to f_*$.

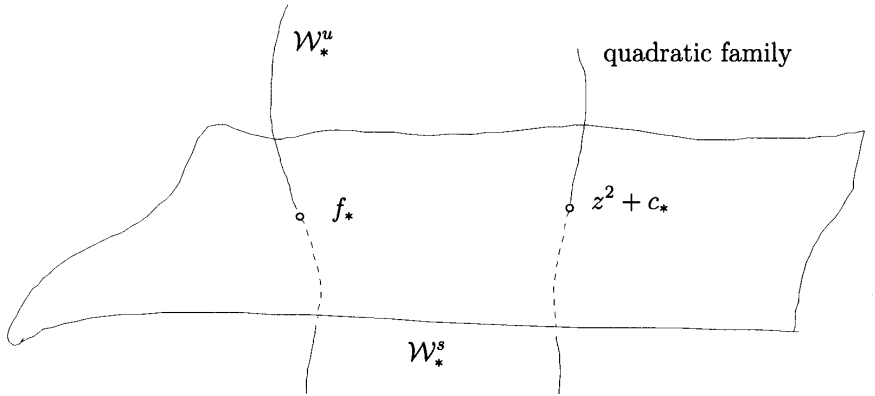

Figure 2. The hyperbolic fixed point of the renormalization operator.

THEOREM 6.5 (Global unstable manifold).

(i) *A point $f \in \mathcal{C}$ belongs to $\mathcal{M}^u(g)$ if and only if there exists a one-sided tower $\mathbf{f} = \{f_0, f_{-1}, \ldots\}$ with stationary combinatorics $\{M, M, \ldots\}$ and a priori bounds. Moreover, in this case $f_{-n} \in \mathcal{W}^u_{\mathrm{loc}}(f_*)$ for all sufficiently big $n$.*

(ii) *The straightening $\mathcal{M}^u \to M_0$ is injective.*

(iii) *For any $\mu > 0$, the set*
$$\mathcal{M}^u_\mu = \{f \in M^u : \exists \text{ a tower } \mathbf{f} \text{ with } f = f_0 \text{ and } \mathrm{mod}(\mathbf{f}) \geq \mu\}$$
*is embedded into a one-dimensional complex analytic manifold $\mathcal{W}^u_\mu(f_*)$ which extends the local manifold $\mathcal{W}^u_{\mathrm{loc}}(f_*)$.*

(iv) *The manifold $\mathcal{W}^u_\mu(f_*)$ is transverse to the foliation $\mathcal{F}$.*

(v) *The germ of the manifold $\mathcal{W}^u_\mu(f_*)$ near $\mathcal{C}$ is invariant under the renormalization.*

*Proof.* (i) Assume that we have an infinitely anti-renormalizable map $f \in \mathcal{C}$ such that $R^{-n}f \to f_*$. Then the germs $f_{-n} = R^{-n}f$ obviously form a one-sided tower with *a priori* bounds.

Vice versa, let $\{f_{-n}\}$ be a one-sided tower with *a priori* bounds. By the Compactness Lemma 4.1, the sequence $\{f_{-n}\}$ is pre-compact. Let us consider its limit set $\Omega$. It is compact and $R$-invariant. Moreover, the map $R : \Omega \to \Omega$ is surjective, since if $f_{-n_k} \to g$, then any limit point of the sequence $\{f_{-n_k-1}\}$ is a preimage of $g$. Hence any point $g \in \Omega$ is included into a two-sided tower with *a priori* bounds. By the Tower Rigidity Theorem, $g$ must coincide with the fixed point $f_*$.

Thus $f_{-n} \to f_*$, and hence the $f_{-n}$ eventually belong to some Banach slice $\mathcal{B}_V \ni f_*$. This Banach slice can be selected so that some iterate $R^N$ keeps it invariant and is hyperbolic at $f_*$. Then $f_{-n}$ must eventually belong to the local unstable manifold of this operator, which coincides (by the Local Unstable Manifold Theorem) with $\mathcal{W}^u_* \equiv \mathcal{W}^u_{\mathrm{loc}}(f_*)$.

(ii) Let us have two maps $f$ and $g$ in $M^u$ which are hybrid equivalent, with $\chi(f) = \chi(g) = c \in M_0$. Then $f_{-n} = R^{-n}f \to f_*$ and $g_{-n} = R^{-n}g \to f_*$. Let $\sigma : M \to M_0$ denote the homeomorphism between the little and big Mandelbrot sets corresponding to $R$. Then the hybrid class of $f_{-n}$ and $g_{-n}$ is $\sigma^{-n}c$, so that $f_{-n}$ is hybrid equivalent to $g_{-n}$. Moreover, the hybrid conjugacy can be selected with uniformly bounded dilatation, since the $\mathrm{mod}(f_{-n})$ and $\mathrm{mod}(g_{-n})$ stay away from 0. Hence the corresponding towers $\mathbf{f} = \{f_{-n}\}$ and $\mathbf{g} = \{g_{-n}\}$ are qc equivalent. By the Tower Rigidity Theorem (ii), these towers coincide up to rescaling, so that $f = g$.

(iii) Note that for $f \in \mathcal{M}^u_\mu$, convergence of the corresponding backward orbit to $f_*$ is uniform: there exists an $N = N(\mu)$ such that $R^{-N}\mathcal{M}^u_\mu \subset \mathcal{W}^u_*$. Indeed, otherwise by compactness one can easily construct a bi-infinite tower with *a priori* bounds which is different from the stationary tower $\{\ldots, f_*, f_*, \ldots\}$.

Let $Q$ stand for the set of $N$ times renormalizable maps in $\mathcal{M}^u_\mu \cap \mathcal{W}^u_*$. Then $R^N Q \supset \mathcal{M}^u_\mu$. Since $R^N$ is transversally nonsingular (by Lemma 5.3), there is a neighborhood $U$ of $Q$ in $\mathcal{W}^u_*$ which is injectively mapped by $R^N$ onto its image. This image is a desired manifold.

(iv) Transversality follows from (ii) and Lemma 4.25.

(v) Invariance follows from the corresponding statement for the local unstable manifold. □

6.4. *Real combinatorics.* Let us now summarize the above information for the case of real combinatorics:

THEOREM 6.6. *Let $M \in \mathcal{N}$ be a real Mandelbrot set and $R = R_M$ be the corresponding renormalization operator. Then*:

(i) *There exists a unique quadratic-like map $f_*$ such that $Rf_* = f_*$; this map is real.*

(ii) *The renormalization operator $R$ is hyperbolic at $f_*$.*

(iii) *The stable manifold $\mathcal{W}^s(f_*)$ coincides with the hybrid class $\mathcal{H}(f_*)$; codim $\mathcal{W}^s(f_*) = 1$.*

(iv) $\dim \mathcal{W}^u_\mu(f_*) = 1$ *and the unstable eigenvalue $\lambda_*$ is positive.*

(v) *For any $\delta > 0$ there exists a $\mu > 0$ such that the unstable manifold $\mathcal{W}^u_\mu(f_*)$ transversally passes through all real hybrid classes $\mathcal{H}_c$ with $c \in [-2, 1/4 - \delta]$.*

*Proof.* There are three points specific for the real situation as compared with the previous complex setting:

a) The complex bounds are established for real maps (see §5.6), so the statement on existence of the fixed point becomes unconditional.

b) Since $\dim \mathcal{W}^u(f_*) = 1$, the unstable eigenvalue $\lambda_*$ is real. To see that it is positive, we need to check that $R$ preserves the orientation of $\mathcal{W}^u(f_*)$ near $f_*$. Sliding to the quadratic family, we see that it would follow from the the property that $\sigma : I \to M_0$ is orientation-preserving, where $I = M \cap \mathbb{R}$. But this property is true by the monotonicity of the kneading invariant in the quadratic family [MT].

c) The bounds of [LS], [LY] are valid for all finitely renormalizable maps such that the last renormalization is not close to the cusp. More precisely, there is a $\mu = \mu(\delta)$ such that for any $N$ and $c \in [-2, 1/4 - \delta]$, there is an $N$ times renormalizable quadratic map $g_N \equiv P_{c_N}$ with $\chi(R^N g_N) = c$ and $\mathrm{mod}(R^k g_N) \geq \mu$, $k = 0, 1, \ldots, N$.

Let us consider finite towers $\mathbf{f}_N = \{f_{N,-k}\}_{k=0}^N$ with $f_{N,-k} = R^{N-k} g_N$. Since they have uniform *a priori* bounds, we can pass to a limit tower $\mathbf{f} = \{f_{-k}\}_{k=0}^\infty$. Then $\mathrm{mod}(Bf) \geq \mu$ and $\chi(f_0) = c$. Hence $\mathcal{M}_\mu^u(f_*)$ passes through all hybrid classes $c$ as above. By the Global Unstable Manifold Theorem (iii–iv), we obtain the last statement of the theorem. □

*Remark.* In particular, the real unstable manifold $\mathcal{W}^u(f_*)$ corresponding to the limit $c_* = \lim c_n$ of the period doublings stretches all way through all real combinatorial types, except "1/4", and is transverse to the bifurcation loci $\mathcal{H}(c_n)$. (And a similar statement can be made for other combinatorics.) This was a part of the Renormalization Conjecture (see [La1]) which previously was established, in the quadratic period-doubling case, with the help of computers, and by Eckmann and Wittwer [EW].

## 7. Hairiness, self-similarity and universality (the stationary case)

7.1. *Proof of the Hairiness Conjecture.* Milnor's Hairiness Conjecture asserts that the Mandelbrot set is becoming dense in small scales near Feigenbaum-like points. Our goal now is to prove this conjecture for stationary combinatorics.

THEOREM 7.1 (Hairiness of the Mandelbrot set). *Let $P_{c_*}$ be a Feigenbaum quadratic polynomial with a priori bounds. Then the magnifications of the Mandelbrot set near $c_*$ converge in the Hausdorff metric on compact sets to the whole complex plane. In particular, this is true for real Feigenbaum points $c_*$.*

*Proof.* Note that by Corollary 10.3, the hairiness property is qc invariant. As the foliation $\mathcal{F}$ is transversally quasi-conformal (Theorem 4.19) and the quadratic family is transverse to $\mathcal{F}$ (Theorem 4.11), it is enough to prove the hairiness property for any transversal. Our choice will be the renormalization unstable manifold.

Let $R$ be the renormalization operator corresponding to the map $P_{c_*}$. Since $P_{c_*}$ has *a priori* bounds, by the results of the previous section, $R$ has a fixed point $f_*$, and is hyperbolic at this point. Let us consider the unstable manifold $\mathcal{W}^u = \mathcal{W}_{\mathrm{loc}}^u(f_*)$ at this point, and the corresponding Mandelbrot set $\mathcal{M}^u = \mathcal{W}^u \cap \mathcal{C}$. By $R^{-1}$ we will denote the branch of the inverse map which maps $\mathcal{W}^u$ into itself. Given a set $\mathcal{X} \subset \mathcal{W}^u$, $R\mathcal{X}$ will mean $R(\mathcal{X} \cap R^{-1}\mathcal{W}^u)$.

Note first that $R\mathcal{M}^u \supset \mathcal{M}^u$. Indeed, if the Julia set of the renormalization $Rf$ is connected then the Julia set of $f$ is connected as well. Hence the Hairiness

Conjecture on $\mathcal{W}^u$ amounts to the following statement:

$$\text{(7.1)} \qquad \bigcup_{m=0}^{\infty} R^m \mathcal{M}^u \ \text{ is dense in } \mathcal{W}^u.$$

On $\mathcal{W}^u$ there is a linearizing coordinate which turns the map $R$ into rescaling by the Feigenbaum universal constant $\lambda_*$. The dist and diam below refer to the distance and the diameter on $\mathcal{W}^u$ in the linearizing coordinate. Accordingly, a "round disk" $D(f, \rho)$ is understood in this sense. Fix some $\varepsilon \in (0, 1)$.

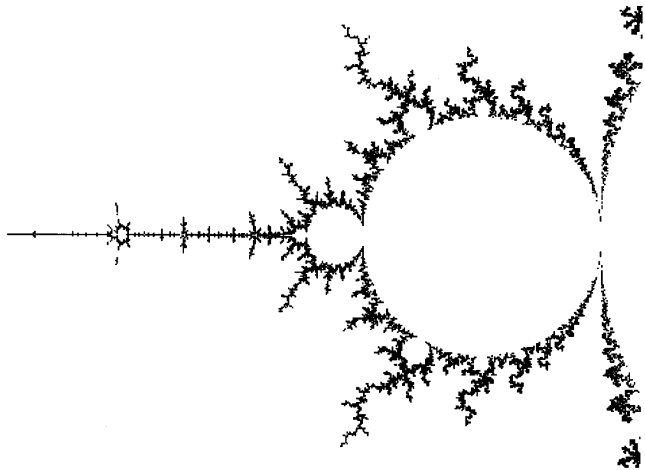

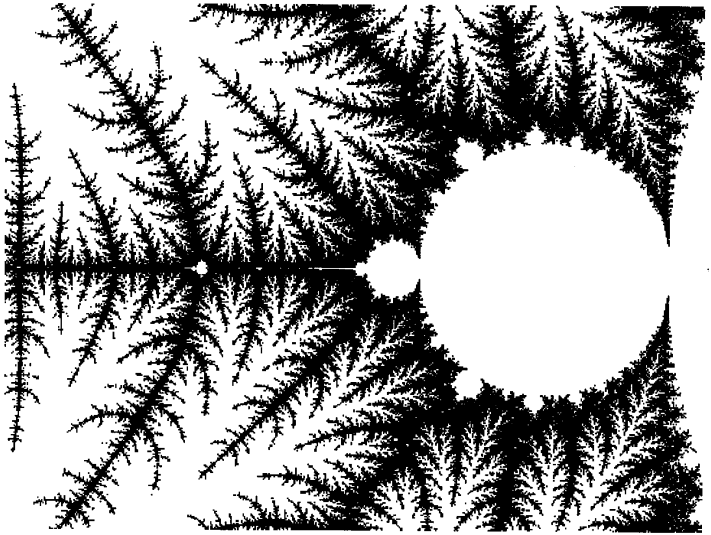

Figure 3. Blow-up of the Mandelbrot set near the Feigenbaum point

Assume that (7.1) fails. Then there exists a round disk

$$\mathcal{V}' = D(f_0, \rho) \Subset \mathcal{W}^u \setminus \bigcup_{m \geq 0} R^m \mathcal{M}^u.$$

Hence for any integer $n \leq 0$,

$$\mathcal{V}'_n \equiv R^n \mathcal{V}' = D(f_n, \lambda_*^n \rho) \Subset \mathcal{W}^u \setminus \bigcup_{m \geq 0} R^m \mathcal{M}^u,$$

where $f_n = R^n f_0$. Let $\mathcal{V}_n \subset \mathcal{V}'_n$ denote the round disks as above of radius $\varepsilon \lambda_*^n \rho$, and let $\mathcal{M}_*^u$ stand for the connected component of $\mathcal{M}^u$ containing $f_*$. Note that $\mathcal{M}_*^u$ is not a single point since the Mandelbrot set $M_0$ does not have isolated points and the holonomy $(M_0, c_*) \to (\mathcal{M}^u, f_*)$ is a local homeomorphism. Now we have:

$$\mathrm{diam} \mathcal{V}_n / \mathrm{dist}(\mathcal{V}_n, \mathcal{M}_*^u) \geq \mathrm{diam} \mathcal{V}_n / \mathrm{dist}(\mathcal{V}_n, f_*) \asymp \varepsilon.$$

Hence

$$\mathrm{diam}_{\mathrm{hyp}}(\mathcal{V}_n | \mathcal{W}^u \setminus \mathcal{M}_*^u) \geq a\varepsilon,$$

where $\mathrm{diam}_{\mathrm{hyp}}(\cdot | \mathcal{U})$ stands for the hyperbolic distance in an open set $\mathcal{U}$, and the constant $a$ is independent of $n$ and $\varepsilon$. It follows that

$$a\varepsilon \leq \mathrm{diam}_{\mathrm{hyp}}(\mathcal{V}_n | \mathcal{W}^u \setminus \mathcal{M}_*^u) \leq \mathrm{diam}_{\mathrm{hyp}}(\mathcal{V}_n | \mathcal{W}^u \setminus \mathcal{M}^u)$$
$$\leq \mathrm{diam}_{\mathrm{hyp}}(\mathcal{V}_n | \mathcal{V}'_n) \asymp \varepsilon.$$

Thus

(7.2) $$\mathrm{diam}_{\mathrm{hyp}}(\mathcal{V}_n | \mathcal{W}^u \setminus \mathcal{M}^u) \asymp \varepsilon$$

with the constant independent of $n$ and $\varepsilon$.

Let us show that if $\varepsilon$ is sufficiently small then over every domain $\mathcal{V}_n$ there is a holomorphic motion $h_F : U' \to \tilde{U}'$ which conjugates $F_0 \equiv f_n : U_0 \to U'_0$ to $F : U_F \to U'_F$, $F \in \mathcal{V}_{-n}$. Here the domains $U_F$ and $U'_F$ will be selected in such a way that the $\mathrm{mod}(U'_F \setminus U_F)$ stay away from 0 (as $n$ changes).

If $\mathcal{W}^u$ is selected sufficiently small then $F(0) \neq 0$ for $F \in \mathcal{W}^u$, and thus $F \in \mathcal{W}^u$ can be normalized so that $F(0) = 1$. Also, by Lemma 4.2, there is a fundamental annulus $A_F = h_F A_0$ holomorphically moving with $F \in \mathcal{W}^u$. Let us consider two preimages of this annulus, $A_F^1 = F^{-1} A_F$ and $A_F^2 = F^{-2} A_F$. If $F$ is sufficiently close to $f_*$ (so that $F^n(0) \notin A_F$ for $n = 1, 2$) then these preimages are holomorphically moving annuli. Moreover, if we select an equipotential foliation in $A_0$, we obtain the holomorphically moving equipotential foliation $\gamma_F(r) = h_F \gamma_0(r)$ in the union of these three annuli. Let us assign the level $r$ to the equipotentials in such a way that the outer boundary of $A_F$ has level 4, the inner boundary has level 2, and $F(\gamma_F(r)) = \gamma_F(r^2)$ for $2^{1/4} \leq r \leq 2$. Note that by our normalization, the union of these three annuli belongs to the twice punctured plane $\mathbb{C} \setminus \{0, 1\}$.

Take a point $z_0 \in A_0^1$. Let it belong to an equipotential $\gamma_0(r)$, $\sqrt{2} \leq r \leq 2$. Let $Q_0 \equiv Q_0(z)$ denote the fundamental annulus bounded by equipotentials $\gamma_0(r^{4/3})$ and $\gamma_0(r^{2/3})$, and let $Q_F = h_F Q_0$ stand for its motion. Let $d > 0$ be the hyperbolic distance (in $\mathbb{C} \setminus \{0,1\}$) from $z_0$ to the boundary of $Q_0$.

We find a moment $l$ such that $F^{l+1}(0) \in A_0$, and take $z_0 = F^l(0)$. Since the map $\mathcal{W}^u \setminus \mathcal{M}^u \to \mathbb{C} \setminus \{0,1\}$, $F \mapsto F^l(0)$, contracts the hyperbolic metric, the hyperbolic distance from $F^l(0)$ to $F_0^l(0)$ in $\mathbb{C} \setminus \{0,1\}$ is at most $A\varepsilon$, provided $F \in \mathcal{V}_n$. For the same reason, for $\zeta \in \partial Q_0$, the hyperbolic distance from $\zeta$ to $h_F \zeta$ is at most $A\varepsilon$ for $F \in \mathcal{V}_n$. It follows that $h_F \zeta \neq F^l(0)$ for $F \in \mathcal{V}_n$, provided $\varepsilon < d/2A$.

Thus for sufficiently small $\varepsilon$, the point $z_F = F^l(0)$, $F \in \mathcal{V}_n$, does not cross the boundary of the fundamental annulus $Q_F$. Hence we have a holomorphic motion of $(\partial Q_F, z_F)$ over $\mathcal{V}_n$. By the $\lambda$-lemma (see Appendix 2), this motion extends to the motion of the whole annulus $(Q_F, z_F)$. Pulling it back by dynamics, we obtain the motion $U_F' \setminus K(F)$ over $\mathcal{V}_n$, where $U_F'$ is the domain bounded by the equipotential $\gamma_F(r^{4/3})$. Using the $\lambda$-lemma again, we extend this motion through the Julia set, which provides us with the desired motion $h_F$ of the domain $U_F'$ over $\mathcal{V}_n$ which conjugates $F_0$ to $F$.

By the $\lambda$-lemma, $h_F$ is $K$-quasi-conformal over the twice smaller disk $\Delta_n = D(f_n, \varepsilon \lambda_*^n \rho/2)$, with an absolute $K$. Take any point $\tilde{f}_0 \in \Delta_0$, and let $\tilde{f}_n = R^n \tilde{f}_0$. We conclude that $\tilde{f}_n$ is $K$-qc conjugate to $f_n$, with an absolute $K$.

Let us now consider the towers (with disconnected Julia sets)
$$\mathbf{f}_m = \{f_{n-m}\}_{n=m}^{-\infty} \quad \text{and} \quad \tilde{\mathbf{f}}_m = \{\tilde{f}_{n-m}\}_{n=m}^{-\infty}$$
(so that $f_{-m}$ is the zero coordinate of $\mathbf{f}_m$). They both converge to the stationary tower $\mathbf{f}_* = \{\ldots, f_*, f_*, \ldots\}$. Moreover, the above $K$-qc conjugacies $h_n$ converge to a qc automorphism $h$ of $\mathbf{f}_*$. By Corollary 5.12,

(7.3) $$h = \text{id}.$$

Let us show that on the other hand, the $h_n$ stay a definite distance away from the id. To this end let us now pass from the unstable manifold $\mathcal{W}^u$ to the vertical fiber $\mathcal{Z} \equiv \mathcal{Z}_{f_*}$. Let $\mathcal{M}^v \subset \mathcal{Z} = \mathcal{Z} \cap \mathcal{C}$ denote the Mandelbrot set in $\mathcal{Z}$. By Theorem 4.18 the foliation $\mathcal{F}$ extends to some Banach neighborhood of $f_*$, and by the $\lambda$-lemma (see Appendix 2) this foliation is transversally quasi-conformal. Moreover, both $\mathcal{W}^u$ and $\mathcal{Z}$ are transverse to the foliation at $f_*$ (by Corollary 4.24, Lemma 4.25 and the Unstable Manifold Theorem). Hence the holonomy from $\mathcal{W}^u$ to $\mathcal{Z}$ is well-defined and quasi-conformal near $f_*$.

Take two points $f_n = R^n f_0 \in \mathcal{V}_n$ and $\tilde{f}_n = R^n \tilde{f}_0 \in \mathcal{V}_n$, $n \leq 0$. Let $\phi_n, \tilde{\phi}_n \in \mathcal{Z}$ correspond to $f_n$ and $\tilde{f}_n$ under the holonomy (they are well-defined for sufficiently big $n \leq 0$). Then by (7.2) and quasi-invariance of the hyperbolic metric (see [LV, Ch. II, §3.3]),

(7.4) $$\text{dist}_{\text{hyp}}(\phi_n, \tilde{\phi}_n | \mathcal{Z} \setminus \mathcal{M}^v) \asymp 1.$$

Take a little disk $\mathcal{D}$ in $\mathcal{Z}$ around $f_*$, and consider a curve $\Gamma = \partial \mathcal{D} \setminus \mathcal{M}^v$ (perhaps disconnected). Then obviously

$$(7.5) \qquad \operatorname{dist}_{\operatorname{hyp}}(\phi_n, \Gamma|\, \mathcal{Z} \setminus \mathcal{M}^u) \to \infty \text{ as } n \to -\infty.$$

Let us consider the map $\xi : \mathcal{Z} \setminus \mathcal{M}^v \to \mathbb{C} \setminus \bar{\mathbb{D}}$ which assigns to $\phi$ the position of the critical value in the external model (see §3.3). By Lemma 4.14, this map is conformal. Hence by (7.4),

$$(7.6) \qquad \operatorname{dist}_{\operatorname{hyp}}(\xi_n, \tilde{\xi}_n|\, \xi(\mathcal{U}_n)) \asymp 1,$$

where $\xi_n = \xi(\phi_n)$, $\tilde{\xi}_n = \xi(\tilde{\phi}_n)$, and $\mathcal{U}_n$ is the component of $\mathcal{Z} \setminus \mathcal{M}^v$ containing $\phi_n$ and $\tilde{\phi}_n$.

Let us consider in $\mathcal{U}_n$ the domain $\hat{\mathcal{U}}_n$ cut off by the arc $\Gamma_n = \Gamma \cap \mathcal{U}_n$ and containing $\phi_n$ and $\tilde{\phi}_n$. The image $\xi(\hat{\mathcal{U}}_n)$ is bounded by the arc $\xi(\Gamma_n)$ and an arc of the unit circle. By (7.4) and (7.5),

$$(7.7) \qquad \operatorname{dist}_{\operatorname{hyp}}(\xi_n, \tilde{\xi}_n|\, \mathbb{C} \setminus \bar{\mathbb{D}}) \asymp 1.$$

We refer now to the proof of Theorem 4.18 which extends the foliation $\mathcal{F}$ beyond the connectedness locus. Consider a quadratic-like map $G_* = \Pi(f_*) \in \mathcal{H}_0$ and the corresponding Riemann mapping $R_* : \mathbb{C} \setminus K(G_*) \to \mathbb{C} \setminus \bar{\mathbb{D}}$. Let us transfer the $\xi$-points by this Riemann mapping as prescribed by (4.14):

$$a_n = R_*^{-1} \xi_n, \quad \tilde{a}_n = R_*^{-1} \tilde{\xi}_n.$$

By (7.7)

$$(7.8) \qquad \operatorname{dist}_{\operatorname{hyp}}(a_n, \tilde{a}_n|\, \mathbb{C} \setminus K(G_*)) \asymp 1.$$

Let us now consider the tubing (4.15) near $f_*$, and transfer the $a$-points according to (4.16):

$$b_n = \Psi_*^{-1} a_n, \quad \tilde{b}_n = \Psi_*^{-1} \tilde{a}_n,$$

where $\Psi_* \equiv \Psi_{G_*}$. Since $\Psi_*$ is quasi-conformal,

$$\operatorname{dist}_{\operatorname{hyp}}(b_n, \tilde{b}_n|\, \mathbb{C} \setminus \bar{\mathbb{D}}) \asymp 1,$$

so that

$$\operatorname{dist}(b_n, \tilde{b}_n) \asymp \operatorname{dist}(b_n, \mathbb{T}) \asymp \operatorname{dist}(\tilde{b}_n, \mathbb{T}).$$

Fix a small $\delta > 0$. Then applying $P_0 : z \mapsto z^2$ to $b_n$ and $\tilde{b}_n$ an appropriate number $l_n$ of times, we can move the above $b$-points a distance of order $\delta$ apart (independently of $n$):

$$(7.9) \qquad \operatorname{dist}(P_0^{l_n} b_n,\, P_0^{l_n} \tilde{b}_n) \asymp \operatorname{dist}(P_0^{l_n} b_n, \mathbb{T}) \asymp \operatorname{dist}(P_0^{l_n} \tilde{b}_n, \mathbb{T}) \asymp \delta.$$

Let us now return to the unstable manifold $\mathcal{W}^u$. Denote the maps and the points corresponding to $f_n$ as follows:

$$G_n = \Pi(f_n) \in \mathcal{H}_0;$$

$S_n : \mathbb{C}\setminus\bar{\mathbb{D}} \to \mathbb{C}\setminus K(G_n)$ is the Riemann mapping with a positive derivative at $\infty$;

$\Psi_n = \Psi_{G_n}$ is the tubing map;

$\eta_n = \xi(f_n) \in \mathbb{C} \setminus \bar{\mathbb{D}}$ is the position of the critical value in the external model;

$d_n = S_n \eta_n$ is the marked points associated with $G_n$

(and denote correspondingly the tilde-objects).

As $f_n$ and $\phi_n$ (respectively $\tilde{f}_n$ and $\tilde{\phi}_n$) lie on the same leaf of the extended foliation, by (4.16):

$$d_n = \Psi_n b_n; \quad \tilde{d}_n = \tilde{\Psi}_n \tilde{b}_n.$$

Since $\Psi_G$ continuously depends on $G$ near $G_*$, $\Psi_n$ and $\tilde{\Psi}_n$ are uniformly close. Hence (7.9) implies that there exist some $\delta'' > \delta' > 0$ depending only on $\delta$ such that

(7.10) $$\delta' < \operatorname{dist}(G_n^{l_n} d_n, \tilde{G}_n^{l_n} \tilde{d}_n) < \delta'',$$

$$\delta' < \operatorname{dist}(G_n^{l_n} d_n, \mathbb{T}) < \delta'', \quad \delta' < \operatorname{dist}(\tilde{G}_n^{l_n} \tilde{d}_n, \mathbb{T}) < \delta''.$$

Let us now transfer the $d$-points to the external model: $\eta_n = S_n^{-1}(d_n)$. By Lemma 4.15, the Riemann mappings $S_n$ converge to the Riemann mapping $S : \mathbb{C} \setminus \bar{\mathbb{D}} \to \mathbb{C} \setminus K(G_*)$. Hence by (7.10),

(7.11) $$\delta' < \operatorname{dist}(g_n^{l_n} \eta_n, \tilde{g}_n^{l_n} \tilde{\eta}_n) < \delta'',$$

$$\delta' < \operatorname{dist}(g_n^{l_n} \eta_n, \mathbb{T}) < \delta'', \quad \delta' < \operatorname{dist}(\tilde{g}_n^{l_n} \tilde{\eta}_n, \mathbb{T}) < \delta'',$$

where the constants $\delta'$ and $\delta''$ are not the same as in (7.10) but satisfy the same properties.

Finally, let us transfer the above $\eta$-points to the original dynamical plane of maps $f_n$ and $\tilde{f}_n$. They correspond to the points $f^{l_n+1}(0)$ and $\tilde{f}^{l_n+1}(0)$ via conformal maps which are defined outside a small neighborhood of the unit disk and by Lemma 4.15 converge to the Riemann mapping $R_* : \mathbb{C}\setminus\bar{\mathbb{D}} \to \mathbb{C}\setminus J(f_*)$. Hence (7.11) yields:

$$\operatorname{dist}(f_n^{l_n+1}(0), \tilde{f}_n^{l_n+1}(0)) \geq \delta' > 0.$$

Since the conjugacy $h_n$ between $f_n$ and $\tilde{f}_n$ carries $f^{l_n+1}(0)$ to $\tilde{f}_n^{l_n+1}(0)$, it stays uniformly away from the identity, as was asserted.

Hence the limiting qc automorphism $h$ of the stationary tower $\mathbf{f}_*$ cannot be identical, contradicting (7.3). This contradiction proves (7.1) and hence the Hairiness Conjecture. □

7.2. *Self-similarity of the Mandelbrot set.* Below we will prove, in the stationary case, the Self-Similarity Theorem stated in the introduction for real-parameter values. The corresponding complex statement is the following:

THEOREM 7.2 (Self-similarity). *Let $M \in \mathcal{N}$ be a real Mandelbrot copy and $\sigma : M \to M_0$ be the homeomorphism of $M$ onto the whole Mandelbrot set $M_0$. Assume that there exists a quadratic-like map $f$ with stationary combinatorics $\tau(f) = \{\ldots, M, M, \ldots\}$ and a priori bounds. Let $c_* = \chi(f)$. Then $c_*$ is a fixed point of $\sigma$, and $\sigma$ is $C^{1+\alpha}$-conformal at $c_*$, with the derivative at $c_*$ equal to the Feigenbaum universal scaling constant $\lambda_* > 1$. Moreover, there exists at most one parameter value $c_*$ satisfying the above assumptions.*

*Remarks.* 1. The real theorem stated in the introduction follows from this complex one and *a priori* bounds (Theorem 5.6).

2. This theorem does not rule out another fixed point $c \in M$ of $\sigma$ (for which the map $P_c$ fails to have *a priori* bounds). However, it rules out other fixed points near $c_*$.

A map $h : (M_1, 0) \to (M_2, 0)$ between two subsets in $\mathbb{C}$ is called $C^{1+\alpha}$-conformal at the origin if there exist a $\tau \neq 0$ such that $h(u) = \tau u(1 + O(|u|^\alpha))$ for $u \in M_1$ near 0.

The foliation $\mathcal{F}$ is transversally $C^{1+\alpha}$-conformal at a point $c \in M$ (or along a leaf $\mathcal{H}_c$) if for any two transversals $\mathcal{S}$ and $\mathcal{T}$ to the leaf $\mathcal{H}_c$, the holonomy $\mathcal{M}_\mathcal{S} \to \mathcal{M}_\mathcal{T}$ between the corresponding Mandelbrot sets is $C^{1+\alpha}$-conformal at the points of intersection with $\mathcal{H}_c$.

LEMMA 7.3 (Transversal conformality at a Feigenbaum point). *Let $c_* \in M$ be a Feigenbaum parameter value satisfying the assumptions of the self-similarity theorem. Then $\mathcal{F}$ is transversally $C^{1+\alpha}$-conformal at $c_*$ with some $\alpha = \alpha(M) > 0$.*

*Proof.* By the Hyperbolicity Theorem, there exists a renormalization hyperbolic fixed point $f_* \in \mathcal{H}_{c_*}$ with the stable manifold $\mathcal{H}_{c_*}$ and the transverse unstable manifold $\mathcal{W}^u \equiv \mathcal{W}^u_{\text{loc}}(f_*)$. Clearly it is enough to check $C^{1+\alpha}$-conformality of the holonomy from a transversal $\mathcal{S}$ via $f \in \mathcal{H}_*$ to the unstable manifold $\mathcal{W}^u$.

By the Stable Manifold Theorem, there exist a quadratic-like representative $f_* : V_* \to V'_*$ and natural numbers $N = N(V_*, f)$, $l = l(V_*)$ such that $(\mathcal{B}_{V_*}, f_*)$ is locally invariant under $R^l$, $R^N f \in \mathcal{B}_V$ and the orbit of $R^N f$ under $R^l$ exponentially converges to $f_*$ in this Banach slice. Moreover, these properties are still valid (with a different $l$) if we take any other representative $f_* : V \to V'$ with $V \subset V_*$.

Let us take a Banach slice $\mathcal{B}_U \ni f$ locally containing the transversal $\mathcal{S}$. Then there is a neighborhood $\mathcal{U} \subset \mathcal{B}_U$ of $f$ which is mapped by $R^N$ into some Banach slice $\mathcal{B}_V$ as above (since for nearby $F \in \mathcal{U}$, the renormalization $R^N F$ is well-defined on any domain $V \subset V_*$). Thus the curve $\mathcal{S}^N = R^N \mathcal{S}$ locally sits in $\mathcal{B}_V$.

But it is sufficient to study the holonomy $h_N$ from $\mathcal{S}^N$ to $\mathcal{W}^u$. Indeed if $h$ denotes the holonomy from $\mathcal{S}$ to $\mathcal{W}^u$ then by the $R$-invariance of the foliation $\mathcal{F}$, $h = h_N \circ R^N$ where $R^N : (\mathcal{S}, f) \to (\mathcal{S}^N, R^N f)$ is a local conformal diffeomorphism (by Lemma 5.3).

Thus the situation reduces to the Banach set up. Let $T = R^l$. Without loss of generality we can assume that $\mathcal{S}$ itself belongs to a Banach neighborhood $\mathcal{V} \subset \mathcal{B}_V$ of $f_*$. Moreover, let us select this neighborhood $\mathcal{V}$ as a box $E^s(\delta) \times E^u(\delta)$, where $E^{s/u}(\delta)$ is the $\delta$-ball in the tangent space $E^{s/u} = \mathrm{T}_{f_*}\mathcal{W}_V^{s/u}$. If $\delta$ is sufficiently small then $T\mathcal{V} \subset \mathcal{B}_V$ and $T$ is hyperbolic on $\mathcal{V}$ in the sense that it satisfies properties H1, H2 of Lemma 2.1 plus the analogous vertical expansion property. Recall also that by Lemma 4.17, the Banach slice $\mathcal{F}_V$ of the foliation $\mathcal{F}$ is a foliation near $f_*$ whose leaves are graphs over $E^s$.

Let $\mathcal{S}^n$ stand for the *truncated iterate* of $\mathcal{S}$, i.e., inductively, let $\mathcal{S}^{n+1} = \mathcal{V} \cap T\mathcal{S}^n$. Then eventually the $\mathcal{S}^n$ can be represented as graphs of analytic functions $\phi_n : E^u(\delta) \to E^s(\delta)$, so that we can assume that this happens from the very beginning. The local unstable manifold $\mathcal{W}^u \cap \mathcal{V}$ can certainly also be parametrized in the same way by some function $\psi$. By the hyperbolicity of $T$ on $\mathcal{V}$, the manifolds $\mathcal{S}^n$, exponentially fast, converge to the unstable manifold:

$$\|\phi_n - \psi\|_{C^1} \leq \kappa \gamma^n, \tag{7.12}$$

where $\gamma \in (0, 1)$ is a strict upper bound on the spectrum of $\mathrm{D}T(g_*)$ lying inside the unit disk. Moreover, $\kappa > 0$ can be *a priori* selected arbitrarily small (just replace $\mathcal{S}$ by some $\mathcal{S}^m$ with $m = m(\kappa)$).

Let us use the projections $p : \mathcal{S}^m \to E^u$ as analytic charts on $\mathcal{S}^m$. By the Koebe theorem, they have distortion $O(\varepsilon)$ in scale $\varepsilon$ with a uniform constant (independent of $m$). To simplify the notation, we will skip $p$, so that for $u, v \in \mathcal{S}^m$, $u - v$ means the difference between the local coordinates: $p(u) - p(v)$.

Let $\mathcal{M} = \mathcal{S} \cap \mathcal{C}$ and $\mathcal{M}^m = \mathcal{S}^m \cap \mathcal{C}$ be the truncated iterates of $\mathcal{M}$.

Select $q > \mu > 1$. Take two points $z_1, z_2 \in \mathcal{M}$ on distances of order $\varepsilon = q^{-n}$ from $a \equiv f$. Push them forward by $T^m$ so that they go to points $\zeta_1, \zeta_2, b = T^m a \in \mathcal{M}^m$ with relative distances of order $\mu^{-n}$. By the Koebe Distortion Theorem, the ratio distortion of $T^m$ at the above three points is of order $\mu^{-n}$:

$$\frac{\zeta_2 - b}{\zeta_1 - b} = \frac{z_2 - a}{z_1 - a}(1 + O(\mu^{-n})). \tag{7.13}$$

Furthermore, $m$ is at least $n \log(q/\mu)/\log \lambda \equiv cn$, where $\lambda$ is an upper bound for the unstable eigenvalue $\lambda_*$ (and $\lambda = \lambda(\delta, \kappa)$ can be made arbitrarily close to $\lambda_*$ by choosing the parameters $\delta$ and $\kappa$ sufficiently small). By (7.12), $\mathcal{S}^m$ is a distance $O(\gamma^m) = O(\rho^n)$ from $\mathcal{W}^u$, where $\rho = \gamma^c < 1$.

It follows that the holonomy from $\mathcal{S}^m$ to $\mathcal{W}^u$ has an exponentially small ratio distortion at $\zeta_1, \zeta_2, b$. Indeed, let us extend the foliation $\mathcal{F}_V$ to a neighborhood of $f_*$ in $\mathcal{B}_V$ (see Theorem 4.18). If $\delta$ is sufficiently small then the bi-disk $E^s(2\delta) \times E^u(2\delta)$ is contained in the domain of the extended foliation. Consider the holonomy $h_m : \mathcal{S}^m \to \mathcal{W}^u$ along the extended foliation. By the $\lambda$-lemma (see Appendix 2), $h_m$ is $K_m$-qc with $K_m = 1 + O(\rho^m)$. Then by the distortion estimates for qc maps (see [LV, Ch. II, Thm. 3.1])

$$\text{(7.14)} \qquad \frac{|h_m\zeta_2 - c_*|}{|h_m\zeta_1 - c_*|} = \frac{|\zeta_2 - b|}{|\zeta_1 - b|}(1 + O(\rho_1^m)),$$

with any $\rho_1 > \rho$.

If the distance between $\zeta_1$ and $\zeta_2$ is commensurable with their distance to $b$ then the same estimate holds for the other two ratios (centered at $\zeta_1$ and $\zeta_2$). Then by Euclidean trigonometry the angles of the triangle $\Delta(b, \zeta_1, \zeta_2)$ differ from the corresponding angles of its image $\Delta(a, h_m\zeta_1, h_m\zeta_2)$ by $O(\rho_1^n)$. But then it is also true without the assumption that the distance between $\zeta_1$ and $\zeta_2$ is commensurable with their distance to $b$ since a small angle can be represented as a difference of two angles satisfying the commensurability assumption. Thus the holonomy $h_m$ preserves the angles at $b$ up to order $O(\rho_1^n)$. Together with (7.14) this yields:

$$\text{(7.15)} \qquad \frac{h_m\zeta_2 - c_*}{h_m\zeta_1 - c_*} = \frac{\zeta_2 - b}{\zeta_1 - b}(1 + O(\rho_1^m)),$$

Finally, let us apply the inverse map $T^{-m} : \mathcal{W}^u \to \mathcal{W}^u$ to the points $h_m(\zeta_i)$. Since the foliation $\mathcal{F}_V$ is $T$-invariant, we obtain the points $h(z_i)$. Moreover, by the Koebe Theorem, the ratio distortion of this transition is $O(\mu_1^{-n})$ with any $\mu_1 > \mu$. Combining this with (7.13), (7.15), we obtain the ratio distortion estimate for the holonomy $h : \mathcal{S} \to \mathcal{W}^u$

$$\frac{hz_2 - c_*}{hz_1 - c_*} = \frac{z_2 - a}{z_1 - a}(1 + O(\mu_1^{-n} + \rho_1^n)).$$

Thus the ratio distortion of $h$ in scale $\varepsilon > 0$ about $a$ is of order $\varepsilon^\alpha$ with some $\alpha > 0$. This implies $C^{1+\alpha}$-conformality. Indeed, take two points $u, v \in \mathcal{M}$ with $|v - a| \leq |u - a| < \varepsilon$, and let $z_0 \equiv u, z_1, \ldots, z_k \equiv v$, be a string of points in $\mathcal{M}$ such that $|z_i| = |z_{i-1}|/2$ for $i < k$ and $|z_{k-1}|/2 \leq |z_k| < |z_{k-1}|$. (Such a string exists since the Mandelbrot set $M_0$ is connected and the holonomy $M_0 \to \mathcal{M}$ is continuous. Hence the set $\mathcal{M}$ intersects every circle around $a \in \mathcal{S}$ with sufficiently small radius.) Then

$$\frac{h(z_i) - a}{z_i - a} : \frac{h(z_{i+1}) - a}{z_{i+1} - a} = 1 + O(\varepsilon^\alpha/2^{\alpha i}).$$

Hence
$$\frac{h(u) - a}{u - a} : \frac{h(v) - a}{v - a} = 1 + O(\varepsilon^\alpha).$$

It follows that there exists a $\lim_{u \to 0}(h(u) - a)/(u - a) \equiv \tau \neq 0$ and
$$\frac{h(u) - a}{u - a} = \tau(1 + O(|u - a|^\alpha)). \qquad \square$$

*Proof of the self-similarity theorem.* By the Combinatorial Rigidity Theorem 5.9, the point $c_*$ is fixed by $\sigma$ and there is only one point $c_*$ satisfying the assumptions of the theorem.

The holonomy $\sigma : M_0 \to \mathcal{M}^u$ locally conjugates $\sigma$ to the renormalization operator $R$. Since this holonomy is $C^{1+\alpha}$-conformal at $c_*$ and $R|\mathcal{W}^u$ is locally conformal, $\sigma$ is $C^{1+\alpha}$-conformal at $c_*$. $\qquad \square$

7.3. *Universality Theorem.* Let us consider a little Mandelbrot copy $M = M_c \in \mathcal{L}$ with $p(M) = p$, and the corresponding homeomorphism $\sigma : M \to M_0$. Then the "tuned copies" $M^n = \sigma^{-n}M$ are centered at superattracting parameter values $c_n = c^{*n}$ with period $p(M^n) = p^n$. The corresponding polynomials $P_n \equiv P_{c_n}$ are $n$ times renormalizable with $R^n P_n \in \mathcal{H}_0$.

THEOREM 7.4 (Universality). *Assume that the polynomials $P_n$ have a common a priori bound*: $\mathrm{mod}(R^k P_n) \geq \mu > 0, \ n = 1, 2, \ldots, \ k = 0, 1, \ldots, n.$ *Then*:

- *The $c_n$ converge exponentially fast to an infinitely renormalizable parameter value $c_*$*:
$$|c_n - c_*| \sim a\lambda_*^{-n}.$$

- *Let $\mathcal{S} = \{f_\mu\}$ be a complex analytic transversal to the hybrid class $\mathcal{H}_{c_*}$ at some $\mu_*$. Then for $\mu$ near $\mu_*$ and all sufficiently big $n$, $\mathcal{S}$ has a unique intersection point $\mu_n$ with the hybrid class $\mathcal{H}_{c_n}$, and the $\mu_n$ converge to $\mu_*$ with the universal exponential rate*:
$$|\mu_n - \mu_*| \sim a(\mathcal{S})\lambda_*^{-n}.$$

*In particular, this yields the Universality Theorem for real parameter values stated in the introduction.*

*Proof.* Take any accumulation point $c_*$ of the sequence $\{c_n\}$. By the uniform *a priori* bounds assumption and the compactness Lemma 4.1, $P_{c_*}$ is an infinitely renormalizable polynomial with combinatorics $\tau = \{M, M, \ldots\}$ and *a priori* bounds. By the Combinatorial Rigidity Theorem, such a $c_*$ is unique. Hence $c_n \to c_*$.

By the Hyperbolicity Theorem, the renormalization operator $R = R_M$ has a unique fixed point $f_*$, and this point is hyperbolic. Let us consider its unstable manifold $\mathcal{W}^u = \{f_\nu\}$. Since the quadratic family $\mathcal{Q}$ and $\mathcal{W}^u$ are both transverse to the hybrid class $\mathcal{H}_* \equiv \mathcal{H}_{c_*}$, for all sufficiently big $n$, there exists a unique parameter value $\nu_n$ near $\nu_*$ corresponding to $c_n$ under the holonomy. Since the holonomy conjugates $\sigma$ and $R$, we have $Rf_{\nu_n} = f_{\nu_{n-1}}$. As $f_*$ is an expanding fixed point for $R|\mathcal{W}^u$ with eigenvalue $\lambda_*$,

$$|\nu_n - \nu_*| \sim q\lambda_*^{-n}, \tag{7.16}$$

and the result follows from the smoothness of the holonomy along the stable leaf (Lemma 7.3). □

7.4. *Connection to the* MLC *problem.* The problem of local connectivity of the Mandelbrot set (MLC) is a central theme in holomorphic dynamics. By works of Yoccoz (see [H]) and the author [L2], [L5], MLC is now established for all real $c$ except those which are infinitely renormalizable with type bounded by some $\bar{p}$. The following criterion links this problem to the compactness of the Mandelbrot set in the unstable manifold.

PROPOSITION 7.5. *Let $M$ be a primitive little Mandelbrot set. Let $c_* \in M_0$ be an infinitely renormalizable parameter value of type $\{M, M, \ldots\}$ with a priori bounds (for example, a real one). Then the following properties are equivalent*:

(i) *The Mandelbrot set $M_0$ is locally connected at $c_*$*;

(ii) *The unstable Mandelbrot set $\mathcal{M}^u$ of the $R_M$-fixed point $f_*$ is compact*;

(iii) *For any $c \in M_0$, there exists a tower $\mathbf{f}_c = \{\ldots \mapsto f_{-1} \mapsto f_0\}$ of type $\{\ldots, M, M\}$ with $\chi(f_0) = c$ and with a priori bounds.*

*If $M$ is satellite then the same statement is true for $M_0$ replaced with $M$.*

*Proof.* It is known that local connectivity of the Mandelbrot set at $c_*$ is equivalent to shrinking the tuned copies $M^n$ to the point (see e.g., [Sch] for a discussion of this kind of relation).

(i) ⇒ (iii). If the tuned copies $M^n$ shrink to $c_*$, then all of them eventually belong to the domain of the holonomy $h : M^n \to \mathcal{M}^{u,n}$. Hence for any $c \in M^n$, there exists a tower $\mathbf{f}_c = \{h(c) = f_0, R^{-1}f_0, \ldots\}$ with *a priori* bounds. Since $\sigma^n$ maps $M^n$ onto the whole Mandelbrot set and $h \circ \sigma^n = R^n \circ h$ (resp. $h \circ \sigma^{n-1} = R^{n-1} \circ h$ in the satellite case), any tower $\mathbf{f}_c$ with $c \in M_0$ (resp. $c \in M$) and *a priori* bounds is realizable.

(iii) ⇒ (ii). Assume that for any $c \in M_0$ (resp. $c \in M$ in the satellite case), there is a tower $\mathbf{f}_c$ with *a priori* bounds. Then by Theorem 6.5, the holonomy $h : M_0 \to \mathcal{M}^u$ (resp. $M \to \mathcal{M}^{u,1}$) is well-defined, and hence has a compact image.

(ii) $\Rightarrow$ (i). Note that the image of $\mathcal{M}^u$ under the straightening $\chi : \mathcal{M}^u \to M_0$ is open in $M_0$. Indeed, the straightening homeomorphically maps a relative neighborhood $\mathcal{U} \subset \mathcal{M}^u$ of $f_*$ onto a relative neighborhood $U \subset M_0$ of $c_*$ and conjugates $R$ to $\sigma$. Since by definition,

$$(7.17) \qquad \mathcal{M}^u = \bigcup_{n \geq 0} R^n(\mathcal{U} \cap \mathrm{Dom}(R^n)),$$

we conclude:

$$\chi(\mathcal{M}^u) = \bigcup_{n \geq 0} \sigma^n(U \cap M^n).$$

But $\sigma^n : M^n \to M_0$ is a homeomorphism, and openness of $\chi(\mathcal{M}^u)$ follows.

On the other hand, if $\mathcal{M}^u$ is compact then the image $\chi(\mathcal{M}^u)$ is closed. As $M_0$ is connected, $\chi(\mathcal{M}^u) = M_0$. Hence

$$(7.18) \qquad \chi(R^{-n}\mathcal{M}^u) = \sigma^{-n} M_0 \equiv M^n.$$

Furthermore, by (7.17) and compactness of $\mathcal{M}^u$, $R^{-N}\mathcal{M}^u \subset \mathcal{U}$ for some $N$. If $\mathcal{U}$ is contained in $\mathcal{W}^u_{\mathrm{loc}}$ then clearly the $R^{-n}\mathcal{U}$ shrink to $f_*$. Hence the $R^{-n}\mathcal{M}^u$ also shrink to $f_*$. By (7.18), the $M^n$ shrink to $c_*$.

(The argument for the satellite case is the same with $M_0$ replaced by $M$.) $\square$

*Remarks.* 1. Note that $\mathcal{M}^u$ is not compact in the satellite case. Indeed, by transverse quasi-conformality of $\mathcal{F}$ (Theorem 4.19), it would be qc equivalent to $M_0$. As $\sigma = \chi \circ R$, the satellite copy $M$ would be qc equivalent to $M_0$ as well, despite the fact that it does not have a cusp at the root point.

2. By the Self-similarity Theorem, the homeomorphism $\sigma : M^{n+1} \to M^n$ is almost linear near $c_*$, so that the local geometry of the little Mandelbrot sets $M^n$ is almost the same. However, it implies that the whole $M^{n+1}$ is almost isometric to $M^n$ only when the $M^n$ shrink to $c_*$ (i.e., when MLC holds at $c_*$).

## 8. The renormalization horseshoe with bounded combinatorics

In this section we will prove the Hyperbolicity Theorem and its consequences for bounded combinatorics.

8.1. *Construction of the horseshoe.* Let us pick a finite family $\mathcal{L} = \{M_k\}_{k=1}^d$ of disjoint Mandelbrot copies, and the corresponding renormalization operator

$$R \equiv R_\mathcal{L} : \mathrm{Dom}(R) \equiv \cup \hat{\mathcal{T}}_k \to \mathcal{QG},$$

where $\mathcal{T}_k \equiv \mathcal{T}_{\hat{M}_k}$. This set up will be carried through the whole section.

We say that a point $f \in \mathcal{QG}$ is *completely non-escaping* if there is a sequence $f_n \in \mathcal{QG}$, $n = 0, \pm 1, \pm 2, \ldots$, such that $Rf_n = f_{n+1}$ and

$$\operatorname{mod}(f_n) \geq \mu = \mu(f) > 0, \ n = 0, \pm 1, \pm 2, \ldots.$$

Consider the set $\mathcal{A} = \mathcal{A}_\mathcal{L}$ of all completely non-escaping points.

Recall that the *natural extension* of a map $\hat{R} : \hat{\mathcal{A}} \to \hat{\mathcal{A}}$ is defined as the lift of $R$ to the space of two-sided orbits $\mathbf{f} = \{f_n\}_{n \in \mathbb{Z}}$, $\hat{R}(\mathbf{f}) = \{Rf_n\}_{n \in \mathbb{Z}}$. Moreover, $\hat{R}$ is a homeomorphism with respect to the weak topology on $\mathcal{A}$. The projection $\phi : \mathbf{f} \mapsto f_0$ to the zero coordinate semi-conjugates $\hat{R}$ to $R$.

Let us also consider the space $\Sigma = \Sigma_\mathcal{L}$ of bi-infinite sequences $\bar{\tau} = \{M_{k(n)}\}_{n=-\infty}^{\infty}$ in symbols $M_k \in \mathcal{L}$ supplied with the weak topology. Let $\omega : \Sigma \to \Sigma$ stand for the left shift on this space (so that $M_{k(1)}$ is the zero coordinate of $\omega(\bar{\tau})$). It is called the *Bernoulli shift*.

Let $\Sigma^+ = \Sigma_\mathcal{L}^+$ stand for the space of one-sided sequences $\{M_{k(n)}\}_{n=0}^{\infty}$ in symbols $M_k \in \mathcal{L}$. Recall that by definition, the combinatorial type $\tau(f) \in \Sigma^+$ of an infinitely renormalizable map $f$ is the itinerary of the one-sided orbit $\operatorname{orb}(f)$ (see §5.5).

LEMMA 8.1. *Assume that there exists a $\nu = \nu_\mathcal{L} > 0$ such that for any one-sided sequence $\tau \in \Sigma^+$ there exists an infinitely renormalizable map $f$ with $\tau(f) = \tau$ and*

(8.1) $$\operatorname{mod}(R^n f) \geq \nu, \quad n = 0, 1, \ldots.$$

*Then the natural extension $\hat{R} : \hat{\mathcal{A}} \to \hat{\mathcal{A}}$ is topologically conjugate to the Bernoulli shift $\omega : \Sigma \to \Sigma$. Thus there exists a continuous map semi-conjugating $\omega$ to $R|\mathcal{A}$.*

*In particular, the statement is valid for a real family $\mathcal{L}$. Moreover, in this case the horseshoe $\mathcal{A}$ is real (i.e., consists of real maps), and the renormalization $R : \mathcal{A} \to \mathcal{A}$ is a homeomorphism topologically conjugated to the Bernoulli shift.*

*Proof.* Let us take a bi-infinite sequence $\bar{\tau} = \{M_k\} \in \Sigma$. By the assumption, for any $l \geq 0$, there is an infinitely renormalizable quadratic-like map $F_l$ with combinatorics $\tau(F_l) = \{M_{-l}, \ldots, M_0, \ldots\}$ and an *a priori* bound $\nu$. Let $f_{0,l} = R^l F_l$. These are infinitely renormalizable quadratic-like germs with common combinatorics $\tau_0 = \{M_0, M_1, \ldots\}$ and $\operatorname{mod}(f_{0,l}) \geq \nu$. Since the set of such maps is compact, we can pass to a quadratic-like limit $f_0 = \lim_{l \to \infty} f_{0,l}$ (along a subsequence) with the same properties.

Let us now do the same thing for every $i \leq 0$. Let $f_{i,l} = R^{l+i} F_l$, and let $f_i = \lim_{l \to \infty} f_{i,l}$ be a limit point. The map $f_i$ has combinatorics $\tau_i = \{M_i, M_{i+1}, \ldots\}$ and $\operatorname{mod}(f_i) \geq \nu$.

Selecting the above converging subsequences by means of the diagonal process, we construct a sequence of infinitely renormalizable quadratic-like maps $\{f_i\}_{i=-\infty}^{\infty}$ such that $Rf_i = f_{i+1}$, $\chi(f_i) = M_i$ and $\mathrm{mod}(f_i) \geq \nu$. This sequence represents a tower $\mathbf{f}$ with combinatorics $\bar{\tau}$ and a moduli bound $\nu$.

Thus any combinatorics $\bar{\tau} \in \Sigma$ is represented by a tower with *a priori* bounds. By the Tower Rigidity Theorem, the tower is unique. This provides us with a bijective map $\Phi : \Sigma \to \hat{\mathcal{A}}$ conjugating the Bernoulli shift $\omega$ to the natural extension $\hat{R}$. Taking the zero coordinate of the tower, we obtain a semi-conjugacy $\phi : \Sigma \to \mathcal{A}$.

To show that $\Phi$ is continuous, we need to check that the coordinate projections $\Sigma \to \mathcal{A}$ are continuous. To be definite, let us take the zero coordinate. If the continuity failed, there would exist two sequences of towers $\mathbf{f}^{(n)} = \{f_k^{(n)}\}_{k \in \mathbb{Z}}$ and $\mathbf{g}^{(n)} = \{g_k^{(n)}\}_{k \in \mathbb{Z}}$ in $\mathcal{A}$ such that $\chi(f_k^{(n)})$ and $\chi(g_k^{(n)})$ belong to the same $M_{i(k)}$ for $-n \leq k \leq n$, but $\mathrm{dist}_{\mathrm{M}}(f_0^{(n)}, g_0^{(n)}) \geq \varepsilon > 0$, where $\mathrm{dist}_{\mathrm{M}}$ is the Montel distance on $\mathcal{A}$. Passing to limits, we would obtain two different towers with the same combinatorics $\bar{\tau} \in \Sigma$ and *a priori* bounds – a contradiction.

In the case of a real family $\mathcal{L}$, the assumption of the lemma is satisfied by Theorem 5.6, and the above construction leads to a real set $\mathcal{A}$. Moreover, by Lemma 5.1, a real bi-infinite tower $\mathbf{f}$ is determined by its zero coordinate $f_0$. Hence $\Phi$ is a homeomorphism. It follows that $R : \mathcal{A} \to \mathcal{A}$ is a homeomorphism as well. □

The set $\mathcal{A} = \mathcal{A}_\mathcal{L}$ will be called the *renormalization horseshoe* (with combinatorics $\mathcal{L}$). *The assumptions of Lemma 8.1 will be the standing assumptions for the rest of this section.*

8.2. *Stable lamination.* Let $f_* \in \mathcal{A}$ and $f \in \mathcal{QG}$ be an infinitely renormalizable map (in the sense of $R$). The orbit of $f$ is *asymptotic* to the orbit of $f_*$ if there exist $\mu > 0$, a sequence of quadratic-like representatives $R^n f_* : V_n \to V_n'$ with $\mathrm{mod}(V_n' \setminus V_n) \geq \mu$, and an $N$ such that for $n \geq N$, $R^n f \in \mathcal{B}_{V_n}$ and

$$\|R^n f - R^n f_*\|_{V_n} \to 0 \text{ as } n \to \infty.$$

If under these circumstances

$$\|R^n f - R^n f_*\|_{V_n} \leq Cq^n,$$

with $C > 0$ and $q \in (0, 1)$, then we say that the orbit of $f$ is *exponentially asymptotic* to the orbit of $f_*$.

Let us say that the orbits of $f_* \in \mathcal{A}$ and $f \in \mathcal{H}(f_*)$ are *uniformly* exponentially asymptotic if the constants $\mu$, $C$, $q$ above are uniform, while $N$ depends only on $\mathrm{mod}(f)$.

The (*global*) *stable set* of a point $f_* \in \mathcal{A}$ is defined as the set of points $f \in \mathcal{QG}$ whose orbits are forward asymptotic to the orbit of $f_*$. (We avoid the usual term "global stable manifold" since in what follows the stable set will not be a manifold but rather a countable union of manifolds.) Let us define the *basin of attraction* of $\mathcal{A}$ as the union of the stable sets of all points $f_* \in \mathcal{A}$.

The following result extends the Stable Manifold Theorem to renormalization operators of bounded type and follows from the works of Sullivan, McMullen and the author in a similar way:

THEOREM 8.2 (Stable lamination). *The basin of attraction of the renormalization horseshoe coincides with the union of the hybrid classes $\mathcal{H}(f_*)$, $f_* \in \mathcal{A}$, and hence forms a lamination in $\mathcal{QG}$ with codimension-one complex analytic leaves. The orbits of $f \in \mathcal{H}(f_*)$ are uniformly exponentially asymptotic to the orbits of $f_* \in \mathcal{A}$.*

*Proof.* Let
$$\mathcal{H}(f_*, \nu) = \{f \in \mathcal{H}(f_*) : \mathrm{mod}(f) \geq \nu\}.$$
We shall show that the orbits of $f \in \mathcal{H}(f_*)$ are *uniformly asymptotic* to the orbits of $f_* \in \mathcal{A}$ in the following sense.

*Statement.* There exist a $\mu > 0$ and a choice of quadratic-like representatives $R^n f_* : V_n \to V_n'$, $f_* \in \mathcal{A}$ with $\mathrm{mod}(V_n' \setminus V_n) \geq \mu$ satisfying the following property. For any $\nu > 0$ and $\varepsilon > 0$, there exists an $N = N(\nu, \varepsilon)$ such that: If $f \in \mathcal{H}(f_*, \nu)$, then for $n \geq N$, $R^n f \in \mathcal{B}_{V_n}(R^n f_*, \varepsilon)$.

We leave the proof to the reader (it is a straightforward adjustment of the proof of the corresponding statement for the stationary combinatorics).

Let us now consider the analytic projection $\Pi : \mathcal{QG} \to \mathcal{H}_0$ (4.4) whose restrictions $\Pi_F : \mathcal{H}(F) \to \mathcal{H}_0$ are diffeomorphisms. Note that by the Product Structure Theorem 4.13, the inverse branches $\Pi_F^{-1} : \mathcal{H}_0 \to \mathcal{H}(F)$ are equicontinuous on compact sets. Let
$$R_F^n = \Pi_{R^n F} \circ R^n \circ \Pi_F^{-1} : \mathcal{H}_0 \to \mathcal{H}_0.$$
Then the above Statement (uniform contraction by $R$), compactness of $\mathcal{A}$, continuity of $\Pi$, and equicontinuity of its inverse branches imply that the family of operators $R_F$ is uniformly contracting as well:
$$R_F^N \mathcal{H}_{0, W(G)}(G, \varepsilon) \subset \mathcal{H}_{0, W(R^N G)}(R_F^N G, \delta),$$
where the $W(G)$ are appropriately selected domains of maps $G \in \Pi(\mathcal{A})$ with $\mathrm{mod}(W(G)) \geq \eta$. By the Schwarz Lemma, this family is uniformly infinitesimally contracting, and hence the iterates of $R$ are uniformly exponentially contracting. $\square$

Thus the hybrid class $\mathcal{H}(f_*)$ of a point $f_* \in \mathcal{A}$ is identified with the connected component of the stable set of $f$. It will also be called the (global) *stable leaf* of $f_*$. For an infinite horseshoe, the global stable set of $f_*$ is the union of infinitely many disjoint stable leaves. It has a dense intersection with $\mathcal{A}$. This is the usual picture for discontinuous hyperbolic maps like the baker transformation.

8.3. *Slow shadowing and hyperbolicity.* In this section we will state a new hyperbolicity criterion for invariant sets of *complex* analytic maps. It says that the lack of hyperbolicity can be detected topologically by the existence of slowly shadowing orbits. It will model a more complicated situation treated in the next section.

Let $L : \mathcal{A} \to \mathcal{A}$ be a continuous transformation of a metric compact set. One says that an $\mathrm{orb}_n(g)$ $\varepsilon$-*shadows* the $\mathrm{orb}_n(f)$ if $\mathrm{dist}(L^k f, L^k g) \leq \varepsilon$, $k = 0, \ldots, n$.

A map $L$ is said to satisfy a "specification property" (compare Bowen [Bo]) if for any $\varepsilon > 0$ there exists an $l = l(\varepsilon)$ such that any orbit $\{L^k f\}_{k=0}^{N-1}$, $f \in \mathcal{A}$, can be $\varepsilon$-shadowed by a periodic orbit of period at most $N + l$.

A basic example of a system satisfying the specification property is the Bernoulli shift $\omega : \Sigma_d \to \Sigma_d$. Note also that this property is preserved under taking quotients: If a map $L : \mathcal{A} \to \mathcal{A}$ satisfies the specification property and $S : \mathcal{X} \to \mathcal{X}$ is a quotient map (i.e., there exists a surjective continuous map $h : \mathcal{A} \to \mathcal{X}$ such that $S \circ h = h \circ L$), then $S$ satisfies the specification property as well.

Let $\mathcal{A}$ be embedded into a complex analytic Banach manifold $\mathcal{U}$, let $\mathcal{V}$ be a neighborhood of $\mathcal{A}$ in $\mathcal{U}$, and let $L : (\mathcal{V}, \mathcal{A}) \to (\mathcal{U}, \mathcal{A})$ be a complex analytic map preserving $\mathcal{A}$.

An orbit of $f \in \mathcal{V}$ *slowly $\varepsilon$-shadows* an orbit of $g \in \mathcal{A}$ if it $\varepsilon$-shadows the latter but is not exponentially asymptotic to it.

Let $\hat{L} : \hat{\mathcal{A}} \to \hat{\mathcal{A}}$ stand for the natural extension of $L$. Given an orbit $\mathbf{f} \in \hat{\mathcal{A}}$, we will denote by $f \equiv f_0 \in \mathcal{A}$ its zero coordinate. A map $L|\mathcal{A}$ is called *uniformly hyperbolic* if:

(i) There is an invariant subbundle $\mathcal{E}^s \subset \mathrm{T}_\mathcal{A}\mathcal{U}$, on which $\mathrm{D}L$ is uniformly exponentially contracting;

(ii) There exists a family of tangent subspaces $E_{\mathbf{f}}^u \subset \mathrm{T}_f\mathcal{U}$ labeled by points $\mathbf{f} \in \hat{\mathcal{A}}$ of the natural extension satisfying the following properties:
  – Transversality: $E_f^s \oplus E_{\mathbf{f}}^u = \mathrm{T}_f\mathcal{U}$;
  – Invariance: $E_{\hat{L}\mathbf{f}}^u = \mathrm{D}L(E_{\mathbf{f}}^u)$;
  – Uniformly exponential expansion: there exist $c > 0$ and $\rho > 1$ such that for any $\mathbf{f} \in \hat{\mathcal{A}}$, $v \in E_{\mathbf{f}}^u$:
  $$\|\mathrm{D}L_f^n(v)\| \geq c\rho^n, \ n = 0, 1, \ldots.$$

(In the invertible case, this definition coincides with the standard one.)

THEOREM 8.3 (Hyperbolicity criterion).  *Let $L : (\mathcal{V}, \mathcal{A}) \to (\mathcal{U}, \mathcal{A})$ be a complex analytic map satisfying the following properties*:

(i) *$L|\mathcal{A}$ is topologically transitive and satisfies the specification property.*

(ii) *There is an invariant complex codimension-one subbundle $\mathcal{E}^s \subset T_\mathcal{A}\mathcal{U}$, on which $DL$ is uniformly exponentially contracting.*

(iii) *$L$ is transversally nonsingular; i.e., the quotient maps*

$$D_* L_f : T_f \mathcal{U}/E^s_f \to T_{Rf}\mathcal{U}/E^s_{Rf}$$

*are invertible.*

(iv) *There is an $\varepsilon > 0$ such that $L|\mathcal{A}$ has no slowly $\varepsilon$-shadowing orbits.*

*Then $L$ is uniformly hyperbolic over $\mathcal{A}$.*

All the assumptions of Theorem 8.3 are satisfied for the renormalization operator of bounded type:

(i) is valid by Lemma 8.1 (renormalization horseshoe);
(ii) is true by the Theorem 8.2 (stable lamination);
(iii) is satisfied by Lemma 5.3;
(iv) is ensured by Theorem 5.9 (combinatorial rigidity).

The only property which fails is that $R$ acts on a Banach manifold. For this reason Theorem 8.3 cannot be directly applied to $R$. However, we can make some iterate of $R$ act fiberwise analytically on a Banach fiber space over $\mathcal{A}$.

The reader can figure a proof of Theorem 8.3 by a straightforward adjustment of the argument below to the manifold setting.

8.4. *Hyperbolicity of $R : \mathcal{A} \to \mathcal{A}$*. In this section we will give a proof of hyperbolicity of the renormalization operator $R$ of bounded type on the renormalization horseshoe $\mathcal{A}$. Let $\hat{R} : \hat{\mathcal{A}} \to \hat{\mathcal{A}}$ be the natural extension of the renormalization and $\phi : \hat{\mathcal{A}} \to \mathcal{A}$ be the natural projection. We will use bold letters $\mathbf{f}$, $\mathbf{g}$ etc. for points in $\hat{\mathcal{A}}$ and the corresponding usual italic letters $f$, $g$ etc. for their projections to $\mathcal{A}$.

THEOREM 8.4 (Hyperbolicity).  *For any $\mathbf{f} \in \hat{\mathcal{A}}$ there is a splitting $T_f \mathcal{QG} = E^s_f \oplus E^u_\mathbf{f}$ with the following properties*:

- $E^s_f = T_f \mathcal{H}(f)$ *and the action of $DR$ is uniformly exponentially contracting on the subbundle $\mathcal{E}^s = \cup_{f \in \mathcal{A}} E^s_f$;*

- $\dim E^u_\mathbf{f} = 1$. *This family of spaces is continuous and $DR$-invariant. The action of $DR$ on it is uniformly exponentially expanding (see §8.3).*

Let us start with the choice of Banach fibers. By compactness of $\mathcal{A}$, there is a choice of domains $W(f) \Subset W'(f)$ of quadratic-like germs $f \in \mathcal{A}$ satisfying the following properties:

W1. $\mathrm{mod}(W'(f) \setminus W(f)) \geq \mu$ with an absolute $\mu > 0$;

W2. There exists an $\eta > 0$ such that if $\mathrm{dist}_\mathrm{M}(g, f) < \eta$ for some $f, g \in \mathcal{A}$ then $g \in \mathcal{B}_f \equiv \mathcal{B}_{W(f)}$.

W3. There exist $\xi > 0$ and $N \in \mathbb{N}$ such that $R^N \mathcal{B}_f(f, \xi) \subset \mathcal{B}_{R^N f}$.

W4. The vertical fibers $\mathcal{Z}_f$ sit locally in $\mathcal{B}_f$, $f \in \mathcal{A}$; hence the vertical lines $E^v(f)$ (4.10) sit in the $\mathcal{B}_f$ as well.

The spaces $\mathcal{B}_f$ are the Banach fibers mentioned above. We will let $\mathcal{B}_f(\delta) \equiv \mathcal{B}_f(f, \delta)$ and
$$\mathrm{Per}_n(R) = \{g \in \mathcal{A} : R^n g = g\}.$$
Now, we consider the stable tangent bundle $\mathcal{E}^s$ over $\mathcal{A}$ with fibers $E^s_f = \mathrm{T}\mathcal{H}_{W(f)}(f)$, and the "normal bundle" $\mathcal{Y}$ over $\mathcal{A}$ with fibers $Y_f = \mathrm{T}_f \mathcal{B}_f / E^s_f$. By the Stable Manifold Theorem, $\dim Y_f = 1$.

Let $R_{*,f} : Y_f \to Y_{Rf}$ stand for the quotient action of the renormalization. Let
$$\gamma(R) = \inf_n \inf_{g \in \mathrm{Per}_n(R)} \|R^n_{*,g}\|^{1/n}_g,$$
where the norm $\|\cdot\|_g$ of $\mathrm{D}R^n_{*,g} : Y_g \to Y_{R^n g}$ is evaluated with respect to the quotient Banach norms on the fibers. Note that by Corollary 6.2 $\gamma(R) \geq 1$.

Let $R^N = L$, where $N$ satisfies Property W3. All further notation involving $L$ will be similar to the corresponding notation for $R$.

LEMMA 8.5. *For any $\lambda \in (0, \gamma(R))$ there is a constant $c = c_\lambda > 0$ such that*
$$\|R^n_{*,f}\| \geq c\lambda^n, \quad f \in \mathcal{A}, \ n \geq 0.$$

*Proof.* It is clearly enough to prove the desired property for $R$ replaced with $L$. By Property W2, $L$ can be locally trivialized: if $f, g \in \mathcal{A}$ and $\mathrm{dist}(f, g) < \delta$, then $L_g : (\mathcal{B}_f, g) \to (\mathcal{B}_{Lf}, Lg)$. Moreover, locally trivializing the normal bundle $\mathcal{Y}$, we make the quotient maps $L_{*,g}$ act on the same space $Y_f$. As $L_{*,g}$ continuously depends on $g$, for any $\varepsilon > 0$ there is a $\delta > 0$ such that $\|L_{*,f} \circ L^{-1}_{*,g} - I\|_f < \varepsilon$, provided $\mathrm{dist}(f, g) < \delta$.

It follows that
$$\|L^n_{*,f}(L^n_{*,g})^{-1} - I\|_f = O(n\varepsilon),$$
provided $\mathrm{orb}_n g$ $\delta$-shadows $\mathrm{orb}_n f$.

By the specification property, any $\mathrm{orb}_n(f)$ can be $\delta$-shadowed by an $\mathrm{orb}_n(g)$ of a periodic point $g$ of period at most $n+l$. Let $f_n = L^n f$, $g_n = L^n g$. It follows that

$$\begin{aligned} \|L^n_{*,f}\|_f &\geq \|L^n_{*,g}\|_f \exp(-Cn\varepsilon) \\ &\geq \|L^{n+l}_{*,g}\|_f \|L^l_{*,g_n}\|_{f_n}^{-1} \exp(-Cn\varepsilon) \geq B^{-l} \gamma(L)^{n+l} \exp(-Cn\varepsilon), \end{aligned}$$

where

$$B = \sup_{f \in \mathcal{A},\, \mathrm{dist}(g,f) < \delta} \|L_{*,g}\|_f$$

and $C$ are independent of $\varepsilon$. As $\gamma(L) \geq 1$, the conclusion follows. □

Given a tangent vector $h \in \mathcal{B}_f$, let $h^s$ and $h^v$ stand for its projections to $E_f^s$ and $E_f^v$ respectively. Let us consider a family of tangent cones

$$C_f^\theta = \{h \in \mathcal{B}_f : \|h^v\| \geq \theta \|h^s\|\}, \ f \in \mathcal{A}.$$

LEMMA 8.6. *For some $N$ and $\theta > 0$, the cone field $C_f^\theta$, $f \in \mathcal{A}$, is $R^N$-invariant. Moreover, there exists a continuous $R^N$-invariant family of complex lines $E_{\mathbf{f}}^u \subset \mathcal{B}_f$, $\mathbf{f} \in \hat{\mathcal{A}}$, complementary to $\mathcal{E}^s$.*

*Proof.* Let us show that for a sufficiently big $N$ and sufficiently small $\alpha > \theta > 0$,

(8.2) $$\mathrm{D}R^N C_f^\theta \subset C_{R^N f}^\alpha,$$

(in particular, the family of cones $C_f^\theta$, $f \in \mathcal{A}$, is $\mathrm{D}R^N$-invariant). Indeed, by the Stable Lamination Theorem and Lemma 8.5, there exist $\lambda \in (0,1)$, $\rho > \lambda$, such that for $S = R^l$

$$\|\mathrm{D}S_f h\| \leq \lambda \|h\|, \ h \in E_f^s; \quad \|\mathrm{D}S_f h\| \geq \rho \|h\|, \ h \in E_f^v.$$

Moreover, since the decomposition $\mathcal{B}_f = E_f^s \oplus E_f^v$ continuously depends on $f$,

$$\|(\mathrm{D}S_f h)^s\| \leq A\|h\|.$$

Take $\theta > 0$ so small that $\lambda + A\theta \equiv \mu < \rho$. Let $h \in \partial C_f^\theta$, so that $\|h^v\| = \theta \|h^s\|$. Then

$$\|(\mathrm{D}S_f h)^s\| \leq \lambda \|h^s\| + A\|h^v\| = \mu \|h^s\|.$$

Hence for $h \in \partial C_f^\theta$ we have:

$$\frac{\|(\mathrm{D}S_f h)^v\|}{\|(\mathrm{D}S_f h)^s\|} \geq \frac{\rho}{\mu} \frac{\|h^v\|}{\|h^s\|},$$

so that $\mathrm{D}S_f(\partial C_f^\theta) \subset C_{Sf}^\alpha$ with $\alpha = (\rho/\mu)\theta$. As the cones are convex, (8.2) follows.

Let us now consider the projective cone $P_f^\theta$, the space of lines in $C_f^\theta$. Let us supply it with the following hyperbolic distance. Take two lines $\Gamma_i \in P_f^\theta$. Select two points $\gamma_i \in \Gamma_i$, join them by the straight line and take the intersection $I = I^\theta(\gamma_1, \gamma_2)$ of this line with the cone $C_f^\theta$. We consider this interval $I$ as a model for the hyperbolic line and supply it with the corresponding hyperbolic metric. Then

$$\mathrm{dist}^\theta(\Gamma_1, \Gamma_2) = \mathrm{dist}_{\mathrm{hyp}}(\gamma_1, \gamma_2 | I).$$

This definition is independent of the choice of representatives $\gamma_1$ and $\gamma_2$ since the transition from one pair to another, $\gamma_1'$ and $\gamma_2'$, is carried by a Möbius transformation $I \to I'$ which preserves the hyperbolic metric.

Moreover, the hyperbolic distance $\mathrm{dist}^\alpha$ in the cone $C^\alpha$ strictly dominates $\mathrm{dist}^\theta$: there exists a $q > 1$ such that

$$\mathrm{dist}^\alpha(\Gamma_1, \Gamma_2) \geq q \mathrm{dist}^\theta(\Gamma_1, \Gamma_2),$$

since $I^\alpha(\gamma_1, \gamma_2)$ has a bounded hyperbolic length in $I^\theta(\gamma_1, \gamma_2)$. As $S : P_f^\theta \to P_{Sf}^\alpha$ is contracting from $\mathrm{dist}^\theta$ to $\mathrm{dist}^\alpha$, the map $S : P_f^\theta \to P_{Sf}^\theta$ is uniformly contracting in $\mathrm{dist}^\theta$. Hence for any two-sided orbit $\mathbf{f} = \{f_k\}_{k \in \mathbb{Z}}$ with itinerary $\mathbf{f} \in \hat{\mathcal{A}}$, the cones $J_{\mathbf{f}}^k \equiv \mathrm{D}S^k(C_{f_{-k}}^\theta)$ exponentially shrink to a single complex line $E_{\mathbf{f}}^u$ as $k \to +\infty$. Obviously, this family of lines is $\mathrm{D}S$-invariant.

If the orbits $\mathbf{f}$ and $\mathbf{g}$ with itineraries $\bar{\tau}$ and $\bar{\eta}$ respectively are so close that their backward pieces of length $k$ $\delta$-shadow one another, then the cones $J_{\mathbf{f}}^k$ and $J_{\mathbf{g}}^k$ are also close and localize well the lines $E_{\mathbf{f}}^u$ and $E_{\mathbf{g}}^u$. This shows continuity of the line field. □

Let $D_f^s(\delta) \subset E_f^s(\delta)$ and $D_{\mathbf{f}}^u(\delta) \subset E_{\mathbf{f}}^u(\delta)$ stand for the $\delta$-disks about $f = \phi(\mathbf{f})$ in the stable/unstable subspaces $E_f^s$ and $E_{\mathbf{f}}^u$ respectively. Let

(8.3) $$D_{\mathbf{f}}(\delta) = D_f^s(\delta) \times D_{\mathbf{f}}^u(\delta)$$

stand for the corresponding $\delta$-bidisks, and let $\partial^u D_{\mathbf{f}}(\delta) = D_f^s(\delta) \times \partial D_{\mathbf{f}}^u(\delta)$ be their horizontal boundaries.

LEMMA 8.7. *If $\gamma(R) = 1$ then $R$ has a slowly shadowing orbit.*

*Proof.* Let us consider the disjoint union $\bar{\mathcal{B}} = \sqcup_{\mathbf{f} \in \hat{\mathcal{A}}} \mathcal{B}_f(\delta_0)$ supplied with the topology induced from $\hat{\mathcal{A}} \times \mathcal{B}$. Let

$$\mathcal{D}(\delta) = \sqcup_{\mathbf{f} \in \hat{\mathcal{A}}} D_{\mathbf{f}}(\delta) \subset \bar{\mathcal{B}}.$$

The renormalization $R^N$ gives rise to an operator $\bar{L} : \mathcal{D}(\delta) \to \bar{\mathcal{B}}$ acting fiberwise. By the Stable Manifold Theorem and Lemma 8.6, this operator is exponentially horizontally contracting and has an invariant cone field in the bidisk family $\mathcal{D}(\delta)$ (provided $\delta$ is sufficiently small).

For $\lambda \in (0, 1)$, let us consider a fiberwise linear contraction $T_\lambda : \bar{\mathcal{B}} \to \bar{\mathcal{B}}$. Perturb $\bar{L}$ by post-composing it with this contraction: $\bar{L}_\lambda = T_\lambda \circ \bar{L}$. If an $\bar{L}$-periodic point $\mathbf{f} \in \hat{\mathcal{A}}$ becomes attracting under this perturbation, then by the Small Orbits Theorem, there is a $t$ and a point $\mathbf{g} \in \partial^u D_{\hat{L}^t \mathbf{f}}(\delta)$ in the basin of $\mathbf{f}$ such that $\bar{L}^k \mathbf{g} \in D_{\hat{L}^{k+t} \mathbf{f}}(\delta)$, $k = 0, 1, \ldots$.

Since $\gamma(R) = 1$, for any $\lambda \in (0,1)$ there is an attracting periodic point $\mathbf{f}_\lambda$ and the corresponding shadowing point $\mathbf{g}_\lambda \in \partial^u D_{\mathbf{f}_\lambda}(\delta)$. But by Lemma 4.1, the set $\mathcal{D}(\delta)$ is compact in $\bar{\mathcal{B}}$. Passing to limits $\mathbf{f} = \lim \mathbf{f}_{\lambda_k} \in \hat{\mathcal{A}}$ and $\mathbf{g} = \lim \mathbf{g}_{\lambda_k} \in D_{\mathbf{f}}(\delta)$ we conclude that $f$ is shadowed by $g \notin \mathcal{W}^s_{\text{loc}}(f)$. □

*Proof of the hyperbolicity Theorem* 8.4. By the Combinatorial Rigidity Theorem, $R|\mathcal{A}$ does not have slowly shadowing orbits. Hence by Lemma 8.7, $\gamma(R) > 1$, so that the periodic points of $R$ are uniformly repelling in the transverse direction. By Lemma 8.5, $R$ is exponentially expanding in the transverse direction. Hence by Lemma 8.6, $R^N|\mathcal{A}$ is exponentially expanding on the family of unstable lines. Together with the Stable Lamination Theorem this yields uniform hyperbolicity of $R^N|\mathcal{A}$.

By Theorem 6.3 (hyperbolicity in the stationary case), the unstable lines $E^u_{\mathbf{f}}$ are uniquely determined for periodic points $\mathbf{f} \in \hat{\mathcal{A}}$ (i.e., independent of the choice of Banach spaces $\mathcal{B}_f$ and the iterate $R^N$), and form a $DR$-invariant family. Since the full family of unstable lines $E^u_{\mathbf{f}}$, $\mathbf{f} \in \hat{\mathcal{A}}$, is continuous, it is uniquely determined and $DR$-invariant everywhere. □

8.5. *Unstable manifolds of $R|\mathcal{A}$*. We will keep the notations of the previous section. In particular, $\mathcal{B}_f$ will stand for the family of Banach slices satisfying properties W1–W4.

THEOREM 8.8 (Local unstable manifolds). *There exists a continuous family of complex one-dimensional analytic manifolds $\mathcal{W}^u_{\text{loc}}(\mathbf{f}) \subset \mathcal{B}_f$ through $f = \phi(\mathbf{f}) \in \mathcal{A}$ satisfying the following properties*:

(i) $\mathcal{W}^u_{\text{loc}}(\mathbf{f})$ *is tangent to $E^u_{\mathbf{f}}$ and transverse to $\mathcal{W}^s_{\text{loc}}(f)$*.

(ii) $\mathcal{W}^u_{\text{loc}}(\hat{R}\mathbf{f}) \subset R\mathcal{W}^u_{\text{loc}}(\mathbf{f})$; *thus the inverse map $R^{-1}_{\mathbf{f}} : \mathcal{W}^u_{\text{loc}}(\hat{R}\mathbf{f}) \to \mathcal{W}^u_{\text{loc}}(\mathbf{f})$ is well-defined*.

(iii) *There exist $\rho \in (0,1)$ and $C > 0$ such that for any $\mathbf{g} \in \mathcal{W}^u_{\text{loc}}(\mathbf{f})$*,
$$\|g_{-n} - f_{-n}\|_{f_{-n}} \leq C\rho^n,$$
*where $f_{-n} = \phi(\hat{R}^{-n}\mathbf{f})$ and $g_{-n}$ are the corresponding preimages of $\mathbf{g}$*.

(iv) *The straightenings $\chi : \mathcal{W}^u_{\text{loc}}(\mathbf{f}) \to M_0$ are uniformly $K$-qc*.

*Proof*. As in the proof of Lemma 8.7, let us lift the iterate $L = R^N$ to the hyperbolic fibered Banach operator $\bar{L} : \mathcal{D}(\delta) \to \bar{\mathcal{B}}$, and let $\hat{L} = \hat{R}^N$. By the Hyperbolicity Theorem, $\bar{L}$ is hyperbolic. By a standard construction, it generates a continuous family of local unstable manifolds $\mathcal{W}^u_{\text{loc}}(\mathbf{f})$ satisfying the above properties.

The construction is the following. Let us consider a family $\mathcal{G}_{\mathbf{f}}$ of complex analytic curves $\Gamma \subset \mathcal{D}_{\mathbf{f}}$ which are represented by the graphs $\psi : D^u_{\mathbf{f}} \to D^s_{\mathbf{f}}$ of analytic functions whose tangent lines belong to the cones $C^{\pi/4}_{\mathbf{f}}$. Supply this family with the uniform norm. Let us define an operator $\bar{L}_*$ on $\mathcal{G}_{\mathbf{f}}$ as the truncated image of the functions: $\psi \mapsto \bar{L}\psi \cap D_{\hat{L}\mathbf{f}}$. Since $\bar{L}$ is horizontally contracting, vertically expanding, and preserves the cone field (for $N$ big enough), it maps $\mathcal{G}_{\mathbf{f}}$ into $\mathcal{G}_{\hat{L}\mathbf{f}}$ uniformly contracting the distance. It follows that there exists a

$$\lim_{n\to\infty} \bar{L}^n_* \psi_{-n} \equiv \mathcal{W}^u_\delta(\mathbf{f}) \in \mathcal{G}_{\mathbf{f}},$$

where $\psi_{-n}$ is an arbitrary function of the family $\mathcal{G}_{\hat{L}^{-n}\mathbf{f}}$. All the desired properties of this family (with $R$ replaced by $L$) are obvious.

Passing from $L = R^N$ back to $R$ we need to ensure property (ii). It is easy to see from the construction that the manifolds $R\mathcal{W}^u_\delta(\mathbf{f})$ and $\mathcal{W}^u_\delta(\hat{R}\mathbf{f})$ represent the same germ at $Rf = \phi(\hat{R}\mathbf{f})$. Moreover, clearly by taking $\delta$ sufficiently small we can make the iterates $R^k$, $k = 0, 1 \ldots, N$, to be well-defined on the $\mathcal{W}^u_\delta(\mathbf{f})$. Let us now define the unstable manifolds as

$$\mathcal{W}^u_{\mathrm{loc}}(\mathbf{f}) = \bigcup_{0 \leq k \leq N} R^k \mathcal{W}^u_\delta(\hat{R}^{-k}\mathbf{f}).$$

Clearly this new family satisfies the properties (i)–(iii).

Property (iv) follows from Theorem 4.19.  □

Similarly, as in the stationary case, we can now globalize the unstable lamination. For $\mathbf{f} \in \mathcal{A}$, let us define the *unstable Mandelbrot set* $\mathcal{M}^u(\mathbf{f})$ as the set of infinitely anti-renormalizable points $g \in \mathcal{C}$ such that there exists a one-sided tower $\mathbf{g} = \{g_{-n}\}_{n \in \mathbb{N}}$ with the following property: $g_{-n} \in \mathcal{B}_{f_{-n}}$ for sufficiently big $n$ and $\|g_{-n} - f_{-n}\|_{f_{-n}} \to 0$ as $n \to \infty$.

THEOREM 8.9 (Global unstable leaves). *Let $\mathbf{f} \in \mathcal{A}$, $\tau_- = \tau_-(\mathbf{f}) = \{M_{-n}\}_{n \in \mathbb{N}}$. Then*

(i) *A point $g \in \mathcal{C}$ belongs to $\mathcal{M}^u(\mathbf{f})$ if and only if there exists a one-sided tower $\mathbf{g} = \{g_0, g_{-1}, \ldots\}$ with combinatorics $\tau_-$ and a priori bounds. Moreover, in this case $g_{-n} \in \mathcal{W}^u_{\mathrm{loc}}(\mathbf{f}_{-n})$ for all sufficiently big $n$, where $\mathbf{f}_{-n} = \hat{R}^{-n}\mathbf{f}$.*

(ii) *The straightening $\mathcal{M}^u(\mathbf{f}) \to M_0$ is injective.*

(iii) *For any $\mu > 0$, the set $\mathcal{M}^u_\mu(\mathbf{f}) = \{g \in \mathcal{M}^u : \exists$ a one-sided tower $\mathbf{g}$ with $g = g_0$ and $\mathrm{mod}(\mathbf{g}) \geq \mu\}$ is embedded into a one-dimensional complex analytic manifold $\mathcal{W}^u_\mu(\mathbf{f})$ which extends the local manifold $\mathcal{W}^u_{\mathrm{loc}}(\mathbf{f})$.*

(iv) *The manifold $\mathcal{W}^u_\mu(\mathbf{f})$ is transverse to the foliation $\mathcal{F}$.*

(v) *The germ of the manifold $\mathcal{W}^u_\mu(\mathbf{f})$ near $\mathcal{C}$ is mapped to the germ of $\mathcal{W}^u_\mu(\hat{R}\mathbf{f})$ under the renormalization.*

(vi) *The straightening $\mathcal{M}^u_\mu(\mathbf{f}) \to M_0$ is $K(\mu)$-quasi-conformal.*

*Proof.* (i) If $\mathbf{g} = \{g_0, g_{-1}, \ldots\}$ is a tower with combinatorics $\tau_-(\mathbf{f})$ and *a priori* bounds, then $g_{-n} \in \mathcal{B}_{f_{-n}}$ for sufficiently big $n$ and $\|g_{-n} - f_{-n}\|_{f_{-n}} \to 0$ as $n \to \infty$. Otherwise we would construct by the diagonal process a bi-infinite tower $\{h_n\}_{n=-\infty}^{\infty}$ with some combinatorics $\bar{\tau} \in \Sigma_\mathcal{L}$ but different from the tower in $\mathcal{A}$ with the same combinatorics. This would contradict the Tower Rigidity Theorem.

The rest of the argument follows the lines of the proof of Theorem 6.5 replacing the fixed point with the inverse orbit. □

The manifolds $\mathcal{W}^u(\mathbf{f}) = \cup_{\mu>0} \mathcal{W}^u_\mu(\mathbf{f})$ constructed above will be called (global) unstable *leaves*. The global unstable set of $\mathbf{f}$ (i.e., the set of points $\mathbf{g}$ whose backward orbits are asymptotic to the backward orbit of $\mathbf{f}$) consists of infinitely many leaves.

8.6. *Real horseshoe.* In the same way as for stationary combinatorics (Theorem 6.6), the above results can be refined in the case of real combinatorics. In the following statement all the sets ($\mathcal{QG}$, $\mathcal{H}(f)$, stable and unstable leaves) mean the real slices of the sets considered above (without change of notation).

THEOREM 8.10 (Real horseshoe). *Let $\mathcal{L} \subset \mathcal{N}$ be a finite family of real Mandelbrot sets and $R = R_\mathcal{L}$ be the corresponding renormalization operator. Then*:

(i) *There exists a real compact invariant set $\mathcal{A} \subset \mathcal{C}$ on which $R$ is topologically conjugate to the shift $\omega : \Sigma_\mathcal{L} \to \Sigma_\mathcal{L}$.*

(ii) *The renormalization operator $R$ is uniformly hyperbolic on $\mathcal{A}$.*

(iii) *The stable leaf $\mathcal{W}^s(f)$, $f \in \mathcal{A}$, coincides with the hybrid class $\mathcal{H}(f)$; $\operatorname{codim} \mathcal{W}^s(f) = 1$.*

(iv) *For any $\delta > 0$ there exists a $\mu > 0$ such that every unstable leaf $\mathcal{W}^u_\mu(f)$, $f \in \mathcal{A}$, transversally passes through all real hybrid classes $\mathcal{H}_c$ with $c \in [-2, 1/4-\delta]$; $\dim \mathcal{W}^u_\mu(f) = 1$. These unstable leaves are pairwise disjoint.*

(v) *The straightening $\mathcal{W}^u_\mu(f) \to [-2, 1/4 - \delta]$ is $K(\delta)$-quasi-symmetric.*

*Remark.* The disjointness of the unstable leaves follows from the injectivity of $R$ (Lemma 5.1).

## 9. Applications of the renormalization horseshoe

In this section we will derive from the above renormalization theory the Hairiness, Self-similarity and Universality Theorems in the case of bounded combinatorics. The exposition will be sketchy as it follows the lines of the stationary case (where the renormalization fixed point is replaced by the orbits of the renormalization horseshoe). At the end we will prove the HD Theorem.

We keep our standing assumptions (of Lemma 8.1) which ensure existence of the hyperbolic horseshoe $\mathcal{A}$.

9.1. *Distortion and linearization.* Given a conformal diffeomorphism $\phi : X \to Y$ between two complex one-dimensional Riemannian manifolds, the *distortion* (or *nonlinearity*) of $\phi$ is defined as

$$n(\phi) = \sup_{z,\zeta \in X} \log \frac{\|D\phi(z)\|}{\|D\phi(\zeta)\|}.$$

By the Koebe Distortion Theorem, if $X$ and $Y$ are hyperbolic planes then $n(\phi)$ is uniformly bounded in any hyperbolic disk of radius $r > 0$ (with a bound depending on $r$ only). Moreover, $n(r) = O(r)$ with an absolute constant, as $r \to 0$.

Let us supply each unstable leaf $\mathcal{W}^u_{\mathrm{loc}}(\mathbf{f})$, $\mathbf{f} \in \hat{\mathcal{A}}$, with the Riemannian metric induced from the Banach space $\mathcal{B}_f \equiv \mathcal{B}_{W(f)}$, $f = \phi(\mathbf{f}) \in \mathcal{A}$ (see §8.4). Let us make a few remarks on this family of metrics:

• Since the $\cup_{\mathbf{f}\in\hat{\mathcal{A}}}\mathcal{W}^u_{\mathrm{loc}}(\mathbf{f})$ sit in a compact part of $\mathcal{QG}$, these metrics are uniformly equivalent (after perhaps a slight shrinking of the domains $W(f)$) to the metric induced from a single Banach space $\mathcal{B}_V$.

• Supply the unit disk $\mathbb{D}$ with the Euclidean metric. Then the uniformizations $\phi_{\mathbf{f}} : \mathcal{W}^u_{\mathrm{loc}}(\mathbf{f}) \to \mathbb{D}$ have a uniformly bounded distortion. Indeed, in the local coordinate systems $\mathcal{B}_{\mathbf{f}} = E^s_{\mathbf{f}} \oplus E^u_{\mathbf{f}}$ the manifolds $\mathcal{W}^u_{\mathrm{loc}}(\mathbf{f})$ have a uniformly bounded vertical slope. Hence the projections

(9.1) $$p_{\mathbf{f}} : \mathcal{W}^u_{\mathrm{loc}}(\mathbf{f}) \to D^u_{\mathbf{f}}(\delta)$$

have a uniformly bounded distortion. But since the space $E^u_{\mathbf{f}}$ is one-dimensional, the Banach disk $D^u_{\mathbf{f}}(\delta)$ is linearly conformally equivalent to the standard Euclidean disk $\mathbb{D}$.

Moreover, as the slopes of the $\mathcal{W}^u_{\mathrm{loc}}(\mathbf{f})$ are uniformly bounded, by the Cauchy inequality, their graphs have uniformly bounded second derivatives as well. Hence the distortion of $p_{\mathbf{f}}$ uniformly linearly vanishes as $\delta \to 0$, i.e. $n(p_{\mathbf{f}}) = O(\delta)$.

It follows that the Koebe Distortion Theorem is valid with respect to the Banach metrics on the leaves $\mathcal{W}^u_{\mathrm{loc}}(\mathbf{f})$. This yields for $R$ the usual distortion estimates. Namely, let $R_{\mathbf{f}}^{-n}$ mean the inverse branches of $R^{-n}|\mathcal{W}_{\text{loc}}^u(\mathbf{f})$ corresponding to the backward orbit $\mathbf{f}$ and $\mathrm{D}^u$ means the differential in the direction of $\mathcal{W}_{\text{loc}}^u(\mathbf{f})$.

LEMMA 9.1 (Distortion).  *There exists a $C > 1$ such that for any $\mathbf{f} \in \hat{\mathcal{A}}$ and $g \in \mathcal{W}_{\text{loc}}^u(\mathbf{f})$,*
$$C^{-1} \leq \frac{\|\mathrm{D}^u R_{\mathbf{f}}^{-n}(f)\|}{\|\mathrm{D}^u R_{\mathbf{f}}^{-n}(g)\|} \leq C.$$
*Moreover, if* $\operatorname{dist}(f, g) \leq \varepsilon$ *then*
$$\frac{\|\mathrm{D}^u R_{\mathbf{f}}^{-n}(f)\|}{\|\mathrm{D}^u R_{\mathbf{f}}^{-n}(g)\|} = 1 + O(\varepsilon).$$

One can go further and linearize $R$ along the unstable lamination:

LEMMA 9.2 (Linearization).  *There is a family of conformal local charts $\psi_{\mathbf{f}} : \mathcal{W}_{\text{loc}}^u(\mathbf{f}) \to \mathbb{C}$, $\mathbf{f} \in \hat{\mathcal{A}}$, with uniformly bounded distortion, and a function $\lambda : \hat{\mathcal{A}} \to \mathbb{C}$ which linearize $R$:*
$$\psi_{\hat{R}\mathbf{f}}(Rg) = \lambda(\mathbf{f})\psi_{\mathbf{f}}(g).$$

*Proof.* Let us start with local charts provided by the family of projections (9.1). Let $\lambda(\mathbf{f})$ denote the derivative of $R|\mathcal{W}_{\text{loc}}^u(\mathbf{f})$ at $f = \phi(\mathbf{f})$ with respect to these local charts at $\mathbf{f}$ and $\hat{R}\mathbf{f}$, and let
$$q_n(\mathbf{f}) = \prod_{m=0}^{n-1} \lambda(\hat{R}^k \mathbf{f})$$
stand for the corresponding $n$-fold derivatives (which form a $\mathbb{C}$-valued multiplicative cocycle).

Consider the backward orbit $\{\mathbf{f}_{-n}\} = \hat{R}^{-n}\mathbf{f}$, and the corresponding backward orbits $\{f_{-n} = \phi(\mathbf{f}_{-n})\}$ and $\{g_{-n}\}$ of some $g \in \mathcal{W}_{\text{loc}}^u(\mathbf{f})$. Let
$$\phi_{\mathbf{f}}(g) = \lim_{n \to \infty} q_n(\mathbf{f}_{-n})\left(p_n(g_{-n}) - p_n(f_{-n})\right),$$
where $p_n \equiv p_{\mathbf{f}_{-n}}$. By the Distortion Lemma, the ratio of two consecutive terms of the above sequence goes to 0 at a uniformly exponential rate, and hence the above limit exists and represents a conformal chart on $\mathcal{W}_{\text{loc}}^u(\mathbf{f})$. It is obvious that these charts provide a desired linearization of $R$. $\square$

9.2. *Hairiness.* The Hairiness Conjecture at a Feigenbaum parameter value $c \in M_0$ of bounded type is stated in the same way as for the stationary combinatorics. Namely, let $r = r(\varepsilon)$ be the maximal number such that $\mathbb{D}(b, r\varepsilon) \subset \mathbb{D}(c, \varepsilon) \setminus M_0$ (i.e., the maximal relative size of gaps in $M_0$ in scale $\varepsilon$). Then

(9.2) $\qquad\qquad\qquad r(\varepsilon) \to 0 \quad \text{as} \quad \varepsilon \to 0.$

*Proof of the Hairiness Conjecture* (*bounded combinatorics*). Assume that (9.2) fails, so that there is a decaying sequence of scales around $c$ in which $M_0$ has definite gaps. Since the holonomy $\mathcal{W}^u_{\text{loc}}(\mathbf{f}) \to M_0$, $\mathbf{f} \in \mathcal{A}$, is uniformly qc (Theorem 8.8), on any unstable manifold $\mathcal{W}^u_{\text{loc}}(\mathbf{f})$ there is a sequence of scales in which the corresponding Mandelbrot set $\mathcal{M}^u(\mathbf{f})$ has definite gaps (where the relative size $r > 0$ of the gaps is uniform over the family of unstable manifolds).

Let us take a point $\mathbf{f} \in \hat{\mathcal{A}}$ and a corresponding definite gap $U \subset \mathcal{W}^u_{\text{loc}}(\mathbf{f})$ in some scale $\varepsilon > 0$. We push it forward by the renormalization, $U_k = R^k U$, $k = 0, 1, \ldots, n$, where $n$ is selected in such a way that the gap $U_n$ has size of order 1. By the Distortion Lemma, at this moment $U_n$ is an "ellipse" of bounded shape whose size is commensurable with the size of the unstable manifold $\mathcal{W}^u_{\text{loc}}(\mathbf{f}_n)$, where $\mathbf{f}_n = \hat{R}^n \mathbf{f}$.

Take now two points $g, \tilde{g} \in U$ whose distance apart is of order $\varepsilon$ and whose distances to the boundary $\partial U$ have the same order. Let $g_k = R^k g$, $\tilde{g}_k = R^k \tilde{g}$, $k = 0, 1, \ldots, n$. By bounded distortion, the distance between $g_k$ and $\tilde{g}_k$, their distances to the boundary $\partial U_k$, and their distances to $f_k = \phi(\mathbf{f}_k)$, are all comparable. In the proof of the Hairiness Conjecture for stationary combinatorics we have shown that this property implies that $g_k$ and $\tilde{g}_k$ are $K$-qc conjugate by a map staying a definite uniform distance away from the identity.

Take now a middle iterate $h_n = g_l$, $l = [n/2]$, as the zero coordinate of the tower $\mathbf{h}_n = \{R^k h_n\}_{k=-l}^{l}$, and pass to a limit as $n \to \infty$. We obtain a bi-infinite tower with bounded combinatorics and *a priori* bounds which admits a nontrivial automorphism contradicting Corollary 5.12. □

9.3. *Self-similarity and universality.*

LEMMA 9.3 (Transverse conformality). *There exists an $\alpha = \alpha(\mathcal{L}) > 0$ such that the foliation $\mathcal{F}$ is transversally $C^{1+\alpha}$-conformal along any stable leaf $\mathcal{H}(g)$, $g \in \mathcal{A}$.*

*Proof.* Let us take a transversal $S$ to $\mathcal{H}(g)$ at some point $f$. It is enough to study the holonomy $h$ from $S$ to an unstable manifold $\mathcal{W}^u_{\text{loc}}(\mathbf{g})$. We consider a family of bidisks $D_{\mathbf{f}} = D^s_f \times D^u_{\mathbf{f}}$ (8.3) and an iterate $L = R^N$ which acts hyperbolically on this family (uniformly contracting in the horizontal direction and uniformly expanding in the vertical).

Let us consider the forward orbit $\mathbf{g}_m = \hat{L}^m \mathbf{g} \in \hat{\mathcal{A}}$ and the corresponding sequence of bidisks $D_{\mathbf{g}_m} \equiv D_m = D^s_m \times D^u_m$. Let $\mathcal{S}^m$ denote the connected component of $D_m \cap L^m \mathcal{S}$ containing $f_m = L^m f$. Then by hyperbolicity of $L$, the $\mathcal{S}^m$ can be eventually represented as graphs of analytic functions $D^u_m \to D^s_m$ with bounded vertical slope. Moreover these graphs are exponentially close to the corresponding unstable manifolds $\mathcal{W}^u_m \equiv \mathcal{W}^u_{\text{loc}}(\mathbf{g}_m)$.

Let $q > \mu > 1$, $\mathcal{M} = \mathcal{S} \cap \mathcal{C}$. Take two points $z_1, z_2 \in \mathcal{M}$ whose distance from $a \equiv f$ is of order $\varepsilon = q^{-n}$, and iterate them forward by $L^m$ until they go to points $\zeta_1, \zeta_2$, whose distance from $b = L^m a$ is of order $\mu^{-n}$. By the Distortion Lemma this transition has distortion of order $\mu^{-n}$, as in (7.13).

Furthermore, the manifolds $\mathcal{S}^m$ and $\mathcal{W}_m^u$ are exponentially close and hence by the $\lambda$-lemma the holonomy between them has exponentially small ratio distortion at points $\zeta_1, \zeta_2, b$, as in (7.15).

Applying $L^{-m}$, we conclude that the holonomy $h$ has ratio distortion of order $\varepsilon^\alpha$ in scale $\varepsilon$ about $a$, which yields $C^{1+\alpha}$-conformality at $a$.    □

THEOREM 9.4 (Self-similarity for bounded combinatorics). *Let $c$ be a Feigenbaum parameter value of bounded type satisfying our standing hypotheses, and $M \in \mathcal{L}$ be the Mandelbrot copy containing $c$. Then the homeomorphism $\sigma : M \to M_0$ is $C^{1+\alpha}$ at $c$.*

*Proof.* This follows from the transverse conformality of the foliation $\mathcal{F}$ (Lemma 9.3) and analyticity of the renormalization $R$ on the unstable lamination.    □

9.4. *Universality for bounded combinatorics.* Let us restrict ourselves to the real case, as the complex statement is obtained by the usual adjustment. We take a finite family $\mathcal{L}$ of real Mandelbrot sets $M_k$ centered at points $a_k \in (-2, 1/4)$. Consider a Feigenbaum parameter value $c_* \in (-2, 1/4)$ with combinatorics $\{M_{i(1)}, M_{i(2)}, \ldots\}$. Also we consider the centers $c_n = a_{i(n)} * \ldots * a_{i(1)}$ of the $n$-fold-tuned Mandelbrot copies.

THEOREM 9.5 (Universality for bounded combinatorics). *The sequence $c_n$ exponentially converges to $c_*$:*

$$(9.3) \qquad b\lambda^{-n} \leq |c_n - c_*| \leq B\Lambda^{-n},$$

*with some $B, b > 0$ and $\Lambda, \lambda > 1$.*

*Let $\mathcal{S} = \{f_\mu\}$ be a real analytic one-parameter family of quadratic-like maps transversally intersecting the hybrid class $\mathcal{H}_{c_*}$ at $\mu_*$. Then for $\mu$ near $\mu_*$ and all sufficiently big $n$, $\mathcal{S}$ has a unique intersection point with the hybrid class $\mathcal{H}_{c_n}$, and*

$$(9.4) \qquad |\mu_n - \mu_*| = a\,|c_n - c_*|(1 + O(q^n)),$$

*where $a = a(\mathcal{S}) > 0$ and $q = q(\mathcal{L}) \in (0, 1)$.*

*Proof.* By Lemma 8.1, there is a map $f_* \in \mathcal{A}$ such that $\chi(f_*) = c_*$. By Theorem 8.10, there exist maps $f_n \in \mathcal{W}^u(f_*)$ such that $\chi(f_n) = c_n$. Then

$$R^m f_n \in \mathcal{W}^u(R^m f_*) \quad \text{and} \quad \chi(R^m f_n) = c_{n-m}.$$

Since $\text{dist}(R^n f_n, R^n f_*)$ is bounded and $R^{-n}$ exponentially contracts the unstable leaves, $f_n \to f_*$ exponentially fast (with the rate depending only on the family $\mathcal{L}$). Since by Lemma 9.3 the foliation $\mathcal{F}$ is transversally $C^{1+\alpha}$, (9.3) follows.

Applying $C^{1+\alpha}$-conformality once more, we obtain (9.4). □

Thus for any bounded combinatorics $\tau = (M_0, M_1, \ldots)$ there is a universal scaling law of transition in generic one-parameter families to the parameters with these combinatorics.

9.5. *Hausdorff dimension.* Let $L_d = \{1, \ldots, d\}$, $d > 1$. Consider a hierarchical family of interval $I_{\bar{i}}^n \subset \mathbb{R}$, where $\bar{i} = (i_0, \ldots, i_{n-1})$, $i_k \in L_d$, $n = 0, 1, \ldots$. Assume that the intervals of a given rank $n$ are pairwise disjoint, while for any $j \in L_d$, $I_{\bar{i},j}^{n+1} \subset I_{\bar{i}}^n$. The components $G_{\bar{i},j}^{n+1}$ of $I_{\bar{i}}^n \setminus \cup_j I_{\bar{i},j}^{n+1}$ are called the *gaps* of rank $n+1$. Let

$$Q = \bigcap_n \bigcup_{\bar{i}} I_{\bar{i}}^n.$$

The set $Q$ is called a *Cantor set with bounded geometry* if the family of configurations $(I_{\bar{i}}^n, \cup I_{\bar{i},j}^{n+1})$ have uniformly bounded geometry; i.e., the intervals $I_{\bar{i},j}^{n+1}$ and the gaps $G_{\bar{i},j}^{n+1}$ are commensurable with $I_{\bar{i}}^n$ (with a constant independent of $n$ and $\bar{i}$). It is a well-known and simple fact that a Cantor set with bounded geometry has Hausdorff dimension strictly in between 0 and 1.

Now, consider a finite family $\mathcal{L} = \{M_k\}$ of real Mandelbrot copies centered at $c_k$. Recall that $I_{\mathcal{L}} \subset (-2, 1/4)$ stand for the set of infinitely renormalizable real parameter values of type specified by this family (see the introduction).

LEMMA 9.6. *The set $I_{\mathcal{L}}$ is a Cantor set with bounded geometry (depending on $\mathcal{L}$).*

*Proof.* Consider the renormalization windows $I_k$ obtained by removing from $M_k \cap \mathbb{R}$ a neighborhood of the cusp $b_k$. For instance, let us remove the intervals $(c_k, b_k]$. On the union of these windows we have the straightened renormalization operator:

$$\sigma = \chi \circ R : \bigcup_{1 \leq k \leq d} I_k \to [-2, 0].$$

For $\bar{i} = (i_0, \ldots, i_{n-1})$, let $I_{\bar{i}}^n = \{c : \sigma^k(c) \in I_{i_k}, \ k = 0, \ldots, n-1\}$. Then

$$I_{\mathcal{L}} = \bigcap_n \bigcup_{\bar{i}} I_{\bar{i}}^n.$$

Let us now transfer these intervals to the unstable lamination. For $f \in \mathcal{A}$, let

$$\mathcal{I}_{\bar{i}}^n(f) = (\chi | \mathcal{W}^u(f))^{-1} I_{\bar{i}}^n.$$

Since the holonomy $\chi : \mathcal{W}^u(f) \to [-2, 0]$ is uniformly quasi-symmetric, the configurations $(\mathcal{W}^u(f), \cup \mathcal{I}_k(f))$ have uniformly bounded geometry (independent of $f \in \mathcal{A}$). But by the Distortion Lemma, the map

$$R^n : (\mathcal{I}_{\bar{i}}^n(f), \bigcup_k \mathcal{I}_{\bar{i},k}^{n+1}(f)) \to (\mathcal{W}^u(R^n f), \bigcup_k \mathcal{I}_k(R^n f))$$

has a bounded distortion. Hence all the configurations $(\mathcal{I}_{\bar{i}}^n(f), \cup_k \mathcal{I}_{\bar{i},k}^{n+1}(f))$ have uniformly bounded geometry, so that the Cantor set $\mathcal{I}(f) = (\chi|\mathcal{W}^u(f))^{-1} I_{\mathcal{L}}$ has bounded geometry. As the holonomy is quasi-symmetric, the desired statement follows. □

Now the HD Theorem stated in the introduction follows.

## 10. Appendix 1: Quasi-conformal maps

The material of this appendix is standard in analysis and dynamics. We add it in order to fix some terminology and notation and to provide for the reader's convenience some basic references.

10.1. *Quasi-conformal maps.* The reader can consult [A], [LV] for the basic theory of quasi-conformal maps.

Let $U, V$ stand for domains in $\mathbb{C}$. We say that a continuous map $f : U \to \mathbb{C}$ belongs to (Sobolev) class $H$ if it has locally square integrable distributional derivatives $\partial h$, $\bar{\partial} h$. A homeomorphism $h : U \to V$ is called quasi-conformal (qc) if it belongs to $H$ and $|\bar{\partial} h / \partial h| \leq k < 1$ almost everywhere. As this local definition is conformally invariant, one can define qc homeomorphisms between Riemann surfaces.

One associates to a qc map an analytic object called a Beltrami differential, namely

$$\mu = \frac{\bar{\partial} h}{\partial h} \frac{d\bar{z}}{dz},$$

with $\|\mu\|_\infty < 1$. The corresponding geometric object is a measurable family of infinitesimal ellipses (defined up to dilation), pull-backs by $h_*$ of the field of infinitesimal circles. The eccentricities of these ellipses are ruled by $|\mu|$, and are uniformly bounded almost everywhere, while the orientation of the ellipses is ruled by the $\arg \mu$. The big axes of these ellipses determine a *line field* on the measurable support of the differential. The dilatation $\mathrm{Dil}(h) \equiv K_h = (1 + \|\mu\|_\infty)/(1 - \|\mu\|_\infty)$ of $h$ is the essential supremum of the eccentricities of these ellipses. A qc map $h$ is called $K$-qc if $\mathrm{Dil}(h) \leq K$.

A remarkable fact is that *any Beltrami differential with $\|\mu_\infty\| < 1$ (or rather a measurable field of ellipses with essentially bounded eccentricities) is* *locally generated by a qc map*, unique up to post-composition with a conformal map. Thus such a Beltrami differential on a Riemann surface $S$ induces a conformal structure quasi-conformally equivalent to the original structure of $S$. Together with the Riemann mapping theorem this leads to the following statement:

THEOREM 10.1 (Measurable Riemann Mapping Theorem). *Let $\mu$ be a Beltrami differential on $\bar{\mathbb{C}}$ with $\|\mu_\infty\| < 1$. Then there is a unique quasi-conformal map $h = h_\mu : \bar{\mathbb{C}} \to \bar{\mathbb{C}}$ which solves the Beltrami equation: $|\bar{\partial} h/\partial h| = \mu$, and is normalized at three points. Moreover, $h_\mu$ holomorphically depends on $\mu$.*

The last statement means that $h_\mu(z)$ holomorphically depends on $\mu$ for every given $z$ (note that $\mu$ is an element of the unit Banach ball of $L^\infty$ which has a natural complex structure) – see [AB] for a thorough discussion.

In what follows, by a conformal structure we will mean a structure associated to a measurable Beltrami differential $\mu$ with $\|\mu\|_\infty < 1$. We will denote by $\sigma$ the standard structure corresponding to the zero Beltrami differential.

Another fundamental property of qc maps is the following:

THEOREM 10.2 (Compactness). *The space of $K$-qc maps $h : \mathbb{C} \to \mathbb{C}$ normalized by $h(0) = 0$ and $h(1) = 1$ is compact in the uniform topology on the Riemann sphere.*

COROLLARY 10.3. *Let $h : (\mathbb{C}, 0, 1) \to (\mathbb{C}, 0, 1)$ be a $K$-qc map. Then for $\varepsilon > 0$, $hD(1, \varepsilon) \supset D(1, \delta)$ with $\delta = \delta(K, \varepsilon) > 0$.*

We will also make use of the following properties:

LEMMA 10.4 (Gluing). *Consider a compact set $Q \subset \mathbb{C}$, two of its neighborhoods $U$ and $V$, and two maps $\phi : U \to \mathbb{C}$ and $\psi : V \setminus Q \to \mathbb{C}$ of class $H$. Assume that these maps match on $\partial Q$, i.e., the map $f : V \to \mathbb{C}$ defined as $\phi$ on $Q$ and as $\psi$ on $V \setminus Q$ is continuous. Then $f \in H$ and the distributional derivatives of $f$ on $Q$ are equal to the corresponding derivatives of $\phi$. In particular, if $\phi$ and $\psi$ are qc homeomorphisms, then $f$ is a qc homeomorphism and $\mathrm{Dil}(f) = \max(\mathrm{Dil}(\phi|Q), \mathrm{Dil}(\psi))$.*

See e.g., [DH2, Lemma 2, p. 303] for a proof (where the lemma is stated for qc homeomorphisms but the proof goes through for general maps of class $H$).

10.2. Qc *classification of quadratic maps.* Let us consider the quadratic family $P_c : z \mapsto z^2 + c$, $c \in \mathbb{C}$, and its Mandelbrot set $M_0 = \{c : J(P_c) \text{ is connected}\}$. Recall that a component $H$ of $\mathrm{int} M_0$ is called *hyperbolic* if the maps $P_c$, $c \in H$, have an attracting cycle. Any hyperbolic component contains a unique superattracting parameter value $c_H$ called its *center*. Nonhyperbolic components of $\text{int} M_0$ are called *queer* (conjecturally they do not exist).

THEOREM 10.5. *The quadratic family is decomposed into the following quasi-conformal classes*:

(i) *the complement of the Mandelbrot set*, $\mathbb{C} \setminus M_0$;

(ii) *punctured hyperbolic components* $H \setminus c_H$;

(iii) *queer components*;

(iv) *single points*.

*The holomorphic deformations in the qc classes can be obtained via holomorphic motions.*

In case (i), the qc deformation is obtained by changing the position of the critical value. In case (ii) it is obtained by changing the multiplier of the attracting point. In case (iii) the deformation is generated by an invariant line field on the Julia set. The last statement says that points on the $\partial M_0$ (and of course the centers of hyperbolic components) are qc rigid. In particular, they do not admit invariant line fields on the Julia set. See [L3], [McM3] for further discussion and references.

## 11. Appendix 2: Complex structures modeled on families of Banach spaces

11.1. *Complex analysis on Banach manifolds.* We assume familiarity with the standard theory of manifolds modeled on Banach spaces (see e.g., [D1], [Lang]). Below we will state a few facts which are specifically complex analytic.

Given a Banach space $\mathcal{B}$, let $\mathcal{B}_r(x)$ stand for the ball of radius $r$ centered at $x$ in $\mathcal{B}$, and $\mathcal{B}_r \equiv \mathcal{B}_r(0)$.

THE CAUCHY INEQUALITY. *Let* $f : (\mathcal{B}_1, 0) \to (\mathcal{D}_1, 0)$ *be a complex analytic map between two unit Banach balls. Then* $\|\mathrm{D}f(0)\| \leq 1$. *Moreover, for* $x \in \mathcal{B}_1$,
$$\|\mathrm{D}f(x)\| \leq \frac{1}{1 - \|x\|}.$$

*Proof* (Yuri Lyubich). Take a vector $v \in \mathcal{B}$ with $\|v\| = 1$ and a linear functional $\psi$ on $\mathcal{D}$ with $\|\psi\| = 1$. Consider an analytic function $\phi : \mathbb{D}_1 \to \mathbb{D}_1$, $\phi(\lambda) = \psi(f(\lambda v))$. As $|\phi(\lambda)| < 1$, the usual Cauchy Inequality yields: $|\phi'(0)| = |\psi(\mathrm{D}f(0)v)| \leq 1$. Since this holds for any normalized $v$ and $\psi$, the former estimate follows by the Hahn-Banach Theorem.

The latter one follows from the former by restricting $f$ to the ball $\mathcal{B}_{1-\|x\|}(x)$. □

The Cauchy Inequality yields:

THE SCHWARZ LEMMA. *Let $r < 1/2$ and $f : (\mathcal{B}_1, 0) \to (\mathcal{D}_r, 0)$ be a complex analytic map between two Banach balls. Then the restriction of $f$ onto the ball $\mathcal{B}_r$ is contracting: $\|f(x) - f(y)\| \leq q\|x - y\|$, where $q = r/(1-r) < 1$.*

*Proof.* By the Cauchy Inequality, $\|Df(x)\| \leq q$ for $x \in \mathcal{B}_r$. Integrating this along the interval $[x, y]$, we obtain the conclusion. □

Let us now state a couple of facts on the intersection properties between analytic submanifolds which provide us a tool to obtain the transversality results.

Let $\mathcal{X}$ and $\mathcal{S}$ be two submanifolds in the Banach space $\mathcal{B}$ intersecting at point $x$. Assume that $\operatorname{codim} \mathcal{X} = \dim \mathcal{S} = 1$. Let us define the *intersection multiplicity* $\sigma$ between $\mathcal{X}$ and $\mathcal{S}$ at $x$ as follows. Select a local coordinate system $(w, z)$ near $x$ in such a way that $x = 0$ and $\mathcal{X} = \{z = 0\}$. Let us analytically parametrize $\mathcal{S}$ near 0: $z = z(t), w = w(t), z(0) = 0, w(0) = 0$. Then by definition, $\sigma$ is the multiplicity of the root of $z(t)$ at $t = 0$.

THE HURWITZ THEOREM. *Under the above circumstances, consider a submanifold $\mathcal{Y}$ of codimension 1 obtained by a small perturbation of $\mathcal{X}$. Then $\mathcal{S}$ has $\sigma$ intersection points with $\mathcal{Y}$ near $x$ counted with multiplicity.*

*Proof.* We use the above local coordinates and parametrization. In these coordinates $\mathcal{Y}$ is a graph of a holomorphic function $z = \phi(w)$ which is uniformly small at some neighborhood of 0 (this is the meaning of $\mathcal{Y}$ being a small perturbation of $\mathcal{X}$). The intersection points of $\mathcal{Y}$ and $\mathcal{S}$ are the roots of the equation $z(t) = \phi(w(t))$. By the classical Hurwitz Theorem, this equation has exactly $\sigma$ roots near the origin, counted with multiplicities if $\phi$ is small enough. □

As usual, a foliation of some analytic Banach manifold is called analytic (smooth) if it can be locally straightened by an analytic (smooth) change of variable.

THE INTERSECTION LEMMA. *Let $\mathcal{F}$ be a codimension-one complex analytic foliation in a domain of a Banach space. Let $\mathcal{S}$ be a one-dimensional complex analytic submanifold intersecting a leaf $\mathcal{L}_0$ of the foliation at a point $x$ with multiplicity $\sigma$. Then $\mathcal{S}$ has $\sigma$ simple intersection points with any nearby leaf.*

*Proof.* Let us select local coordinates $(w, z)$ near $x$ so that $x$ corresponds to 0, and the leaves of the foliation near 0 are given by the equations $\mathcal{L}_\varepsilon = \{z = \varepsilon\}$.

Let $z = z(t), w = w(t)$ be an analytic parametrization of $\mathcal{S}$, with $t = 0$ corresponding to $x = 0$. Then $z(t) = at^\sigma(1 + O(t))$, $a \neq 0$, has a root of multiplicity $\sigma$ at 0. Clearly there is an analytic local chart $\tau = \tau(t)$ in which the curve is parametrized as exact power: $z(\tau) = \tau^\sigma$. Then for small $\varepsilon \neq 0$, the equation $z(\tau) = \varepsilon$ has $\sigma$ simple roots near 0: $\tau_i = \varepsilon^{1/\sigma}$. □

COROLLARY 11.1. *Under the circumstances of the above lemma, $\mathcal{S}$ is transverse to $\mathcal{L}_0$ at $x$ if and only if it has a single intersection point near $x$ with all nearby leaves.*

Let $X \subset \mathbb{C}$ be a subset of the complex plane. A *holomorphic motion* of $X$ over a Banach domain $(\Lambda, 0)$ is a a family of injections $h_\lambda : X \to \mathbb{C}$, $\lambda \in \Lambda$, with $h_0 = \mathrm{id}$, holomorphically depending on $\lambda \in \mathcal{B}_1$ (for any given $z \in X$). The graphs of the functions $\lambda \mapsto h_\lambda(z)$, $z \in X$, form a foliation $\mathcal{F}$ (or rather a lamination as it is partially defined) in $\Lambda \times \mathbb{C}$ with complex codimension-one analytic leaves. This is a "dual" viewpoint on holomorphic motions.

We will now state a basic fact about holomorphic motions usually referred to as the "λ-lemma". It consists of two parts: extension and quasi-conformality which will be stated separately. The consecutively improving versions of the Extension Lemma appeared in [L1], [MSS], [ST], [BR], [Sl]. The final result is due to Slodkowski:

THE λ-LEMMA (EXTENSION). *A holomorphic motion $h_\lambda : X_* \to X_\lambda$ of a set $X_* \subset \mathbb{C}$ over a topological disc $D$ admits an extension to a holomorphic motion $H_\lambda : \mathbb{C} \to \mathbb{C}$ of the whole complex plane over $D$.*

The point of the following simple lemma as compared with the previous deep one is smoothness of the extension. The parameter space is allowed to be infinitely dimensional.

LEMMA 11.2 (Local extension). *Let us consider a compact set $Q \subset \mathbb{C}$ and a smooth holomorphic motion $h_\lambda$ of a neighborhood $U$ of $Q$ over a Banach domain $(\Lambda, 0)$. Then there is a smooth holomorphic motion $H_\lambda$ of the whole complex plane $\mathbb{C}$ over some neighborhood $\Lambda' \subset \Lambda$ of 0 whose restriction to $Q$ coincides with $h_\lambda$.*

*Proof.* We can certainly assume that $\bar{U}$ is compact. Take a smooth function $\phi : \mathbb{C} \to \mathbb{R}$ supported in $U$ and let

$$H_\lambda = \phi \, h_\lambda + (1 - \phi)\mathrm{id}.$$

Clearly $H$ is smooth in both variables, holomorphic in $\lambda$, and identical outside $U$. As $H_0 = \mathrm{id}$, $H_\lambda : \mathbb{C} \to \mathbb{C}$ is a diffeomorphism for $\lambda$ sufficiently close to 0, and we are done. □

Given two complex one-dimensional transversals $\mathcal{S}$ and $\mathcal{T}$ to the lamination $\mathcal{F}$ in $\mathcal{B}_1 \times \mathbb{C}$, we have a holonomy $\mathcal{S} \to \mathcal{T}$. We say that this map is locally quasi-conformal if it admits local quasi-conformal extensions near any point.

Given two points $\lambda, \mu \in \mathcal{B}_1$, let us define the hyperbolic distance $\rho(\lambda, \mu)$ in $\mathcal{B}_1$ as the hyperbolic distance between $\lambda$ and $\mu$ in the one-dimensional complex slice $\lambda + t(\mu - \lambda)$ passing through these points in $\mathcal{B}_1$.

THE $\lambda$-LEMMA (QUASI-CONFORMALITY). *Holomorphic motion $h_\lambda$ of a set $X$ over a Banach ball $\mathcal{B}_1$ is transversally quasi-conformal. The local dilatation $K$ of the holonomy from $p = (\lambda, u) \in \mathcal{S}$ to $q = (\mu, v) \in \mathcal{T}$ depends only on the hyperbolic distance $\rho$ between the points $\lambda$ and $\mu$ in $\mathcal{B}_1$. Moreover, $K = 1 + O(\rho)$ as $\rho \to 0$.*

*Proof.* If the transversals are vertical lines $\lambda \times \mathbb{C}$ and $\mu \times \mathbb{C}$ then the result follows from the classical $\lambda$-lemma [MSS] by restriction of the motion to the complex line joining $\lambda$ and $\mu$.

Furthermore, the holonomy from the vertical line $\lambda \times \mathbb{C}$ to the transversal $\mathcal{S}$ is locally conformal at point $p$. To see this, we select holomorphic coordinates $(\theta, z)$ near $p$ in such a way that $p = 0$ and the leaf via $p$ becomes the parameter axis. Let $z = \psi(\theta) = \varepsilon + \ldots$ parametrize a nearby leaf of the foliation, while $\theta = g(z) = bz + \ldots$ parametrizes the transversal $\mathcal{S}$.

Let us do the rescaling $z = \varepsilon\zeta, \theta = \varepsilon\nu$. In these new coordinates, the above leaf is parametrized by the function $\Psi(\nu) = \varepsilon^{-1}\psi(\varepsilon\nu)$, $|\nu| < R$, where $R$ is a fixed parameter. Then $\Psi'(\nu) = \psi'(\varepsilon\nu)$ and $\Psi''(\nu) = \varepsilon\psi''(\varepsilon\nu)$. By the Cauchy Inequality, $\Psi''(\nu) = O(\varepsilon)$. Moreover, $\psi$ uniformly goes to 0 as $\psi(0) \to 0$. Hence $|\Psi'(0)| = |\psi'(0)| \leq \delta_0(\varepsilon)$, where $\delta_0(\varepsilon) \to 0$ as $\varepsilon \to 0$. Thus $\Psi'(\nu) = \delta_0(\varepsilon) + O(\varepsilon) \leq \delta(\varepsilon) \to 0$ as $\varepsilon \to 0$ uniformly for all $|\nu| < R$. It follows that $\Psi(\nu) = 1 + O(\delta(\varepsilon)) = 1 + o(1)$ as $\varepsilon \to 0$.

On the other hand, the manifold $\mathcal{S}$ is parametrized in the rescaled coordinates by a function $\nu = b\zeta(1 + o(1))$. Since the transverse intersection persists, $\mathcal{S}$ intersects the leaf at the point $(\nu_0, \zeta_0) = (b, 1)(1 + o(1))$ (so that $R$ should be selected bigger than $\|b\|$). In the old coordinates the intersection point is $(\theta_0, z_0) = (\varepsilon, b\varepsilon)(1 + o(1))$.

Thus the holonomy from $\lambda \times \mathbb{C}$ to $\mathcal{S}$ transforms the disc of radius $|\varepsilon|$ to an ellipse with small eccentricity, which means that this holonomy is asymptotically conformal. As the holonomy from $\mu \times \mathbb{C}$ to $\mathcal{T}$ is also asymptotically conformal, the conclusion follows. □

Quasi-conformality is apparently the best regularity of holomorphic motions which is satisfied automatically. A holomorphic motion is called *smooth* (or *real analytic* etc.) if it has the corresponding regularity in *both* variables.

11.2. *Inductive limits.* Let $(\mathbb{V}, \succ)$ be a partially ordered set. In this section we assume that $\mathbb{V}$ is *directed* i.e., any two elements have a common majorant. We assume that $\mathbb{V}$ has a countable base; i.e., there is a countable subset $\mathbb{W} \subset \mathbb{V}$ such that any $V \in \mathbb{V}$ has a majorant $W \in \mathbb{W}$.

Let us consider a family of Banach spaces $\mathcal{B}_V$ labeled by the elements of $\mathbb{V}$. An $\varepsilon$-ball in $\mathcal{B}_V$ centered in an $f \in \mathcal{B}_V$ will be denoted $\mathcal{B}_V(f, \varepsilon)$. Elements of the $\mathcal{B}_V$ will be called "maps"(keep in mind further applications to quadratic-like maps). For every pair $U \succ V$, we have a continuous linear *injection* $i_{U,V} : \mathcal{B}_V \to \mathcal{B}_U$. We assume the following properties:

C1. *Density*: the image $i_{U,V}\mathcal{B}_V$ is dense in $\mathcal{B}_U$;

C2. *Compactness*: the map $i_{U,V}$ is compact; i.e., the images of balls, $i_{U,V}\mathcal{B}_V(f, R)$, are pre-compact in $\mathcal{B}_U$.

LEMMA 11.3.

- If $U, W \succ V$, $f \in \mathcal{B}_V$, $R > 0$, then the metrics $\rho_U$ and $\rho_W$ induced on the ball $\mathcal{B}_V(f, R)$ from $\mathcal{B}_U$ and $\mathcal{B}_W$ are equivalent.

- Let $U \succ V$, and $\phi_i : (\mathcal{B}_U, \mathcal{B}_V) \to (\mathbb{C}, \mathbb{C})$ be a family of linear functionals continuous on the both spaces. Let us consider the common kernels of these functionals in the corresponding spaces: $\mathcal{L}_U \subset \mathcal{B}_U$ and $\mathcal{L}_V = \mathcal{L}_U \cap \mathcal{B}_V$. Then $\operatorname{codim}(\mathcal{L}_U | \mathcal{B}_U) = \operatorname{codim}(\mathcal{L}_V | \mathcal{B}_V)$.

*Proof.* • It is clearly enough to check the case when $W \succ U$. Assume that there exists a sequence $f_n \in \mathcal{B}_V(0, R)$ such that $\|f_n\|_W \to 0$ while $\|f_n\|_U$ stay bounded away from 0. By compactness of $i_{U,V}$, we can pass to a limit $f_n \to f$ in $\mathcal{B}_U$ along a subsequence. Then $f \neq 0$, while $i_{W,U} f = 0$ contradicting that $i_{W,U}$ is injective.

• If a family of functional as above is linearly dependent on $\mathcal{B}_V$ then by the density property it is linearly dependent on $\mathcal{B}_U$ as well. This yields the second statement. □

For any $W \succ V$, let us identify any $f \in \mathcal{B}_V$ with its image $i_{W,V} f \in \mathcal{B}_W$ and span the equivalence relation generated by these identifications. Thus $f \in \mathcal{B}_V$ and $g \in \mathcal{B}_U$ are equivalent if there is a common majorant $W \succ (U, V)$ such that $i_{W,V} f = i_{W,U} g$ (then by injectivity this holds for any common majorant). The equivalence classes will be called *germs*. The space of germs is called the *inductive limit* of the Banach spaces $\mathcal{B}_V$ and is denoted by $\mathcal{B} = \lim \mathcal{B}_V$.

Every space $\mathcal{B}_V$ is naturally injected into the space of germs, and will be considered as a subset of the latter. Given a subset $\mathcal{X} \subset \mathcal{B}$, the intersection $\mathcal{X}_V \equiv \mathcal{X} \cap \mathcal{B}_V$ will be called a *(Banach) slice* of $\mathcal{X}$.

Let us supply $\mathcal{B}$ with the *inductive limit topology*. In this topology, a set $\mathcal{X} \subset \mathcal{B}$ is claimed to be open if all its Banach slices $\mathcal{X}_V$ are open. The axioms of topology are obviously satisfied, and the linear operations are obviously continuous (note that the product topology on $\mathcal{B} \times \mathcal{B}$ coincides with the natural inductive limit topology). Thus $\mathcal{B}$ is a topological vector space. Since points are obviously closed in this topology, $\mathcal{B}$ is Hausdorff (see [Ru]). The following lemma summarizes some useful general properties of the inductive limits.

LEMMA 11.4.

(i) *In the inductive limit topology, $f_n \to f$ if and only if all the maps $f_n$ and $f$ belong to the same Banach slice $\mathcal{B}_V$ and $f_n \to f$ in the intrinsic topology of $\mathcal{B}_V$. Any cluster point $f$ of a set $K \subset \mathcal{B}$ is a limit of a sequence $\{f_n\} \subset K$.*

(ii) *A set $\mathcal{X} \subset \mathcal{B}$ is open if and only if it is sequentially open.*

(iii) *If $X$ is a metric space and $\phi : (X, a) \to (\mathcal{B}, g)$ is a continuous map then there is a neighborhood $D \ni a$ and an element $V \in \mathbb{V}$ such that $\phi D \subset \mathcal{B}_V$.*

(iv) *A set $\mathcal{K} \subset \mathcal{B}$ is compact if and only if it is sequentially compact. Moreover, $\mathcal{K}$ sits in some Banach space $\mathcal{B}_V$ such that the induced metric on $\mathcal{K}$ is compatible with its topology.*

(v) *A map $\phi : \mathcal{B} \to T$ to a topological space $T$ is continuous if and only if every restriction $\phi|\mathcal{B}_V$ is continuous. The map $\phi$ is continuous if and only if it is sequentially continuous.*

*Proof.* (i) Since the inclusions $\mathcal{B}_V \to \mathcal{B}$ are continuous, any sequence $\{f_n\} \subset \mathcal{B}_V$ converging to $f$ in $\mathcal{B}_V$ converges to $f$ in $\mathcal{B}$ as well.

Let us assume that $\{f_n\}$ converges to $f$ in $\mathcal{B}$ but does not sit in any Banach slice. Then we can select a subsequence which hits any Banach slice at most finitely many times and never hits $f$ itself. By definition of the inductive limit topology, the complement of this sequence is a neighborhood of $f$ – a contradiction.

Similarly, if $\{f_n\}$ is not bounded in any Banach slice, then we can select a subsequence $f_{n(k)}$ such that $\|f_{n(k)}\|_{U_i} \geq i$ for $i = 1, \ldots, k$, where $\{U_k\}$ is a countable base in $\mathbb{V}$. This subsequence has a discrete intersection with any Banach space $\mathcal{B}_{U_k}$ (in the Banach topology), and hence $\mathcal{B} \setminus \{f_{n(k)}\}$ is open – a contradiction.

Thus the whole sequence $\{f_n\}$ sits in some Banach space $\mathcal{B}_V$ and is bounded there. Hence it is compact in any $B_U$ with $U \succ V$. But any limit point $g \in \mathcal{B}_U$ of $\{f_n\}$ must coincide with $f$. Hence $f_n \to f$ in the Banach topology of $\mathcal{B}_U$.

If the latter statement concerning cluster points fails then $f$ is not a cluster point for the slices $\mathcal{K}_V$. Then we can construct a neighborhood $\mathcal{U} \subset \mathcal{B}$ of $f$ missing $\mathcal{K}$ in the same way as above.

(ii) Generally, any open set is sequentially open (which means that its complement is closed with respect to taking limits of converging sequences). Vice versa, if a set is sequentially open, then clearly its Banach slices are sequentially open. By definition, the set itself is open.

(iii) Otherwise there would be a sequence $x_n \to a$ such that the maps $f_n = \phi(x_n)$ do not sit in a common space $\mathcal{B}_V$ despite the fact that $f_n \to g$.

(iv) Let us show that any compact set $\mathcal{K} \subset \mathcal{B}$ sits in some Banach slice. Otherwise there would be a sequence $\{f_n\} \subset \mathcal{K}$, no subsequence of which sits in a common Banach slice (since $\mathbb{V}$ has a countable base). But by the first point of this lemma, such a subsequence has no cluster points. Similarly one can see that $\mathcal{K}$ is a bounded subset in some Banach slice $\mathcal{B}_V$.

Let $W \succ V$. Then $\mathcal{K}$ is compact in the topology of $\mathcal{B}_W$. Since this topology is finer than the inductive limit topology, they must coincide on $\mathcal{K}$.

Exactly the same argument can be applied to sequentially compact sets. Since compactness and sequential compactness are equivalent in metric spaces, the desired statement follows.

(v) Take an open set $X \subset \mathcal{T}$. By the definition of topology, the preimage $\phi^{-1}X$ is open if and only if all its slices $(\phi|\mathcal{B}_V)^{-1}X$ are open. This implies the former statement, which yields the latter. □

The above metrics on compact sets defined in (iv) will be called *Montel metrics* $\mathrm{dist}_M$. (We will not specify the Banach space $\mathcal{B}_V$ from which the metric is induced.)

*Remarks.* 1. Any continuous curve $\gamma : (\mathbb{R}, 0) \to (\mathcal{B}, g)$ locally sits in some space $\mathcal{B}_V$.

2. Given a continuous transformation $R : (\mathcal{B}, f) \to (\mathcal{T}, g)$ between two spaces of germs over $\mathbb{V}$ and $\mathbb{U}$ respectively, for any $V \in \mathbb{V}$ there exist an $\varepsilon > 0$ and an element $U \in \mathbb{U}$ such that $R(\mathcal{B}_V(f, \varepsilon)) \subset \mathcal{B}_U$.

3. The third statement of the above lemma shows that $\mathcal{B}$ is not a Freshe space, i.e., it is not metrizable, and thus does not have a local countable base of neighborhoods. However, as we see, the sequential description of basic topological properties (cluster points, compactness, continuity etc.) is adequate.

4. Note that the Banach slices $\mathcal{B}_V$ are dense in the space of germs $\mathcal{B}$. Thus their intrinsic topology is not induced from $\mathcal{B}$.

Let us define a *sublimit* of the directed family $\mathcal{B}_V$, $V \in \mathbb{V}$, as the inductive limit of Banach spaces $\mathcal{B}_V$ corresponding to a directed subset $\mathbb{U} \subset \mathbb{V}$ (which is not necessarily exhausting).

All linear operators $A : \mathcal{B} \to \mathcal{T}$ between spaces of germs are assumed to be continuous. Let us supply this space with the following convergence structure. A sequence of linear operators $A_n : \mathcal{B} \to \mathcal{T}$ converges to an operator $A$ if for any $V \in \mathbb{V}$ and $W, U \in \mathbb{U}$, $W \succ V$, such that $A(\mathcal{B}_V) \subset \mathcal{T}_U$, we have: $A_n(\mathcal{B}_V) \subset \mathcal{T}_W$ for all sufficiently big $n$ and the restrictions $A_n : \mathcal{B}_V \to \mathcal{B}_W$ converge to $A : \mathcal{B}_V \to \mathcal{B}_W$ in the uniform operator topology.

11.3. *Main example: analytic germs.* Let $\mathbb{V}$ be the directed set of topological discs $V \ni 0$ with piecewise smooth boundary, with $U \succ V$ if $U \Subset V$. Let $\mathcal{B}_V$ denote the affine Banach space of normalized analytic functions of the form $f(z) = c + z^2 + \ldots$ on $V \in \mathbb{V}$ continuous up to the boundary supplied with sup-norm $\|\cdot\|_V$.

*Remark.* To make this example consistent with the previous discussion, one can make $\mathcal{B}_V$ linear by putting the origin into the map $f(z) = z^2$. Or one can rather make the previous discussion in the category of affine Banach spaces.

For $U \succ V$, define the injection $i_{U,V} : \mathcal{B}_V \to \mathcal{B}_U$ by restricting the functions. Since polynomials are dense in $\mathcal{B}_U$, this inclusion has a dense image. Moreover, by Montel's Theorem, the balls of $\mathcal{B}_V$ are relatively compact in $\mathcal{B}_U$. Thus this family of Banach spaces satisfies assumptions C1–C2 from subsection 11.2, so that we can form the inductive limit $\mathcal{B} = \lim \mathcal{B}_V$. The elements of this space are *analytic germs* at 0.

Let us say that two metrics $\rho$ and $d$ on the same space $\mathcal{K}$ are *Hölder equivalent* if there exist constants $C > 0$ and $\delta > 0$ such that

$$C^{-1}\rho(x,y)^{1/\delta} \leq d(x,y) \leq C\rho(x,y)^{\delta}.$$

The following classical statement is a version of the Hadamard Three Circles Theorem.

LEMMA 11.5. *Consider three domains $V \Subset W \Subset U$. Then the metrics $\|\cdot\|_V$ and $\|\cdot\|_W$ induced on the unit Banach ball of $\mathcal{B}_U$ are Hölder equivalent. Moreover, the Hölder exponent goes to 1 as $V \to W$ in the Hausdorff metric.*

*Proof.* Take a holomorphic function $f$ on $U$ with $\|f\|_U \leq 1$. Let $\|f\|_V = \varepsilon$.

Let us consider a positive harmonic function $h$ on the annulus $U \setminus V$ with boundary values 0 and 1 on its outer and inner boundaries respectively. Then

(11.1) $$\log|f| \leq h \log \varepsilon$$

on the boundary of the annulus. Since $\log |f|$ is subharmonic, (11.1) also holds inside the annulus. Putting $\delta = \inf_{z \in \partial W} h(z)$, we conclude that $|f|_W \leq \varepsilon^\delta = \|f\|_V^\delta$.

Moreover, it is clear from the above formula for the exponent $\delta$ that $\delta$ is close to 1 if $V$ is Hausdorff close to $W$. □

Thus the Montel metric $\text{dist}_M$ on compact sets of germs is well defined up to Hölder equivalence. In other words, compact subsets $\mathcal{K} \subset \mathcal{B}$ bear a natural Hölder structure.

11.4. *Analytic maps.* Let us consider an inductive limit $\mathcal{B}$ over $\mathbb{V}$. By definition, a function $\phi : \mathcal{B} \to \mathbb{C}$ is complex analytic if all the restrictions $\phi|\mathcal{B}_V$ are complex analytic in the Banach sense.

Let us consider a continuous map $R : \mathcal{V} \to \mathcal{B}'$, where $\mathcal{V}$ is an open subset of $\mathcal{B}$ and $\mathcal{B}'$ is an inductive limit space over $\mathbb{V}'$. It is called *differentiable* at a point $f \in \mathcal{B}$ if there is a real linear operator $A \equiv \mathrm{D}R(f) : \mathcal{B} \to \mathcal{B}'$ such that any Banach restriction $R_V : \mathcal{B}_V \to \mathcal{B}'_U$ with $f \in \mathcal{B}_V$ is differentiable at $f$, and $\mathrm{D}R_V(f) = A|\mathcal{B}_V$.

As usual, a map $R : \mathcal{V} \to \mathcal{B}'$ is called *smooth* if it is differentiable at every point $f \in \mathcal{V}$ and the differential $\mathrm{D}R(f)$ is (sequentially) continuous in $f$ (which amounts to the smoothness of all Banach restrictions). A map $R : \mathcal{V} \to \mathcal{B}'$ is called *analytic* if it is smooth, and the differentials $\mathrm{D}R(f)$ are linear over $\mathbb{C}$.

11.5. *Varieties.* Let us have a family of Banach spaces $\mathcal{B}_V$ labeled by elements $V$ of some set $\mathbb{V}$, and open sets $\mathcal{U}_V \subset \mathcal{B}_V$. Let us have a set $\mathcal{QG}$ and a family of injections $j_V : \mathcal{U}_V \to \mathcal{QG}$. The images $\mathcal{S}_V \equiv j_V \mathcal{U}_V$ are called Banach slices in $\mathcal{QG}$. The images $j_V \mathcal{V}_V \subset \mathcal{S}_V$ of open sets $\mathcal{V}_V \subset \mathcal{U}_V$ are called *Banach neighborhoods*. We assume the following properties (compare with C1 and C2):

P1: *countable base and compactness.* There is a countable family of slices $\mathcal{S}_i$ with the following property: For any $f \in \mathcal{QG}$ and any slice $\mathcal{S}_V \ni f$, there exists a Banach neighborhood $\mathcal{V}_V \subset \mathcal{S}_V$ compactly contained in some $\mathcal{S}_i$.

P2: *analyticity.* If some Banach neighborhood $j_V \mathcal{V}_V \subset \mathcal{S}_V$ is also contained in another slice $\mathcal{S}_U$, then the transit map $i_{U,V} = j_U^{-1} \circ j_V : \mathcal{V}_V \to \mathcal{U}_U$ is analytic.

P3: *density.* The differential $\mathrm{D}i_{U,V}(f)$ of the above transit map has a dense image in $\mathcal{B}_U$.

We endow $\mathcal{QG}$ with the finest topology which makes all the injections $j_V$ continuous by declaring a set $\mathcal{V} \subset \mathcal{QG}$ open if and only if all its Banach slices $j_V^{-1}\mathcal{V}$ are open. Lemma 11.4 should be modified a little in this more general situation:

LEMMA 11.6.  *In $\mathcal{QG}$, $f_n \to f$ if and only if the sequence $\{f_n\}$ sits in a finite union of the Banach slices, and the corresponding subsequences converge* *to $f$ in the Banach metric. All other statements of Lemma* 11.4 *are valid in $\mathcal{QG}$ as well with the modification that a single Banach slice in* (iii) *and* (iv) *should be replaced with a finite union of Banach slices.*

*Proof.* This is similar to the proof of Lemma 11.4. Let us just comment on the induced metric on a compact set $\mathcal{K} \subset \mathcal{QG}$. Such a set is covered with finitely many Banach balls $\mathcal{B}_i$, $i = 1, \ldots, N$, which bear the Banach $\text{dist}_i$. Let $R > \max \text{diam} \mathcal{B}_i$, where the diameter is evaluated with respect to the corresponding metric.

Given two points $f, g \in \mathcal{K}$, we define $\text{dist}(f, g)$ as follows. If there is a linking sequence $f = f_0, f_1, \ldots, f_n = g$ such that $f_k$ and $f_{k+1}$ belong to the same Banach ball $B_{j(k)}$ then

$$\text{dist}_\text{M}(f, g) = \inf \sum \text{dist}_{j(k)}(f_k, f_{k+1}),$$

where the infimum is taken over all possible linking sequences (note that in this case $\text{dist}(f, g) < RN$). If no linking sequence exists then $\text{dist}_\text{M}(f, g) = RN$. □

Similar to the inductive limit case, the above metrics on compact subsets of $\mathcal{QG}$ will be called Montel metrics.

We say that a topological space $\mathcal{QG}$ as above is a *complex analytic variety* modeled on a family of Banach spaces. A subset $\mathcal{QG}^\#$ will be called a *slice* of $\mathcal{QG}$ if it is a union of some family of Banach neighborhoods $j_V \mathcal{V}_V$. It naturally inherits from $\mathcal{QG}$ complex analytic structure.

Let $\mathbb{V}_f = \{V \in \mathbb{V} : f \in \mathcal{S}_V\}$. Define the *tangent cone* to $\mathcal{QG}$ at $f$ as follows:

$$T_f \mathcal{QG} = \bigsqcup_{\mathcal{V} \in \mathbb{V}_f} \mathcal{B}_\mathcal{V} \backslash \sim,$$

where the equivalence relation $\sim$ is generated by identifications $h \in \mathcal{B}_V$ with $\text{D}i_{UV}(h) \in \mathcal{B}_U$, $U \succ V$. Note that it is generally not a linear space but it is the union of linear (Banach) slices $T_f \mathcal{S}_V \approx \mathcal{B}_V$.

Let us call a point $f \in \mathcal{QG}$ *regular* if any two Banach neighborhoods $\mathcal{U} \subset \mathcal{S}_U$ and $\mathcal{V} \subset \mathcal{S}_V$ around $f$ are contained in a common slice $\mathcal{S}_W$. At such a point the tangent cone $\text{T}_f \mathcal{QG}$ is a linear space identified with the inductive limit of the Banach spaces,

$$\text{T}_f \mathcal{QG} = \lim_{U \in \mathbb{V}_f} \text{T}_f \mathcal{S}_U.$$

Tangent cones at regular points will be called *tangent spaces*. If all points of $\mathcal{QG}$ are regular then it will be called a complex analytic *manifold* (modeled on a family of Banach spaces).

A map $R : \mathcal{QG}^1 \to \mathcal{QG}^2$ is called *analytic* if it locally transfers any Banach slice $\mathcal{S}_U$ to some slice $\mathcal{S}_V$, and its Banach restriction $j_V^{-1} \circ R \circ j_U$ is analytic.

An analytic map has a well-defined differential $DR(f) : T_f \mathcal{QG}^1 \to T_{Rf}\mathcal{QG}^2$ continuously depending on $f$ whose Banach restrictions are linear.

Let $\mathcal{M}$ be a complex analytic manifold modeled on a family of Banach spaces. An analytic map $i : \mathcal{M} \to \mathcal{QG}$ is called *immersion* if for any $m \in \mathcal{M}$ the differential $Di(m)$ is a linear homeomorphism onto its image. The image $\mathcal{X}$ of an injective immersion $i$ is called an *immersed submanifold*. It is called an (*embedded*) *submanifold* if additionally $i$ is a homeomorphism onto $\mathcal{X}$ supplied with the induced topology. For example, if there is an analytic projection $\pi : \mathcal{QG} \to \mathcal{M}$ such that $\pi \circ i = \mathrm{id}$ then $\mathcal{X}$ is a submanifold in $\mathcal{M}$. By definition, the dimension of $\mathcal{X}$ is equal to the dimension of $\mathcal{M}$.

If $i : (\mathcal{M}, m) \to (\mathcal{X}, f) \subset (\mathcal{QG}, f)$ is an immersion, then the tangent space $T_f \mathcal{X}$ is defined as the image of the differential $Di(m)$. If the point $f$ is regular then $T_f \mathcal{X}$ is a linear subspace in $T_f \mathcal{QG}$, so that we have a well-defined notion of the codimension of $\mathcal{X}$ at $f$. Moreover, if $\mathcal{M}$ is a Banach manifold (in particular, a finite-dimensional manifold) then $\mathcal{X}$ locally sits in a Banach slice of $\mathcal{QG}$.

We say that a submanifold $\mathcal{X} \subset \mathcal{QG}$ is regular (of codimension $d$) if all its points $f \in \mathcal{X}$ are regular (and $\mathcal{X}$ has codimension $d$ at all its points).

As usual, two regular submanifolds $\mathcal{X}$ and $\mathcal{Y}$ in $\mathcal{QG}$ are called *transverse* at a point $g \in \mathcal{X} \cap \mathcal{Y}$ if $T_g \mathcal{X} \oplus T_g \mathcal{Y} = T_g \mathcal{QG}$.

State University of New York, Stony Brook, NY
*E-mail address*: mlyubich@math.sunysb.edu

### References

[A]   L. Ahlfors, *Lectures on Quasiconformal Mappings*, Van Nostrand Math. Studies, No. 10, Von Nostrand Co. Inc., New York, 1966.

[AB]  L. Ahlfors and L. Bers, Riemann's mapping theorem for variable metrics, Annals of Math. **72** (1960), 385–404.

[B]   G. D. Birkhoff, Surface transformations and their dynamical applications. Acta Math. **43** (1920), 1–119, in Vol. II of the Collected Mathematical Papers, AMS's 1950 Edition.

[Be]  L. Bers, Spaces of Riemann surfaces as bounded domains, Bull. Amer. Math. Soc. **66** (1960), 98–103.

[Bo]  R. Bowen, Some systems with unique equilibrium states, Math. Systems Theory **8** (1975), 193–202.

[BR]  L. Bers and H. L. Royden, Holomorphic families of injections, Acta Math. **157** (1986), 259–286.

[CE]  P. Collet and J.-P. Eckmann, *Iterated Maps on the Interval as Dynamical Systems*, Birkhäuser, Boston, 1980.

[CEL] P. Collet, J.-P. Eckmann, and O. E. Lanford III, Universal properties of maps on an interval, Comm. Math. Physics **76** (1980), 211–254.

[CG]  L. Carleson and T. W. Gamelin, *Complex Dynamics*, Springer-Verlag, New York, 1993.

[CT]  P. Coullet and C. Tresser, Itérations d'endomorphismes et groupe de renormalisation. J. Phys. Colloque C 539, C5-25 (1978).

[Cv] P. Cvitanović, *Universality in Chaos*, Adam Hilger Ltd., Bristol, 1984.
[D1] A. Douady, Le problème des modules pour sous-espaces analytiques compacts d'un espace analytique donné, Ann. Inst. Fourier **16** (1966), 1–95.
[D2] ———, Chirurgie sur les applications holomorphes, *Proc. ICM*-86, *Berkeley*, A.M.S., Providence, RI (1987), 724–738.
[D3] ———, Does a Julia set depend continuously on the polynomial?, in *Complex Dynamical Systems*, Proc. Symp. Appl. Math. **49**, Amer. Math. Soc., 1994.
[DH1] A. Douady and J. H. Hubbard, Étude dynamique des polynômes complexes, Publication Mathématiques d'Orsay, **84-02**, 1984; **85-04**, 1985.
[DH2] ———, On the dynamics of polynomial-like maps, Ann. Sci. École Norm. Sup. **18** (1985), 287–343.
[DH3] ———, A proof of Thurston's topological characterization of rational functions, Acta Math. **171** (1993), 263–297.
[DGP] B. Derrida, A. Gervois, and Y. Pomeau, Universal metric properties of bifurcations of endomorphisms, J. Phys. A. **12** (1979), 269–296.
[E1] H. Epstein, Fixed points of composition operators II, Nonlinearity **2** (1989), 305–310.
[E2] ———, Fixed points of the period-doubling operator, Lecture notes, Lausanne, 1992.
[EE] J.-P. Eckmann and H. Epstein, Bounds on the unstable eigenvalue for period doubling, Comm. Math. Phys. **128** (1990), 427–435.
[EW] J.-P. Eckmann and P. Wittwer, A complete proof of the Feigenbaum conjectures, J. Statist. Phys. **46** (1987), 455–475.
[F1] M. J. Feigenbaum, Quantitative universality for a class of nonlinear transformations, J. Statist. Phys. **19** (1978), 25–52.
[F2] ———, The universal metric properties of nonlinear transformations, J. Statist. Phys. **21** (1979), 669–706.
[dF] E. de Faria, Proof of universality for critical circle maps, Thesis, CUNY, 1992.
[dFM] E. de Faria and W. de Melo, Rigidity of critical circle maps, I and II, Preprints IMS at Stony Brook, # 1997/16 and # 1997/17.
[GSK] A. I. Golberg, Ya. G. Sinai, and K. M. Khanin, Universal properties of a sequence of period-tripling bifurcations, Russian Math. Surveys **38** (1983), 187–188.
[H] J. H. Hubbard, Local connectivity of Julia sets and bifurcation loci: three theorems of J.-C. Yoccoz, in *Topological Methods in Modern Mathematics*, A Symposium in Honor of John Milnor's 60th Birthday, Publish or Perish, Houston, TX, 1993.
[La1] O. E. Lanford III, A computer-assisted proof of the Feigenbaum conjectures, Bull. Amer. Math. Soc. **6** (1982), 427–434.
[La2] ———, Renormalization group methods for circle mappings, *Nonlinear Evolution and Chaotic Phenomena*, NATO Adv. Sci. Inst., Ser. B, Phys. **176** (1988), 25–36.
[Lang] S. Lang, *Differentiable Manifolds*, Addison-Wesley Publishing Co., Inc., Reading, MA, 1972.
[L1] M. Lyubich, Some typical properties of the dynamics of rational maps, Russian Math. Surveys **38** (1983), 154–155.
[L2] ———, Dynamics of quadratic polynomials, I,II, Acta Math. **178**(1997), 185–297.
[L3] ———, Renormalization ideas in conformal dynamics, *Current Developments in Mathematics*, 1995, (eds.: R. Bott et al.), International Press, Cambridge, MA, 1994, 155–190.
[L4] ———, Dynamics of quadratic polynomials, III: Parapuzzle and SBR measures, preprint IMS at Stony Brook (http://www.math.sunysb.edu), #1996/5, to appear in Astérisque.
[L5] ———, Almost every real quadratic map is either regular or stochastic, preprint IMS at Stony Brook (http://www.math.sunysb.edu), #1997/8.

[LY]     M. Lyubich and M. Yampolsky, Dynamics of quadratic polynomials: complex bounds for real maps, Ann. Inst. Fourier. **47** (1997), 1219–1255.
[LS]     G. Levin and S. van Strien, Local connectivity of the Julia set of real polynomials, Annals of Math. **147** (1998), 471–541.
[LV]     O. Lehto and K. I. Virtanen, *Quasiconformal Mappings in the Plane*, Springer-Verlag, New York, 1973.
[M]      J. Milnor, Self-similarity and hairiness in the Mandelbrot set, in *Computers in Geometry and Topology*, Lect. Notes in Pure and Appl Math. **114** (1989), 211–257.
[MT]     J. Milnor and W. Thurston, On iterated maps of the interval, *Dynamical Systems*, LNM **1342** Springer-Verlag, New York (1988), 465–563.
[McM1]   C. McMullen, *Complex Dynamics and Renormalization*, Annals of Math. Studies **135**, Princeton University Press, Princeton, NJ, 1994.
[McM2]   ———, *Renormalization and* 3-*manifolds which Fiber over the Circle*, Annals of Math. Studies **142**, Princeton University Press, Princeton, NJ, 1996.
[McM3]   ———, The classification of conformal dynamical systems, *Current Developments in Mathematics*, 1995, International Press, Cambridge, MA, 1994 (eds. R. Bott et al.), 323–360.
[Ma]     M. Martens, The periodic points of renormalization, Annals of Math. **147** (1998), 543–584.
[MSS]    R. Mañé, P. Sad, and D. Sullivan, On the dynamics of rational maps, Ann. Sci. École Norm. Sup. **16** (1983), 193–217.
[MvS]    W. de Melo and S. van Strien, *One-dimensional Dynamics*, Springer-Verlag, New York, 1993.
[PM]     R. Perez-Marco, Fixed points and circle maps, Acta Math. **179** (1997), 243–294.
[Ru]     W. Rudin, *Functional Analysis*, McGraw-Hill, Inc. New York, 1991.
[Sh]     M. Shub, *Global Stability of Dynamical Systems*, Springer-Verlag, New York, 1968.
[S1]     D. Sullivan, Quasiconformal homeomorphisms in dynamics, topology and geometry, Proc. ICM-86, Berkeley **II**, A.M.S., Providence, RI (1987), 1216–1228.
[S2]     ———, Bounds, quadratic differentials, and renormalization conjectures, AMS Centennial Publications **II**, *Mathematics into Twenty-first Century*, 417–466, 1992.
[ST]     D. Sullivan and W. Thurston, Extending holomorphic motions, Acta. Math. **157** (1986), 243–257.
[Sl]     Z. Slodkowski, Holomorphic motions and polynomial hulls, Proc. Amer. Math. Soc. **111**(1991), 347–355.
[Sch]    D. Schleicher, Structure of the Mandelbrot set, preprint, 1995.
[SN]     L. Sario and M. Nakai, *Classification Theory of Riemann surfaces*, Springer-Verlag, New York, 1970.
[T]      C. Tresser, Aspects of renormalization in dynamical systems theory, in *Chaos and Complexity*, 11–19, Editions Frontiers, 1995.
[TC]     C. Tresser and P. Coullet, Itérations d'endomorphismes et groupe de renormalisation, C. R. Acad. Sci. Paris **287A** (1978), 577–580.
[VSK]    E. B. Vul, Ya. G. Sinai, and K. M. Khanin, Feigenbaum universality and the thermodynamical formalism, Russian Math. Surveys **39** (1984), 1–40.